\documentstyle[11pt,twoside,amsfonts,amssymb]{article}
  
\begin{titlepage}

\title{The Spectrum of Kleinian Manifolds}

\author{Ulrich Bunke\thanks{Mathematisches Institut, Universit\"at G\"ottingen, Bunsenstr. 3-5, 37073 G\"ottingen, GERMANY, E-mail:bunke@uni-math.gwdg.de} and
Martin
Olbrich\thanks{Mathematisches Institut, Universit\"at G\"ottingen, Bunsenstr. 3-5, 37073 G\"ottingen, GERMANY, E-mail:bunke@uni-math.gwdg.de  }
}

\end{titlepage}

\newcommand{\proof}{{\it Proof.$\:\:\:\:$}}

\newcommand{\kaaa}{{\frak k}}

\newcommand{\taaa}{{\frak t}}
\newcommand{\haaa}{{\frak h}}

\newcommand{\R}{{\Bbb R}}
\newcommand{\Z}{{\Bbb Z}}
\newcommand{\C}{{\Bbb C}}
\newcommand{\HH}{{\Bbb H}}
\newcommand{\OO}{{\Bbb O}}

\newcommand{\gaaa}{{\frak g}}
\newcommand{\maaa}{{\frak m}}
\newcommand{\aaaa}{{\frak a}}
\newcommand{\naaa}{{\frak n}}

\newcommand{\res}{{\mathrm{ res}}}

\newcommand{\cZ}{{\cal Z}}
\newcommand{\cH}{{\cal H}}

\newcommand{\singsupp}{{\mathrm{ singsupp}}}

\newcommand{\cI}{{\cal I}}

\newcommand{\cC}{{\cal C}}
\newcommand{\mod}{{\mathrm{ mod}}}

\newcommand{\cO}{{\cal O}}

\newcommand{\cU}{{\cal U}}
\newcommand{\Hom}{{ \mathrm{ Hom}}}
\newcommand{\vol}{{\mathrm{ vol}}}

\newcommand{\vp}{\varphi}

\newcommand{\End}{{ \mathrm{ End}}}

\newcommand{\im}{{ \mathrm{ im}}}

\newcommand{\cF}{{\cal F}}
\newcommand{\Ree}{{\mathrm{ Re }}}

\newcommand{\inter}{{\mathrm{ int}}}
\newcommand{\clo}{{\mathrm{ clo}}}

\newcommand{\ee}{{\mathrm{ e}}}

\newcommand{\tr}{{ \mathrm{ tr}}}

\newcommand{\ad}{{ \mathrm{ ad}}}

\newcommand{\coker}{{\mathrm{ coker}}}
\newcommand{\id}{{ \mathrm{ id}}}

\newcommand{\nat}{{\Bbb  N}}
\newcommand{\supp}{{ \mathrm{ supp}}}

\newcommand{\aca}{{\aaaa_\C^\ast}}

\newcommand{\cR}{{\cal R}}
\def\hB{\hspace*{\fill}$\Box$\newline\noindent}

 \newcommand{\cG}{{\cal G}}

\newtheorem{prop}{Proposition}[section]
\newtheorem{lem}[prop]{Lemma}
\newtheorem{ddd}[prop]{Definition}
\newtheorem{theorem}[prop]{Theorem}
\newtheorem{kor}[prop]{Corollary}

\newtheorem{fact}[prop]{Fact}
\def\imath{i}
\newcommand{\ii}{i}

\begin{document}
\setcounter{page}{1}
\maketitle

\tableofcontents
 
\parskip3ex

\section{Introduction}

Let $G$ be a real simple linear connected Lie group of real rank one, and let $\Gamma\subset G$ be a convex-cocompact, non-cocompact, torsion-free, discrete
subgroup. This paper is devoted to the decomposition of the right regular
representation of $G$ on $L^2(\Gamma\backslash G)$ into irreducibles, the
so called Plancherel decomposition. We also allow twists by finite-dimensional unitary representations $(\varphi,V_\vp)$ of $\Gamma$, i.e., 
we investigate the right regular representation of $G$ on the space 
$$ L^2(\Gamma\backslash G,\vp):=\{f:G\rightarrow V_\vp\:|\: f(gx)=\vp(g)f(x)\ \forall g\in\Gamma,\: x\in G,\:\ 
\int_{\Gamma\backslash G} |f(x)|^2\:dx<\infty\}\ .$$  
Let $K\subset G$ be a maximal compact subgroup.
Then $X:=G/K$ is a Riemannian symmetric space of negative curvature which
is the universal covering of the locally symmetric space $Y:=\Gamma\backslash X$. Our assumption on $\Gamma$ implies that $Y$ has infinite volume and no
cusps. We call such a locally symmetric space a Kleinian manifold. Let $(\gamma,V_\gamma)$ be a finite-dimensional unitary representation
of $K$. Then we form the homogeneous vector bundle $V(\gamma):=G\times_KV_\gamma$ over $X$ and the locally
homogeneous vector bundle $V_Y(\gamma,\vp):=\Gamma\backslash (V(\gamma)\otimes V_\vp)$ over $Y$.
Let $\gaaa$ denote the Lie algebra of $G$, 
$\cU(\gaaa)$ the universal enveloping algebra of $\gaaa$ and $\cZ$ its center.
Via the left regular action of $\cU(\gaaa)$ on $C^\infty(X,V(\gamma))$ 
any  $A\in \cZ$ gives rise to a $G$-invariant differential
operator. This operator descends to $C^\infty(Y,V_Y(\gamma,\vp))$.
Here the Casimir operator $\Omega_G\in\cZ$, which defines an essentially
selfadjoint elliptic operator of second order acting on $L^2(Y,V_Y(\gamma,\vp))$, is of particular interest.

Our initial motivation was to obtain the spectral decomposition of
the space of sections $L^2(Y,V_Y(\gamma,\vp))$ of the bundle $V_Y(\gamma,\vp)$
over $Y$ with respect
to the Casimir operator and other locally invariant differential operators.
However, the isomorphism 
$$L^2(Y,V_Y(\gamma,\vp))\cong 
[L^2(\Gamma\backslash G,\vp)\otimes V_\gamma]^K$$ 
implies that the Plancherel
decomposition of $L^2(\Gamma\backslash G,\vp)$ is more or less equivalent to the desired spectral
decompositions for all the bundles  at once (for details see
Section \ref{didi}). 
Our main results are the Plancherel theorem Theorem \ref{pl2} and its consequences for spectral decompositions obtained in Theorem \ref{pl3}.

The structure of the Plancherel decomposition depends on the critical exponent $\delta_\Gamma$ of $\Gamma$ (see Definition \ref{dodo}). For technical reasons
we have to exclude  discrete subgroups  of
the isometry group of $X=\OO H^2$ with
$\delta_\Gamma\ge 0$ (because our method of the meromorphic continuation of Eisenstein series involves the "embedding trick" (see below) which does not work in this
case).
Then Theorem \ref{pl2} provides
a decomposition
$$ L^2(\Gamma\backslash G,\vp) = L^2(\Gamma\backslash G,\vp)_{ac}\oplus L^2(\Gamma\backslash G,\vp)_{cusp}\oplus L^2(\Gamma\backslash G,\vp)_{\res}\oplus L^2(\Gamma\backslash G,\vp)_{U}\ .
$$  
Here $L^2(\Gamma\backslash G,\vp)_{ac}$ decomposes further into a sum of direct integrals
corresponding to the unitary principal series representations of $G$, each
occuring with infinite multiplicity, $L^2(\Gamma\backslash G,\vp)_{cusp}$ decomposes into discrete series representations of $G$, each
discrete series representation of $G$ occurs with infinite multiplicity. 
These two parts, which are in a sense the main contribution to  $L^2(\Gamma\backslash G,\vp)$, have essentially the same structure as in the case
of the trivial group, i.e., in the Plancherel theorem for $L^2(G)$ due to
Okamoto \cite{okamoto65}, Hirai \cite{hirai66}, and Harish-Chandra \cite{harishchandra662}, \cite{HC76III}. The remaining two parts 
$L^2(\Gamma\backslash G,\vp)_{\res}$ and $L^2(\Gamma\backslash G,\vp)_{U}$
can only be non-trivial if
$\delta_\Gamma\ge 0$. They consist of a direct sum of non-discrete series representations of $G$ with real infinitesimal
character occuring with finite multiplicity. The subsbscript $\res$ stands
for residual spectrum. Indeed, the space $L^2(\Gamma\backslash G,\vp)_{\res}$
is generated by residues of Eisenstein series. The "stable" subspace $L^2(\Gamma\backslash G,\vp)_{U}$ is of similar nature but is orthogonal to the residues of Eisenstein
series. It contains representations of integral infinitesimal character only.
The understanding of its significance deserves further study.
 
Concerning the spectral decomposition of the operator $-\Omega_G$ acting
on $L^2(Y,V_Y(\gamma,\vp))$ it follows that the absolute continuous
spectrum consists of finitely many branches $[c_i,\infty]$ of infinite
multiplicity, where the constants $c_i\in\R$ are computable in terms of
$\gamma$, that there is no singular continuous spectrum, and that the discrete spectrum is finite. 
The set of eigenvalues with infinite multiplicity coincides with the discrete
spectrum of $-\Omega_G$ on $L^2(X,V(\gamma))$.
Eigenvalues
of finite multiplicity can only occur if $\delta_\Gamma\ge 0$. 
The corresponding eigenspaces split into a residual and a "stable" part, where
the residual part is generated by
residues of Eisenstein series. The "stable" part can occur only if $\gamma$ is non-trivial.
Note that the Plancherel theorem also provides a finer decomposition of the generalized eigenspaces
 of the Casimir
operator with respect to the algebra of locally invariant diffrential operators
$D(G,\gamma)$. For instance, the Casimir operator may have
embedded eigenvalues which are isolated with respect to $D(G,\gamma)$.
For more information see Section \ref{didi}.

Spectral decompositions of $L^2(Y,V_Y(\gamma,\vp))$ 
(respectively partial results; sometimes also cusps are allowed) were previously obtained for 
\begin{itemize} 
\item trivial $\gamma$ and surfaces by Patterson \cite{patterson75}, compare also \cite{fay77} and \cite{elstrodt75}
\item higher dimensional real hyperbolic manifolds and trivial $\gamma$  by Lax-Phillips \cite{laxphillips82}, 
\cite{laxphillips84}, \cite{laxphillips85}, Perry \cite{perry87}, Mazzeo-Melrose \cite{mazzeomelrose87}, Mandouvalos \cite{mandouvalos88} 
\item differential forms on real hyperbolic manifolds by Mazzeo-Phillips
\cite{mazzeophillips90}   
\item differential forms on complex-hyperbolic manifolds by Epstein-Melrose-Mendoza \cite{epsteinmelrosemendoza91}, Epstein-Melrose \cite{epsteinmelrose90}.
\end{itemize}
As in most of these papers the crucial step towards a spectral decomposition is the construction and
meromorphic continuation of Eisenstein series and the scattering matrix (at least
up to the critical axis). Besides the papers just cited this problem for trivial $\gamma$
and real hyperbolic manifolds $Y$ is
treated e.g. in
\cite{patterson76}, \cite{patterson761}, \cite{patterson89}, \cite{colindeverdiere85}, \cite{mandouvalos86}, \cite{agmon87},  \cite{mandouvalos89}, \cite{perry89}. 
Our approach differs philosophically,
if not technically, from these papers. After two sections of preparatory character Sections \ref{pdext}-\ref{uni} contain the development of our geometric version of
"scattering theory". The emphasis is on analysis on the sphere
at infinity, i.e., the geodesic boundary $\partial X$ of $X$. The advantage
of this approach becomes manifest if one goes beyond the case of the trivial
representation $\gamma$. Indeed, as experience shows, "meromorphic objects"
(Eisenstein series, scattering matrices, Selberg zeta functions etc.) correspond to families of bundles on the boundary rather than to bundles on the (locally) symmetric space itself. The key notions which we are going to discuss here are the extension map $ext$ and invariant distributions supported on the limit set. Very similar ideas  
appear in the work of van den Ban-Schlichtkrull and Delorme on the Plancherel
formula for reductive symmetric spaces (see e.g. \cite{vandenban88}, 
\cite{brylinskidelorme92}, \cite{carmonadelorme94}, \cite{vandenbanschlichtkrull97}).

Let $P\subset G$ be a minimal parabolic subgroup. Then $\partial X$ can be
viewed as the homogeneous space $G/P$. Set $M:=K\cap P$. Let $\sigma\in\hat M$ be an irreducible representation of $M$. Finite-dimensional irreducible 
representations of $P$ then come in families $\{\sigma_\lambda\}_{\lambda\in\aca}$ 
parametrized by a one-dimensional complex vector space $\aca$. They determine
families of homogeneous vector bundles $\{V(\sigma_\lambda)\}_{\lambda\in\aca}$
over $\partial X$.
Sections of these bundles carry representations of $G$, the so-called principal series representations. We are interested in the space of $\Gamma$-invariant distribution sections ${}^\Gamma C^{-\infty}(\partial X,V(\sigma_\lambda))$ as
well as in its twisted version ${}^\Gamma C^{-\infty}(\partial X,V(\sigma_\lambda,\vp))$, where $\vp$ is a finite-dimensional representation
of $\Gamma$. It turns out to be useful to allow non-unitary twists $\vp$, too. $\partial X$ is the union of the domain of discontinuity $\Omega$
of $\Gamma$ and the limit set $\Lambda$. The most interesting invariant
distributions are those which are supported on the limit set. By a slight abuse
of notation we denote the space of such distributions by ${}^\Gamma C^{-\infty}(\Lambda,V(\sigma_\lambda,\vp))$. A prominent example of an invariant distribution supported on the limit set is given by the Patterson-Sullivan measure which is
an element in ${}^\Gamma C^{-\infty}(\Lambda,V(1_{\delta_\Gamma}))$. Here $1$
stands for the trivial representation of $M$, and the critical exponent $\delta_\Gamma$ is viewed as an element of $\aca$ in a natural way.
Eventually, it will turn out that for unitary $\vp$ the representations appearing in $L^2(\Gamma\backslash G,\vp)_{\res}
\oplus L^2(\Gamma\backslash G,\vp)_{U}$ can be parametrized by the
set $\{(\sigma,\lambda)\in\hat M\times \aca\:|\:\Ree(\lambda)\ge 0, {}^{\Gamma} C^{-\infty}(\Lambda,V(\sigma_\lambda,\varphi))\not=0\}$.

In order to construct elements of ${}^\Gamma C^{-\infty}(\partial X,V(\sigma_\lambda,\vp))$ we proceed as follows. We consider the boundary
at infinity of $Y$, the compact
manifold $B:=\Gamma\backslash \Omega$. The bundle $V(\sigma_\lambda,\vp)$ on
$\partial X$ induces a corresponding bundle $V_B(\sigma_\lambda,\vp)$ on $B$
such that 
$$C^{-\infty}(B,V_B(\sigma_\lambda,\vp))\cong {}^\Gamma C^{-\infty}(\Omega,V(\sigma_\lambda,\vp))\ .$$
Using this isomorphism we want to construct the extension map
$$ext: C^{-\infty}(B,V_B(\sigma_\lambda,\vp))\rightarrow {}^\Gamma C^{-\infty}(\partial X,V(\sigma_\lambda,\vp))$$
which extends an invariant distribution on $\Omega$ across the limit set $\Lambda$. This is possible if $\lambda$ belongs to some right half-plane depending on the critical exponent $\delta_\Gamma$. In fact, we define
$ext$ to be the adjoint of an operator $\pi_*$ which associates to each smooth section of $V(\tilde\sigma_{-\lambda},\tilde\vp)$ the $\Gamma$-average of its restriction to
$\Omega$. 
Because of the convergence of the
Poincare series this average exists for $\Ree(\lambda)$ sufficiently large and
depends holomorphically on $\lambda$ in that region. The first task (Section \ref{extttttt}) is to obtain
a meromorphic continuation of $ext$ to all of $\aca$.

Classically, for trivial $\gamma$, Eisenstein series are obtained by averaging the Poisson kernel
$P_\lambda(x,b)$, $x\in X$, $b\in\Omega\subset\partial X$, over $\Gamma$.
Then the pairing of this Eisenstein series with a distribution $\phi$ on $B$ yields
the eigenfunction $E(\lambda,\phi)$ on $Y$, also called Eisenstein series.
Using the extension map we can rewrite 
\begin{equation}\label{narr}
E(\lambda,\phi)=P_\lambda\circ ext(\phi)\ ,
\end{equation}
where $P_\lambda$ is the Poisson transform. It is Equation (\ref{narr}) which we will use to define
Eisenstein series for general bundles. Thus, in our approach the extension map $ext$ is the primary object. Once $ext$ is understood, the Eisenstein series will not cause any essential additional difficulties. 

There is a second important object closely related to
$ext$, the scattering matrix $S_\lambda$ which we define as follows
\begin{equation}\label{scatext}
S_\lambda:= res\circ J_\lambda\circ ext: C^{-\infty}(B,V_B(\sigma_\lambda,\vp))\rightarrow C^{-\infty}(B,V_B(\sigma_{-\lambda},\vp))\ ,
\end{equation}
where
$$res:{}^\Gamma C^{-\infty}(\partial X,V(\sigma_\lambda,\vp))\rightarrow C^{-\infty}(B,V_B(\sigma_\lambda,\vp))$$
is induced by restriction from $\partial X$ to $\Omega$, and
$$ J_\lambda: C^{-\infty}(\partial X,V(\sigma_\lambda,\vp))\rightarrow C^{-\infty}(\partial X,V(\sigma_{-\lambda},\vp))$$
is the scattering matrix for the trivial group $\Gamma$ known in representation theory as Knapp-Stein intertwining operator . 

Initially, $ext$ and $S_\lambda$ are defined on some right half-plane. Their
meromorphic continuation proceeds in three surprisingly simple steps which we are going to sketch now. 

For step one 
assume that $ext$ is defined on a half-plane $\{\Ree(\lambda)>-\epsilon\}$ for some $\epsilon>0$. Then we show the functional equations 
\begin{eqnarray}
S_\lambda&=&S_{-\lambda}^{-1}\label{fun1}\\
 ext&=&J_{-\lambda}\circ ext\circ S_\lambda\  \label{fun2}
\end{eqnarray}
for $|\Ree(\lambda)|<\epsilon$. Under the additional hypotheses $\sigma=1$ we use meromorphic Fredholm theory in order to
show that $S_{-\lambda}^{-1}$ extends meromorphically to a much bigger half-plane. The main point here is that $J_\lambda$ can be used to construct a nice family of
parametrices for $S_{-\lambda}$. Now (\ref{fun1}) and (\ref{fun2}) give the
continuation of $S_\lambda$ and $ext$, respectively, to this half-plane. 

The remaining two steps are purely algebraic
in nature. In the second step we embed $G$, hence $\Gamma$, into the isometry group of a symmetric space of sufficiently large dimension. 
In this higher dimensional situation we can apply step one. By the first two steps we obtain the meromorphic continuation of $ext$ and $S_\lambda$ to a half-plane which is independent of $\delta_\Gamma$. 
This "embeddding trick" has appeared already in \cite{mandouvalos86} in the context
of the meromorphic continuation of Eisenstein series.
Unfortunately, it is not applicable to the exceptional case $X=\OO H^2$. 

In the third step we use tensoring with finite-dimensional
$G$-representations in order to embed the bundle $V_B(\sigma_\lambda,\vp)$
into a bundle of the form $V_B(1_\mu,\vp^\prime)$ for a suitable representation
$\vp^\prime$ of $\Gamma$ and $\mu\in\aca$ belonging to a region where $ext$
is already known to be meromorphic. In this step it is crucial to allow non-unitary
twists. Note that this method of meromorphic continuation is independent of any spectral
theoretic argument.

In Section \ref{invvv} we show how $ext$ can be used
to construct all invariant distributions $\phi\in {}^\Gamma C^{-\infty}(\partial X,V(\sigma_\lambda,\vp))$, in particular those supported on the
limit set. Indeed, it follows from $res\circ ext =\id$ that
at points $\lambda\in\aca$, where $ext$ has a pole, the leading singular part
of its Laurent expansion at $\lambda$ maps $C^{-\infty}(B,V_B(\sigma_\lambda,\vp))$ to ${}^\Gamma C^{-\infty}(\Lambda,V(\sigma_\lambda,\vp))$. This construction may be viewed
as a generalization of the construction of the Patterson-Sullivan measure.
The main result of this section is Theorem \ref{disfin} stating the discreteness
of the set of "resonances" $\{\lambda\in \aca\:|\:{}^{\Gamma} C^{-\infty}(\Lambda,V(\sigma_\lambda,\varphi))\not=0\}$ and that for each
$\lambda\in \aca$ the space ${}^{\Gamma} C^{-\infty}(\Lambda,V(\sigma_\lambda,\varphi))$ is finite-dimensional. In contrast to
the meromorphic continuation of $ext$ the proof of this theorem also
requires analysis on $X$, in particular a detailed knowledge of asymptotics of Poisson transforms and a variant of Green's formula. Note that the asymptotic formulas
we use are simple consequences of asymptotic expansions
of matrix coefficients known from representation theory. They are strong enough
to imply also the asymptotics of Eisenstein series which will play a decisive
role in the final proof of the Plancherel theorem. It also follows that the scattering matrix defined by (\ref{scatext}) is really a scattering matrix in the
sense that it
determines the relation between the two leading asymptotic terms of the Eisenstein series. In Section \ref{uni} we
return to the assumption that $\vp$ is unitary which allows us to gain more
detailed information concerning the location of the singularities of $ext$.
In particular, $ext$, and hence the Eisenstein series, are regular at non-zero imaginary $\lambda$. 

Note that only a part of the results obtained up to Section \ref{invvv} is
really needed for the derivation of the Plancherel theorem. However, we believe
that the present version of scattering theory also provides an adequate foundation
for more ambitious tasks in analysis on Kleinian manifolds as there are trace formulas,
the investigation of Selberg and Ruelle zeta functions or Paley-Wiener theorems. For instance, in \cite{bunkeolbrich963} we have used this approach in the case of real
hyperbolic manifolds in order to prove a conjecture of Patterson concerning
the relation between invariant distributions supported on the limit set and
the singularities of the (untwisted) Selberg zeta function. 
In a paper which is in preparation
we extend this approach to scattering theory to geometrically finite groups $\Gamma$.

The spectral theoretic part of the paper starts with Section \ref{abs}.
The problem we have to solve is twofold: first to produce a certain amount
of eigensections and wave packets of them, and second to show that they
span the whole Hilbert space. While the first task is almost standard once the Eisenstein series are constructed the second requires additional arguments.
In particular, one has to show the absence of the singular continuous spectrum.
Usually, the limiting absorption principle (e.g. \cite{perry87}) or commutator  methods (see e.g. \cite{froesehislopperry91}) are employed at this point.
Here we use a completely different method proposed by Bernstein \cite{bernstein88}. This method, brought to our attention by Delorme's proof of the Plancherel theorem for reductive symmetric spaces
\cite{delorme98}, rests on a theory of appropriate Schwartz spaces
for $Y$ (or $\Gamma\backslash G$). It leads to the notion of tempered
eigensections or, switching to representation theoretic language, tempered invariant distribution vectors of unitary representations of $G$. 
These notions will be discussed in Section \ref{abs}. The crucial point is that a priori only tempered eigensections can occur in the spectral decomposition. In Section
\ref{relsec} we relate tempered invariant distribution vectors to invariant distributions on $\partial X$, 
in particular those supported on the limit set. Combining the results of
Sections \ref{invvv} and \ref{uni} with knowledge of the structure
of the unitary dual $\hat G$ of $G$ we obtain a classification of tempered invariant distribution vectors. This classification enables us to complete
the exhaustion part of the proof of the Plancherel theorem in Section \ref{didi}.

Some of the readers may have noticed an earlier version of this paper which
appeared as an e-print more than two years ago. At that time we were not able
to continue $ext$ to all of $\aca$, but only to the complement of a set of
integer points. In particular, in order to conclude finiteness of the discrete
spectrum we were forced to combine several not very natural arguments. Now the
discreteness result Theorem \ref{disfin} which is a consequence of the meromorphy
of $ext$ gives (among other things) a very satisfactory understanding of the
finiteness of the discrete spectrum. In addition, the point of view in the present
version is more representation theoretic. In our opinion, this makes the fine
structure of the spectrum much more transparent.

\noindent
{\it Acknowledgement: We thank R. Mazzeo and P. Perry for discussing
of parts of this work.}

\section{Geometric preparations}\label{fiert}

Let $G$ be a connected, linear, real simple Lie group of rank one, $G=KAN$ be an Iwasawa decomposition
of $G$, $\gaaa=\kaaa\oplus\aaaa\oplus\naaa$ be the corresponding Iwasawa decomposition of the Lie algebra $\gaaa$,
$M:=Z_K(A)$ be the centralizer of $A$ in $K$, and $P:=MAN$
be a minimal parabolic subgroup. 
The group $G$ acts isometrically on the rank-one symmetric space $X:=G/K$.
Let $\partial X:=G/P=K/M$ be its geodesic boundary. We consider $\bar{X}:=X\cup\partial X$ as a compact manifold with boundary.

By the classification of symmetric spaces with strictly negative
sectional curvature $X$ is one of the following spaces:
\begin{itemize}
\item a real hyperbolic space $\R H^n,\ n\ge 1$,
\item a complex hyperbolic space $\C H^n,\ n\ge 2$,
\item a quaternionic hyperbolic space $\HH H^n,\ n\ge 2$,
\item or the Cayley hyperbolic plane $\OO H^2$, 
\end{itemize}
and $G$ is a linear group finitely covering the orientation-preserving isometry 
group of $X$. 

We consider a torsion-free discrete subgroup $\Gamma \subset G$ such that $\partial X$
admits 
a $\Gamma$-invariant partition $\partial X =\Omega\cup \Lambda$, where  $\Omega\not=\emptyset$ is open and 
$\Gamma$ acts freely and cocompactly on 
$X\cup\Omega$. The closed subset $\Lambda$ is called the limit set of $\Gamma$.
The locally symmetric space $Y:=\Gamma\backslash X$ is a complete Riemannian manifold of infinite volume without cusps. It can be compactified by adjoining
the geodesic boundary $B:=\Gamma\backslash \Omega$.

A subgroup $\Gamma$ satisfying this assumption is often called convex-cocompact
since it acts cocompactly on the convex hull of the limit set.
The quotient $Y$ can be called a Kleinian manifold generalizing
the corresponding notion for three-dimensional hyperbolic manifolds. 

By $\aca$ we denote the comlexified dual of $\aaaa$.
If $\lambda\in \aca$, then we set $a^\lambda:=\ee^{\langle \lambda,\log(a)\rangle}\in\C$. Let $\alpha$ be the short root of
$\aaaa$ in $\naaa$. We set 
$A_+:=\{a\in A\:|\:a^\alpha\ge 1\}$.
Define $\rho\in \aaaa^*$ as usual by $\rho(H):=\frac{1}{2}\tr(\ad(H)_{|\naaa})$, $\forall H\in\aaaa$. We have\\  
\centerline{\begin{tabular}{|c||c|c|c|c|}
\hline
$X$&$\R H^n$&$\C H^n$&$\HH H^n$&$\OO H^n$\\
\hline
$\rho$&$\frac{n-1}{2}\alpha$&$n\alpha$&$(2n+1)\alpha$&$11\alpha$\\
\hline
\end{tabular}\ . }

We adopt the following conventions about the notation for points of $X$ and $\partial X$.
A point $x\in \partial X$ can equivalently be denoted by a subset $kM\subset K$
or $gP\subset G$ representing this point in $\partial X=K/M$ or $\partial X=G/P$.
If $F\subset \partial X$, then $FM:=\bigcup_{kM\in F}kM\subset K$.
Analogously, we can denote $b\in X$ by $gK\subset G$, where $gK$ represents
$b$ in $X=G/K$.

\begin{lem}\label{no}
For any compact $F\subset \Omega$ we have $\sharp(\Gamma\cap (FM)A_+K)<\infty$.
\end{lem}
\proof
Note that $(FM)A_+K\cup F\subset X\cup\Omega$ is compact. Thus its intersection with the orbit $\Gamma K$ of the origin of $X$ is finite. \hB

The geometry of the action of $\Gamma$ on $X\cup\partial X$ can
be studied in a uniform way using the various decompositions of $G$. Any element $g\in G$ has a Cartan decomposition $g=k_ga_gh$, $k_g,h\in K$, $a_g\in A_+$, where $a_g$ and $k_gM\in K/M$ are uniquely determined by $g$.
Let $g=\kappa(g)a(g)n(g)$, $\kappa(g)\in K$, $a(g)\in A$, $n(g)\in N$ be
defined with respect to the given Iwasawa decomposition. 
The function
$G\times K\ni (g,k)\mapsto a(g^{-1}k)$ descends to $X\times\partial X$.

Given a normalization of the invariant distance $d$ on $X$ we identify
$A_+$ with $[1,\infty)$ such that
$a=\ee^{d(eK,aK)}$. Then for any $g\in G$ and $k\in K$ we have
\begin{equation}\label{mufti}
a_g=\ee^{d(eK,gK)}\ \  \mbox{ and }\ \  a(g^{-1}k)=\ee^{\pm d(eK,HS_{gK,kM})}\ ,
\end{equation}
where $HS_{gK,kM}=kNk^{-1}gK$ is the horosphere passing through $gK\in X$ and $kM\in\partial X$. The sign $\pm$ is positive (negative) if $eK$ lies inside
(outside) the corresponding horoball. 

\begin{ddd}\label{dodo}
The critical exponent $\delta_\Gamma\in \aaaa^*$ of $\Gamma$ is the smallest element such that $\sum_{g\in\Gamma} a_g^{-(\lambda+\rho)}$ converges for all $\lambda\in\aaaa^*$ with $\lambda>\delta_\Gamma$. If $\Gamma$ is the trivial group,
then we set $\delta_\Gamma:=-\infty$.
\end{ddd}
Equation (\ref{mufti}) shows that this definition of $\delta_\Gamma$  differs from the
usual one by a $\rho$-shift, only.

The critical exponent $\delta_\Gamma$ has been extensively studied, in particular by Patterson \cite{patterson762}, Sullivan \cite{sullivan79}, and Corlette \cite{corlette90}. From these papers we know that $\delta_\Gamma\in [-\rho,\rho)$, if $\Gamma$ is non-trivial.
Moreover, $\delta_\Gamma+\rho=\dim_H(\Lambda)\alpha$, where $\dim_H(\Lambda)$ denotes the Hausdorff dimension of the limit set with respect to the natural class of sub-Riemannian metrics on $\partial X$
(the Hausdorff dimension of the empty set is
by definition $-\infty$). If $X$ is a quaternionic hyperbolic space or the
Cayley hyperbolic plane, then $\delta_\Gamma$ can not be arbitrary close to
$\rho$. In these cases we have $\delta_\Gamma\le (2n-1)\alpha$ and $\delta_\Gamma\le 5\alpha$, respectively \cite{corlette90}.

We now collect some facts concerning the relation between the Cartan and
the Iwasawa decomposition. First of all Equation (\ref{mufti}) implies
\begin{equation}\label{grufti}
a(g^{-1}k) \le  a_g \quad \mbox{ for all } g\in G,\ k\in K\ .
\end{equation}
 
\begin{lem}\label{no1}
Let $k_0M \in \partial X$. For any compact $W\subset (\partial X\setminus k_0M)M$ and any neighbourhood $U\subset K$ of $k_0M$ satisfying $W\subset (\partial X\setminus \clo(U)M)M$ there exists a constant $C>0$, such that   
\begin{displaymath}
C a_g \le  a(g^{-1}k) \le  a_g 
\end{displaymath}
for all $g=k_ga_gh\in KA_+K$ with $k_g \in W$, and all $k\in U$.
\end{lem}
\proof 
The upper bound is given by (\ref{grufti}). We prove the lower bound.
Let $w\in N_K(M)$ represent the non-trivial element of the Weyl group of $(\gaaa,\aaaa)$.
Set   $\bar{N}=\theta(N)$,
where $\theta$ is the Cartan involution of $G$ fixing $K$.
Since the set $W^{-1}\clo(U)M$ is compact and disjoint from $M$, there is a precompact open $V\subset \bar{N}$ such that $W^{-1}\clo(U)M\subset w \kappa(V)M$.

Let $k\in U$.
Then we have $k_g^{-1}k=w\kappa(\bar{n})m $ for $\bar{n}\in V$, $m\in M$.
Using that $a(\bar{n})\ge 1$ for all $\bar{n}\in \bar{N}$ (see e.g. \cite{helgason84}, Ch. IV, Cor. 6.6.)  we obtain
\begin{eqnarray*}
a(g^{-1}k)&=&a(a_g^{-1}k_g^{-1}k)\\
&=& a(a_g^{-1}w\kappa(\bar{n})m) \\
&=&a(a_g \kappa(\bar{n}))\\
&=&a(a_g\bar{n}n(\bar{n})^{-1}a(\bar{n})^{-1})\\
&=&a(a_g\bar{n}a_g^{-1}) a(\bar{n})^{-1} a_g\\
&\ge&a(\bar{n})^{-1} a_g\ .
\end{eqnarray*}
Since $V$ is precompact we have $C:=\inf_{\bar{n}\in V}a(\bar{n})^{-1}>0$.
It follows that $a(g^{-1}k )\ge Ca_g\ .$
\hB

As a corollary we obtain a certain converse of the triangle inequality.

\begin{kor}\label{plumps}
Let $k_0$, $W$, and $U$ be as in Lemma \ref{no1}. Then there exists a constant $C>0$ such that for all $a\in A$, $k\in U$, and $g\in G$ with $k_g\in W$
$$ a_{g^{-1}ka}\ge Ca_ga\ .$$
\end{kor}  
\proof
Combining Equation (\ref{grufti}) with Lemma \ref{no1} we obtain
$$ a_{g^{-1}ka}=a_{(g^{-1}ka)^{-1}}\ge a(g^{-1}ka)=a(g^{-1}k)a\ge Ca_ga\ .$$
\hB

A word concerning normalizations: The basic object will be a fixed invariant Riemannian
metric on $X$. Here the reader has the freedom to choose his favoured
one. The exponential map then induces a metric on $\aaaa$, hence on $\aaaa^*$. Throughout the paper we will often isometrically identify $\aaaa$ with $\R$,
$\aca$ with $\C$, $A_+$ with $[1,\infty)$ such that (\ref{mufti}) holds.
Haar measures will be normalized as follows: The measures on $A$, $\aaaa^*$ are fixed by
the above metric. Compact groups will always have total
mass $1$. The Haar measure of $G$ is given by $dg=dk\:dx$, where $dk$ is the
Haar measure of $K$ and $dx$ is the Riemannian measure on $X=G/K$. Finally,
we will 
normalize the Haar measure on $\bar N$ such that
$$\int_{\bar N} a(\bar n)^{-2\rho} \:d\bar n =1\ .$$

\section{Analytic preparations}\label{anaprep}

Let $(\sigma,V_\sigma)$ be a finite-dimensional unitary representation of $M$. For $\lambda\in \aca$  we  form the representation $\sigma_\lambda$
of $P$ on $V_{\sigma_\lambda}:=V_\sigma$, which is given by
$\sigma_\lambda(man):=\sigma(m)a^{\rho-\lambda}$.
Let $V(\sigma_\lambda):=G\times_P V_{\sigma_\lambda}$ be the associated
homogeneous bundle over $\partial X=G/P$. 
It induces a bundle on $B=\Gamma\backslash\Omega$ defined by $V_B(\sigma_\lambda):=\Gamma\backslash V(\sigma_\lambda)_{|\Omega}$.

Let $\tilde{\sigma}$ be the dual representation to $\sigma$.
Then there are  natural pairings 
\begin{eqnarray*}
V(\tilde{\sigma}_{-\lambda})\otimes V(\sigma_\lambda)&\rightarrow& \Lambda^{max}T^*\partial X\\
V_B(\tilde{\sigma}_{-\lambda})\otimes V_B(\sigma_\lambda)&\rightarrow& \Lambda^{max}T^*B\ .
\end{eqnarray*}
The orientation of $\partial X$ induces one of $B$.
Employing these pairings and integration with respect to the fixed
orientation we obtain identifications
\begin{eqnarray*}
C^{-\infty}(\partial X,V(\sigma_\lambda))&=&C^{\infty}(\partial X,V(\tilde{\sigma}_{-\lambda}))^\prime\\
C^{-\infty}(B,V_B(\sigma_\lambda))&=&C^{\infty}(B,V_B(\tilde{\sigma}_{-\lambda}))^\prime\ .
\end{eqnarray*}

As a $K$-homogeneous bundle we have a canonical identification $V(\sigma_\lambda)\cong K\times_M V_\sigma$. Thus $\bigcup_{\lambda\in\aca} V(\sigma_\lambda)\rightarrow \aca\times \partial X$
has the structure of a trivial holomorphic family of bundles.

Let $\pi^{\sigma,\lambda}$ denote the representation of $G$ on the space of sections of $V(\sigma_\lambda)$ given by the left-regular representation. 
Then $\pi^{\sigma,\lambda}$ is called a principal series representation
of $G$. Note that there are different globalizations of this
representation which are distinguished by the regularity
of the sections (smooth, distribution  etc.).

For any small open subset $U\subset B$
and diffeomorphic lift $\tilde{U}\subset \Omega$  the restriction $V_B(\sigma_\lambda)_{|U}$
is canonically isomorphic to $V(\sigma_\lambda)_{|\tilde{U}}$.
Let $\{U_\alpha\}$ be a cover of $B$ by open sets as above.
Then  
$$\bigcup_{\lambda\in\aca} V_B(\sigma_\lambda)\rightarrow \aca\times B$$
can be given the structure of a holomorphic family of bundles by gluing the trivial families
$$\bigcup_{\lambda\in\aca}V_B(\sigma_\lambda)_{|U}\cong \bigcup_{\lambda\in\aca}V(\sigma_\lambda)_{|\tilde{U}}$$ together
using the holomorphic families of gluing maps induced by $\pi^{\sigma,\lambda}(g)$, $g\in\Gamma$.

The structure of a holomorphic family of bundles allows us to consider holomorphic or smooth or continuous families
of sections $\aca\ni\mu\mapsto f_\mu\in C^{\pm\infty}(\partial X,V(\sigma_\mu))$,
$\aca\ni\mu\mapsto f_\mu\in C^{\pm\infty}(B,V_B(\sigma_\mu))$, respectively.

Let $(\vp,V_\vp)$ be a finite-dimensional representation of $\Gamma$.
We form the bundle $V(\sigma_\lambda,\varphi):=V(\sigma_\lambda)\otimes V_\vp$ on $\partial X$  carrying the tensor product action of $\Gamma$ and define $V_B(\sigma_\lambda,\varphi):=\Gamma\backslash (V(\sigma_\lambda)\otimes V_\vp)_{|\Omega}$. Then we have the spaces of sections $C^{\pm\infty}(\partial X,V(\sigma_\lambda,\varphi))$ and $C^{\pm\infty}(B,V_B(\sigma_\lambda,\varphi))$ as
well as the various notions of $\aca$-parametrized families of sections. 

When dealing with holomorphic families of vectors in topological vector
spaces we will employ the following functional analytic facts. Let $\cF,\cG,\cH \dots$ be complete locally convex
topological vector spaces.
A locally convex vector space is called a Montel space if its   
closed bounded subsets are compact. 
A Montel space is reflexive, i.e., the canonical map into its bidual is an isomorphism.
Moreover, the dual space of a Montel space is again a Montel space.
\begin{fact} The space of smooth sections of a vector bundle
and its topological dual are Montel spaces.
\end{fact}
We equip $\Hom(\cF,\cG)$ with the topology of uniform convergence
on bounded sets.
  Let $V\subset \C$ be open.
A map $f:V\rightarrow \Hom(\cF,\cG)$ is called holomorphic
if for any $z_0\in V$ there is a sequence $f_i\in \Hom(\cF,\cG)$ such that
 $f(z)=\sum_{n=0}^\infty f_i (z-z_0)^i$ converges for all $z$ close to $z_0$.
Let $f:V\setminus \{z_0\} \rightarrow \Hom(\cF,\cG)$ be holomorphic
and  $f(z)=\sum_{n=-N}^\infty f_i (z-z_0)^i$ for all  $z\not=z_0$ close to $z_0$.
Then we say that  $f$ is meromorphic and has a pole of order $N$ at $z_0$.
If $f_i$, $i=-N,\dots,-1$, are finite dimensional, then $f$
has, by definition, a finite-dimensional singularity.
We call a subset  $A\subset \cF\times \cG^\prime$ sufficiently large if for    
$B\in \Hom(\cF,\cG)$ the condition 
 $\langle \phi,B \psi \rangle=0$, $ \forall (\psi,\phi)\in A$,  implies  $B=0$.
Proofs of
the following facts can be found in \cite{bunkeolbrich963}, Section 2.2.

\begin{fact}\label{holla}
The following assertions are equivalent :
\begin{enumerate}
\item (i)  $\:\:f:V\rightarrow \Hom(\cF,\cG)$ is holomorphic.
\item (ii) $\:\:\: f$ is continuous and 
there is a sufficiently large subset $A\subset \cF\times \cG^\prime$
such that for all $(\psi,\phi)\in A$ the function $V\ni z\mapsto \langle \phi,f(z)\psi\rangle $
is holomorphic.
\end{enumerate}
\end{fact}
 
\begin{fact}\label{seq}
Let $f_i:V\rightarrow Hom(\cF,\cG)$ be a sequence
of holomorphic maps. Moreover let $f :V\rightarrow Hom(\cF,\cG)$ be continuous such that for a sufficiently large subset
$A \subset \cF\times \cG^\prime$ the functions $\langle \phi,f_i \psi\rangle $,
$(\psi,\phi)\in A$,    
converge locally uniformly
in $V$ to $\langle\phi,f \psi\rangle$. Then $f$ is holomorphic, too.
\end{fact}
 
\begin{fact}\label{adjk}
Let $f:V\rightarrow \Hom(\cF,\cG)$
 be continuous.
Then the adjoint $f^\prime:V\rightarrow \Hom(\cG^\prime,\cF^\prime)$
is continuous.
If $f$ is holomorphic, then so is $f^\prime$. 
\end{fact}

\begin{fact}\label{comp}
Assume that $\cF$ is a Montel space.
Let $f:V\rightarrow \Hom(\cF,\cG)$
and $f_1:V\rightarrow \Hom(\cG,\cH)$
be continuous.
Then $f_1\circ f : V\rightarrow \Hom(\cF,\cH)$
is continuous.
If $f,f_1$ are holomorphic, so is $f_1\circ f$. 
\end{fact}

Let $\cH$ be a Hilbert space and $\cF\subset \cH$
be a Fr\'echet space such that the embedding is continuous and compact.
In the application we have in mind $\cH$ will be some $L^2$-
space of sections of a vector bundle over a compact closed manifold
and $\cF$ be the Fr\'echet space of smooth sections of this bundle.
The continuous maps $\Hom(\cH,\cF)$ will  be called smoothing operators.

Let $V\subset \C$ be open and connected, and $V\ni z\rightarrow R(z)\in \Hom(\cH,\cF)$
be a meromorphic family of smoothing operators  with at most finite-dimensional singularities.
Note that $R(z)$ is a meromorphic family of compact operators on $\cH$ in a natural way.

\begin{lem} \label{merofred}
If $1-R(z)$ is invertible at some point $z\in V$ where $R(z)$ is regular,
then 
$$(1-R(z))^{-1}=1-S(z)\ ,$$ 
where $V\ni z\rightarrow S(z)\in \Hom(\cH,\cF)$
is a meromorphic family of smoothing operators  with at most finite-dimensional singularities.
\end{lem}
\proof
We apply Reed-Simon IV \cite{reedsimon78}, Theorem XIII.13
in order to conclude that $(1-R(z))^{-1}$ is a meromorphic
family of operators on $\cH$ having at most finite-dimensional singularities.
Making the ansatz $(1-R(z))^{-1}= 1-S(z)$, where apriori $S(z)$ is a meromorphic
familiy of bounded operators on $\cH$ with finite-dimensional singularities,
we obtain $S=-R-R\circ S$. This shows that $S$ is a meromorpic family in $\Hom(\cH,\cF)$.
\hB

\section{Push-down and extension}\label{pdext}

Distribution sections of $V_B(\sigma_\lambda,\vp)$ can be identified with
$\Gamma$-invariant sections of $V(\sigma_\lambda,\vp)_{|\Omega}$. 
In order to extend these distributions across the limit set in an invariant way we will construct the extension map 
$$ext: C^{-\infty}(B,V_B(\sigma_\lambda,\varphi))\rightarrow {}^\Gamma C^{-\infty}(\partial X,V(\sigma_\lambda,\varphi))\ .$$
We first introduce its "adjoint"
which is the push-down
$$\pi_\ast:C^\infty(\partial X,V(\sigma_\lambda,\varphi))\rightarrow  C^\infty(B,V_B(\sigma_\lambda,\varphi))\ .$$ 
 
Using the identification  $C^\infty(B,V_B(\sigma_\lambda,\varphi))= {}^\Gamma C^\infty(\Omega,V(\sigma_\lambda,\varphi))$
we define $\pi_\ast$ by
\begin{equation}\label{summ}
\pi_*(f)(kM)= \sum_{g\in\Gamma} (\pi(g)f)(kM),\quad kM\in\Omega \ ,
\end{equation}
if the sum converges. Here $\pi(g)$ is the action
induced by $\pi^{\sigma,\lambda}(g)\otimes\varphi(g)$.

Note that the universal enveloping algebra 
$\cU(\gaaa)$ is a filtered algebra. Let $\cU(\gaaa)_m$, $m\in\nat_0$,
be the space of elements of degree less or equal than $m$.
Choose some norm on $V_\vp$. It induces a norm $|.|$ on $V_{\sigma}\otimes V_\vp$. For any $m$ and bounded subset $A\subset \cU(\gaaa)_m$
we define the seminorm $\rho_{m,A}$ on $C^\infty(\partial X,V(\sigma_\lambda,\varphi))$
by $$\rho_{m,A}(f):=\sup_{X\in A, k\in K}|f(\kappa(kX))|\ ,$$
where we consider $f$ as a function on $K$ with values in $V_{\sigma}\otimes V_\vp$.
These seminorms define the Fr\'echet topology of $C^\infty(\partial X,V(\sigma_\lambda,\varphi))$
(in fact a countable set of such seminorms is sufficient).

In order to describe the Fr\'echet topology on $C^\infty(B,V_B(\sigma_\lambda,\varphi))$
we fix an open cover $\{U_\alpha\}$ of $B$ such that each $U_\alpha$ has a diffeomorphic
lift $\tilde{U}_\alpha \subset\Omega$.
Then we have canonical isomorphisms 
$$C^\infty(\tilde{U}_\alpha ,V(\sigma_\lambda,\varphi))\cong C^\infty(\ U_\alpha,V_B(\sigma_\lambda,\varphi))\ .$$
For any $U\in \{U_\alpha\}$ we define the topology of $C^\infty(\tilde{U},V(\sigma_\lambda,\varphi))$
using the seminorms 
$$\rho_{U,m,A}(f):=\sup_{X\in A, k\in \tilde{U}M }|f(\kappa(kX))|\ ,$$
where  $m\in\nat_0$ and $A\subset \cU(\gaaa)_m$ is bounded.
Since $C^\infty(B,V_B(\sigma_\lambda,\varphi))$ maps to $C^\infty( U_\alpha,V_B(\sigma_\lambda,\varphi))$
by restriction for each $\alpha$ we obtain a system of seminorms defining the Fr\'echet topology of $C^\infty(B,V_B(\sigma_\lambda,\varphi))$.

Since $\Gamma$ is finitely generated we can find an element $\mu\in\aaaa^*$, $\mu \ge 0$, and a constant $C$ such that 
\begin{equation}\label{wau}
||\varphi(g)||\le Ca_g^{\mu} \quad \mbox{ for all } g\in \Gamma\ .
\end{equation}
 
\begin{ddd}
Let $\delta_\varphi\in \aaaa^*$ be the infimum of all $\mu\in\aaaa^*$ satisfying Equation (\ref{wau}) for some $C$. It is independent of the chosen norm. We call $\delta_\varphi$ the exponent of $(\varphi,V_\vp)$.
\end{ddd}
Note that unitary representations have zero exponents. If $\tilde\varphi$ is the dual representation of $\varphi$, then $\delta_{\tilde\varphi}=\delta_\varphi$. Furthermore, if
$\varphi$ is the restriction of a finite-dimensional representation of $G$ with highest $\aaaa$-weight $\mu$, then $\delta_\varphi=\mu$. 

\begin{lem}\label{anal1}
If $\Ree(\lambda)<-(\delta_\Gamma+\delta_\varphi)$, then the sum (\ref{summ}) converges for $f\in C^\infty(\partial X,V(\sigma_\lambda,\varphi))$ and defines
a  continuous map $$\pi_\ast: C^\infty(\partial X,V(\sigma_\lambda,\varphi))\rightarrow  C^\infty(B,V_B(\sigma_\lambda,\varphi))\ .$$
Moreover, $\pi_*$ depends holomorphically on $\lambda$.
\end{lem}
\proof 
Consider $U\in \{U_\alpha\}$. We want to estimate  
$$C^\infty(\partial X,V(\sigma_\lambda,\varphi))\ni f\mapsto res_{\tilde{U}}\circ \pi(g) f \in 
C^\infty(\tilde{U},V(\sigma_\lambda,\varphi))\ .$$
Let $\Delta:\cU(\gaaa)\rightarrow \cU(\gaaa)\otimes \cU(\gaaa)$
be the coproduct and write $\Delta(X)=\sum_i X_i\otimes Y_i$.
Fix $l\in \nat_0$ and a bounded set $A\in \cU(\gaaa)_l$.
Then for any $\epsilon>0$ there is another bounded set $A_1\subset \cU(\gaaa)_l$
depending on $A$ such that
\begin{eqnarray}\label{kati}
\rho_{U,l,A}(res_{\tilde{U}M}\circ \pi(g) f)&=&\sup_{X\in A,k\in \tilde{U}}|(\pi(g)f)(\kappa(kX))|\nonumber\\
&=& \sup_{X\in A,k\in \tilde{U}M}|\sum_ i a(g^{-1}\kappa(kX_i))^{\lambda-\rho}(\id\otimes\varphi(g))
f(\kappa(g^{-1}\kappa(kY_i)))|\nonumber\\
&\le & 
a_g^{\delta_\varphi+\epsilon}\sup_{X\in A_1,k\in \tilde{U}M}| a(g^{-1}kX)^{\lambda-\rho}| \sup_{X\in A_1,k\in \tilde{U}M}|   f(\kappa(g^{-1}kX))| \ .
\end{eqnarray}

The Poincar\'e-Birkhoff-Witt theorem gives a  decomposition $\cU(\gaaa)=\cU(\bar{\naaa})\cU(\maaa)\cU(\aaaa)\oplus \cU(\gaaa)\naaa$.
Let $q:\cU(\gaaa)\rightarrow \cU(\bar{\naaa})\cU(\maaa)\cU(\aaaa)$
be the associated projection. Then for $g\in G$ and $X\in\cU(\gaaa)$ we have 
$\kappa(gX)=\kappa(gq(X))$, $a(gX)=a(gq(X))$.

Let $U_1\subset \Omega$ be an open neighbourhood of $\tilde{U}$.
Then by Lemma \ref{no} the intersection $\Gamma\cap U_1MA_+K$ is finite.
Let $W:=(\partial X\setminus U_1)M$. As in the proof of Lemma \ref{no1}
we can find a compact $A_+$-invariant set $V\subset \bar{\naaa}$
such that $W^{-1}\tilde{U}M\subset w\kappa(V)M$.
For $g=ha_gh^\prime\in WA_+K$ and $k\in \tilde{U}M$ we obtain $h^{-1}k=w\kappa(\bar{n})m$ for some $\bar{n}\in V$, $m\in M$.

Let $X\in \cU(\gaaa)$. Then
\begin{eqnarray*}
\kappa(g^{-1}kX)&=&\kappa(h^{\prime -1}a_g^{-1}h^{-1}kX)\\
&=&h^{\prime -1}\kappa(a_g^{-1}w\kappa(\bar{n})mX)\\
&=&h^{\prime -1}w\kappa(a_g\bar{n}n(\bar{n})^{-1}a(\bar{n})^{-1}mX) \\
&=&h^{\prime -1}w\kappa(a_g\bar{n}a_g^{-1} a_g[n(\bar{n})^{-1}a(\bar{n})^{-1}mXm^{-1}a(\bar{n})n(\bar{n})]a_g^{-1})m\\
&=&h^{\prime -1}w\kappa(a_g\bar{n}a_g^{-1} a_gq(n(\bar{n})^{-1}a(\bar{n})^{-1}mXm^{-1}a(\bar{n})n(\bar{n}))a_g^{-1})m\ .
\end{eqnarray*}
Since $V$ is compact the sets  $n(V)^{-1}a(V)^{-1}MA_1Ma(V)n(V)=:A_2\subset \cU(\gaaa)_l$
and  $q(A_2)$ are bounded. 
Conjugating $q(A_2)$ with $A_+$ gives
clearly another bounded set $A_3\subset \cU(\gaaa)_l$.
We can find a bounded set $A_4\subset \cU(\gaaa)_l$ such that
$\kappa(a_g\bar{n}a_g^{-1}A_3)\subset \kappa(\kappa(a_g\bar{n}a_g^{-1})A_4)$
for all $a_g\in A_+$.
This implies for $g\in WA_+K$ that
\begin{equation}\label{kkll1}\sup_{X\in A_1,k\in \tilde{U}M}|   f(\kappa(g^{-1}kX))|\le \rho_{l,A_4}(f)\ .\end{equation}
We also have 
\begin{eqnarray*}
a(g^{-1}kX)&=&a(h^{\prime -1}a_g^{-1}h^{-1}kX)\\
&=& a(a_g^{-1}w\kappa(\bar{n})mX)\\
&=&a(a_g \kappa(\bar{n})mXm^{-1})\\
&=&a(a_g\bar{n}n(\bar{n})^{-1}a(\bar{n})^{-1}mXm^{-1})\\
&=&a(a_g\bar{n}a_g^{-1}a_g n(\bar{n})^{-1}a(\bar{n})^{-1}mXm^{-1}a(\bar{n})n(\bar{n})a_g^{-1}) a(\bar{n})^{-1} a_g\ .\\
&=&a(a_g\bar{n}a_g^{-1}a_g q(n(\bar{n})^{-1}a(\bar{n})^{-1}mXm^{-1}a(\bar{n})n(\bar{n}))a_g^{-1}) a(\bar{n})^{-1} a_g\ . 
\end{eqnarray*}
As in the proof of Lemma \ref{no1} there is a constant $C<\infty$ such that 
$$ |a(a_g\bar{n}a_g^{-1}a_g q(n(\bar{n})^{-1}a(\bar{n})^{-1}mXm^{-1}a(\bar{n})n(\bar{n}))a_g^{-1})^{\lambda-\rho}| a(\bar{n})^{ \rho- \lambda } <C$$
for all $a_g\in A_+$, $\bar{n}\in V$, $m\in M$, and $X\in A_1$.
It follows that
\begin{equation}\label{hj1}\sup_{X\in A_1,k\in \tilde{U}M}| a(g^{-1}kX)^{\lambda-\rho}|\le C a_g^{\lambda-\rho}\end{equation}
for almost all $g\in \Gamma$.
The estimates (\ref{kati}), (\ref{kkll1}) and (\ref{hj1}) together imply that the sum
$$C^l(\partial X,V(\sigma_\lambda,\varphi))\ni f\mapsto \sum_{g\in\Gamma} res_{\tilde{U}}\circ \pi(g) f \in  C^l(\tilde{U},V(\sigma_\lambda,\varphi))$$
converges for $\Ree(\lambda)<-(\delta_\Gamma+\delta_\varphi)$ 
and defines a continuous map of Banach spaces.
This map depends holomorphically on $\lambda$
by Fact \ref{seq}.

Combining  these considerations for all  
$U\in \{U_\alpha\}$ and $l\in\nat_0$ we obtain that 
$$\pi_*:C^\infty(\partial X,V(\sigma_\lambda,\varphi))\rightarrow C^\infty(B,V_B(\sigma_\lambda,\varphi))$$
is defined and continuous for  $\Ree(\lambda)<-(\delta_\Gamma+\delta_\varphi)$. 
Moreover it is easy to see that $\pi_*$ depends holomorphically
on $\lambda$. \hB
\begin{ddd}\label{defofext}
For $\Ree(\lambda)>\delta_\Gamma+\delta_\varphi$ we define 
the extension map 
$$ext: C^{-\infty}(B,V_B(\sigma_\lambda,\varphi))\rightarrow {}^\Gamma C^{-\infty}(\partial X,V(\sigma_\lambda,\varphi))$$
to be the adjoint of 
$$\pi_*:C^{ *}(\partial X,V(\tilde{\sigma}_{-\lambda},\tilde\varphi)) \rightarrow C^{*}(B,V_B(\tilde{\sigma}_{-\lambda},\tilde\varphi))\ .$$
\end{ddd}
This definition needs a justification. In fact, 
by Lemma \ref{anal1} the extension exists, is continuous, and  by Fact \ref{adjk} it  depends holomorphically on $\lambda$. Moreover, it is easy to see that the range of the adjoint of $\pi_*$ consists of $\Gamma$-invariant vectors.

We
now consider the left-inverse of $ext$, the restriction 
$$res:{}^\Gamma C^{-\infty}(\partial X,V(\sigma_\lambda,\varphi))\rightarrow C^{-\infty}(B,V_B(\sigma_\lambda,\varphi))\ .$$
The space ${}^\Gamma C^{-\infty}(\Omega,V(\sigma_\lambda,\varphi))$ 
of $\Gamma$-invariant vectors in $C^{-\infty}(\Omega,V(\sigma_\lambda,\varphi))$
can be canonically identified with the corresponding space 
$C^{-\infty}(B,V_B(\sigma_\lambda,\varphi))$.
Composing this identification with the restriction
$res_\Omega:C^{-\infty}(\partial X,V(\sigma_\lambda,\varphi))\rightarrow C^{-\infty}(\Omega,V(\sigma_\lambda,\varphi))$
we obtain the required restriction map $res$. 
 
\begin{lem}\label{lopi}
There exists a continuous map 
$$\widetilde{res}: C^{-\infty}(\partial X ,V(\sigma_\lambda,\varphi))\rightarrow C^{-\infty}(B,V_B(\sigma_\lambda,\varphi))\ ,$$  
which depends holomorphically on $\lambda$ and coincides with $res$ on
${}^\Gamma C^{-\infty}(\partial X,V(\sigma_\lambda,\varphi))$.
\end{lem}
\proof
Let $\tilde\varphi$ be the dual representation of $\varphi$. We exhibit $\widetilde{res}$ as the adjoint of a continuous map
$$\pi^*:C^\infty(B,V_B(\tilde{\sigma}_{-\lambda},\tilde\varphi))\rightarrow C^\infty (\partial X ,V(\tilde{\sigma}_{-\lambda},\tilde\varphi))$$
which depends holomorphically on $\lambda$. Then the lemma follows from Fact
\ref{adjk}.

Let $\{U_\alpha\}$ be a finite open cover of $B$  such that each $U_\alpha$ has a diffeomorphic
lift $\tilde{U}_\alpha \subset\Omega$. Choose a subordinated partition of unity
$\{\chi_\alpha\}$. Pulling $\chi_\alpha$ back to $\tilde{U}_\alpha$ and extending the resulting function by $0$ we obtain a function $\tilde \chi_\alpha\in  C^\infty (\partial X)$. We define 
$$\pi^*(f):=\sum_\alpha \tilde\chi_\alpha f, \quad f\in C^\infty(B,V_B(\tilde{\sigma}_{-\lambda},\tilde\varphi))\ ,$$
where we consider $f$ as an element of 
${}^\Gamma C^{-\infty}(\Omega,V(\tilde{\sigma}_{-\lambda},\tilde\varphi))$.
Then we set $\widetilde{res}:=(\pi^*)^\prime$. \hB

\begin{lem}\label{iden}
We have 
$res\circ ext = \id$. 
\end{lem}
\proof
Recall the definition of $\pi^*$ from the proof of Lemma \ref{lopi}. Then
$\pi_*\pi^*$ is the identity on $C^\infty(B,V_B(\tilde{\sigma}_{-\lambda},\tilde\vp))$.
We obtain
$$res \circ ext=\widetilde{res}\circ ext=  (\pi^*)^\prime \circ (\pi_*)^\prime
=(\pi_*\pi^*)^\prime=\id\ .$$ 
\hB

Let $C^{-\infty}(\Lambda,V(\sigma_\lambda,\varphi))$ denote the space of distribution sections of $V(\sigma_\lambda,\varphi)$ with support in the limit set $\Lambda$.
Since $\Lambda$ is $\Gamma$-invariant
$C^{-\infty}(\Lambda,V(\sigma_\lambda,\varphi))$ is a subrepresentation of 
the representation $\pi^{\sigma,\lambda}\otimes\varphi$ of $\Gamma$ on $C^{-\infty}(\partial X,V(\sigma_\lambda,\varphi))$.
Note that ${}^\Gamma C^{-\infty}(\Lambda,V(\sigma_\lambda,\varphi))$ is exactly the kernel of $res$.

\begin{lem}\label{mainkor}
If  $\:{}^\Gamma C^{-\infty}(\Lambda ,V(\sigma_\lambda,\varphi))=0$ and if $ext$ is defined, then  we have $ext\circ res =\id$.
\end{lem}
\proof 
The assumption implies that $res$ is injective.
By Lemma \ref{iden} we have
$res(ext\circ res - \id)=0$. \hB

In order to apply this lemma we have to check its assumption.
In the course of the paper we will prove several vanishing results for
${}^\Gamma C^{-\infty}( \Lambda ,V(\sigma_\lambda,\varphi))$.
Our first will be
\begin{theorem}\label{pokm}
If $\Ree(\lambda)>\delta_\Gamma+\delta_\varphi$, then   
${}^\Gamma C^{-\omega}( \Lambda ,V(\sigma_\lambda,\varphi))=0$.
\end{theorem}

The proof of this theorem is the first of many instances in this paper where we employ the Poisson transform and its asymptotics.

Let $\gamma$ be a finite-dimensional representation of $K$ on $V_\gamma$
and $T\in \Hom_M(V_\sigma,V_\gamma)$.
We will view sections of $V(\gamma)$ as functions from $G$ to
$V_\gamma$ satisfying the usual $K$-invariance condition. 
\begin{ddd}\label{defofpoi} The Poisson transform 
$$P:=P^T_\lambda:C^{-\infty}(\partial X,V(\sigma_\lambda))\rightarrow C^\infty(X,V(\gamma))$$
is defined by  
$$(P \phi)(g):=\int_K a(g^{-1}k)^{-(\lambda+\rho)}\gamma(\kappa(g^{-1}k)) T \phi(k) dk\ .$$
Here $\phi\in C^{-\infty}(\partial X, V(\sigma_\lambda))$ and the integral is a formal notation
meaning that the distribution $\phi$ has to be applied to the smooth integral kernel.
\end{ddd}

The theory of the Poisson transform is well-known in the case $\sigma=\gamma=1$ (see e.g. \cite{schlichtkrull84}). The general case has been worked out in \cite{vdV}, \cite{yang94}, \cite{olbrichdiss}. The Poisson transform $P$ intertwines
the left-regular representations of $G$ on $C^{-\infty}(\partial X,V(\sigma_\lambda))$ and $C^\infty(X,V(\gamma))$. If $\sigma$ is irreducible, then the image of $P$ consists of joint eigensections with respect to the action of the center $\cZ$ of the universal enveloping algebra of $\gaaa$, where $\cZ$ acts by the infinitesimal character of $\pi^{\sigma,\lambda}$. 
For any finite-dimensional representation $(\varphi, V_\vp)$ of $\Gamma$ we denote
the tensor product $V(\gamma)\otimes V_\vp$ by $V(\gamma,\varphi)$. Then the transformation $P^T_\lambda\otimes\id$, which we will also denote by $P$, intertwines the $\Gamma$-modules $C^{-\infty}(\partial X,V(\sigma_\lambda,\varphi))$
and $C^\infty(X,V(\gamma,\varphi))$.

For the proof of Theorem \ref{pokm} we can assume that $\sigma$ is irreducible.
The asymptotic properties of the Poisson transform will be discussed in detail
at the beginning of section \ref{invvv}. At this point we need the following
two facts. We fix a minimal $K$-type $\gamma$ of the principal series representation $C^{\infty}(\partial X,V(\sigma_\lambda))$ (see \cite{knapp86}, Ch. XV for all that) and choose $0\not=T\in\Hom_M(V_\sigma,V_\gamma)$. Assume that $\Ree(\lambda)>0$. Then for $\phi\in C^{-\infty}(\partial X,V(\sigma_\lambda,\varphi))$ and $f\in C^{\infty}(\partial X,V((\tilde\gamma_{|M})_{-\lambda},\tilde\varphi))$ we have
\begin{equation}\label{asy1}
\lim_{a\to\infty}a^{\rho-\lambda}\int_K \langle P\phi(ka),f(k)\rangle\ dk= c_\sigma(\lambda)\langle \phi,{}^tTf\rangle\ ,
\end{equation}
where $c_\sigma(\lambda)\not=0$ (see Corollarry \ref{lim} below). In particular, $P$ is injective. For any compact $U\subset K\setminus\supp(\phi) M$ there are constants $\epsilon>0$, $C<\infty$ such that for all $k\in U$ and $a\in A_+$
\begin{equation}\label{asy2}  
|P\phi(ka)|\le C a^{\Ree(\lambda)-\rho-\epsilon} \ .
\end{equation}
This follows from Lemma \ref{poisasm}.

We now prove Theorem \ref{pokm} under the additional hypotheses $\Ree(\lambda)>0$. Choose a minimal $K$-type $\gamma$ and an injective $T\in \Hom_M(V_\sigma,V_\gamma)$ as above. Let $\phi\in {}^\Gamma C^{-\infty}(\Lambda,V(\sigma_\lambda,\varphi))$ and $f\in C^{\infty}(\partial X,V((\tilde\gamma_{|M})_{-\lambda},\tilde\varphi))$. We extend $f$ to a section
$\tilde f\in C^\infty(X,V(\tilde\gamma,\tilde\varphi))$ by $\tilde f (k_1ak_2)=\chi(a)\tilde\gamma(k_2)^{-1}f(k_1)$, where $\chi\in C^\infty(A_+)$ is some cut-off function satisfying $\chi\equiv 0$ in a neighbourhood of $1\in A_+$ and $\chi(a)=1$ for $a\gg 1$. 
Let $F$ be a fundamental domain of the action of $\Gamma$ on $X$ such that $\clo(F)\cap\partial X\subset\Omega$ is compact. There is a constant $c\in\R$ such that
\begin{eqnarray*}
\lefteqn{\lim_{a\to\infty}a^{\rho-\lambda}\int_K \langle P\phi(ka),f(k)\rangle\ dk}\\
&=&
c\lim_{a\to\infty}a^{-(\lambda+\rho)} \int_{\{x\in G\:|\: a\le a_x \le 2a\}}  \langle P\phi(x),\tilde f(x)\rangle\ dx\\
&=&c\lim_{a\to\infty}a^{-(\lambda+\rho)}\sum_{g\in\Gamma} \int_{\{x\in G\:|\: a\le a_{gx} \le 2a\}\cap FK}\langle (\id\otimes\varphi(g))P\phi(x),\tilde f(gx)\rangle\ dx\ .
\end{eqnarray*}
Choose $\epsilon>0$ such that the inequality (\ref{asy2}) holds and $\epsilon^\prime:=\Ree(\lambda)-(\delta_\Gamma+\delta_\varphi+2\epsilon)>0$. We now use Corollary \ref{plumps} in order to estimate
for $g\in\Gamma$
\begin{eqnarray*}
\lefteqn{\Big |a^{-(\lambda+\rho)} \int_{\{x\in G\:|\: a\le a_{gx} \le 2a\}\cap FK}\langle (\id\otimes\varphi(g))P\phi(x),\tilde f(gx)\rangle\ dx \:\Big |}\\
&\le& C_0\sup_{k\in K}|f(k)|\: a^{-(\Ree(\lambda)+\rho)}\:a_g^{\delta_\varphi +\epsilon}
\int_{\{x\in G\:|\: a\le a_ga_x\le C_1a\}\cap FK} |P\phi(x)|\: dx\\
&\le& C_2\:  a^{-(\Ree(\lambda)+\rho)}\:a_g^{\delta_\varphi +\epsilon} \: (a_g^{-1}a)^{\Ree(\lambda)+\rho-\epsilon}\\
&=& C_2\: a_g^{-(\delta_\Gamma+\rho+\epsilon^\prime)}\: a^{-\epsilon}\ .
\end{eqnarray*}
The constants $C_i$, $i=0,1,2$ are independent of $g\in\Gamma$.
Since $\sum_{g\in\Gamma} a_g^{-(\delta_\Gamma+\rho+\epsilon^\prime)}$ converges we obtain
$$ |a^{\rho-\lambda}\int_K \langle P\phi(ka),f(k)\rangle\ dk|\le C_3 \: a^{-\epsilon}\ .$$
It follows that
$$ \lim_{a\to\infty}a^{\rho-\lambda}\int_K \langle P\phi(ka),f(k)\rangle\ dk
=0\ .$$
Since $f$ was arbitrary and ${}^tT$ is surjective we conclude from (\ref{asy1}) that $\phi=0$.

It remains to consider the case $\Ree(\lambda)\le 0$. Let $(\pi_\mu, W_\mu)$ be the
irreducible finite-dimensional spherical representation (a representation which contains the trivial $K$-type) with highest weight $\mu\in\aaaa^*$. Then
the highest $\aaaa$-weight space of $W_\mu$ is $P$-invariant, and it is isomorphic
to $V_{1_{-(\mu+\rho)}}$ as a $P$-module. Hence we have an embedding of $P$-representations $V_{\sigma_\lambda}\hookrightarrow 
V_{\sigma_{\lambda+\mu}}\otimes W_{\mu\:|P}$. This leads to an embedding of
the $G$-homogeneous bundles over $\partial X$ 
$$V(\sigma_\lambda)\hookrightarrow V(\sigma_{\lambda+\mu}\otimes\pi_\mu)\cong V(\sigma_{\lambda+\mu})\otimes W_\mu$$ 
as well as of the corresponding spaces of sections. Note that $G$ acts on
$V(\sigma_{\lambda+\mu})\otimes W_\mu$ with the tensor product action.  Tensoring
with $V_\vp$ we obtain an injective intertwining operator of $\Gamma$-modules
\begin{equation}\label{wippe}
i_\mu: C^{-\infty}(\partial X,V(\sigma_{\lambda},\varphi))\hookrightarrow
C^{-\infty}(\partial X,V(\sigma_{\lambda+\mu},\varphi\otimes \pi_\mu))\ ,
\end{equation}
which is local.
Thus, if $\phi\in {}^\Gamma C^{-\infty}(\Lambda,V(\sigma_\lambda,\varphi))$, then $i_\mu(\phi)\in  {}^\Gamma C^{-\infty}(\Lambda,V(\sigma_{\lambda+\mu},\varphi\otimes \pi_\mu))$. Moreover,
$\Ree(\lambda+\mu)>\delta_\Gamma+\delta_\varphi+\mu
=\delta_\Gamma+\delta_{\varphi\otimes \pi_\mu}$. Now choose $\mu$ large enough
such that $\Ree(\lambda+\mu)>0$. It follows from what we have shown above that $i_\mu(\phi)=0$
and hence $\phi=0$. This finishes the proof of Theorem \ref{pokm}.

\section{Meromorphic continuation of $ext$}\label{extttttt}

The extension $ext$ forms a holomorphic
family of maps depending on $\lambda\in\aca$
(we have omitted this dependence in order to simplify
the notation) which is defined now for $\Ree(\lambda)>\delta_\Gamma+\delta_\varphi$.
In the present section we construct the meromorphic continuation of $ext$
to all of $\aca$ for $X\not=\OO H^2$. Our main tool is the scattering matrix
which we introduce below. The scattering matrix for the trivial
group $\Gamma=\{1\}$ is the Knapp-Stein intertwining operator
of the corresponding principal series representation.
We first recall basic properties of the Knapp-Stein intertwining operators.
Then we define the scattering matrix
using the extension and the Knapp-Stein intertwining operators. 
We simultaneously obtain the meromorphic continuations of the scattering matrix and the extension map. 

If $\sigma$ is a representation of $M$, then its Weyl-conjugate $\sigma^w$,
acting on the same vector space $V_\sigma$, is defined by 
$\sigma^w(m):=\sigma( w^{-1} mw)$,
where $w\in N_K(\aaaa)$ is a representative of the non-trivial element of the Weyl group $W(\gaaa,\aaaa)\cong \Z_2$.
The Knapp-Stein intertwining operators form  meromorphic families of $G$-equivariant 
operators (see  \cite{knappstein71} and Lemmas \ref{iny1} and \ref{jjj} below)
$$\hat{J}^w_{\sigma,\lambda}:C^{*}(\partial X,V(\sigma_\lambda)) \rightarrow C^{*}(\partial X ,V(\sigma^w_{-\lambda})),\quad *= -\infty,\infty\ .$$
Here  $\:\hat{}\:$ indicates that $\hat{J}^w_{\sigma,\lambda}$ is unnormalized.

In order to fix our conventions we give a definition of $\hat{J}^w_{\sigma,\lambda}$
as an integral operator acting on smooth functions for $\Ree(\lambda)<0$.
Its continuous extension to distributions is obtained by duality. 
For $\Ree(\lambda)\ge 0$ it is defined by meromorphic continuation.

Consider $f\in C^{\infty}(\partial X,V(\sigma_\lambda))$
as a function on $G$ with values
in $V_{\sigma_\lambda}$ satisfying $f(gp)=\sigma_\lambda(p)^{-1}f(g)$ for all $p\in P$.
For $\Ree(\lambda)<0$ the intertwining operator is defined by the
convergent integral
\begin{equation}\label{furunkel}
(\hat{J}^w_{\sigma,\lambda}f)(g):=\int_{\bar{N}} f(g w\bar{n}) d\bar{n}\ .
\end{equation}

Our first goal is to show that the intertwining operators form
a meromorphic family of operators in the sense defined in Section \ref{anaprep}.
It is an easy consequence of the approach to the intertwining operators
developed by Vogan-Wallach (see \cite{wallach92}, Ch. 10). The only additional point
we have to verify is that in the domain of convergence of (\ref{furunkel}) the
operators $\hat{J}^w_{\sigma,\lambda}$ indeed form a continuous family.
\begin{lem}\label{iny1}
For $\Ree(\lambda)<0$ the intertwining operators 
$$\hat{J}^w_{\sigma,\lambda}:C^{\infty}(\partial X,V(\sigma_\lambda)) \rightarrow C^{\infty}(\partial X ,V(\sigma^w_{-\lambda}))$$
form a holomorphic family of continuous operators.
\end{lem}
\proof
Let $X_i$, $i=1,\dots,\dim(\kaaa)$, be an orthonormal base of $\kaaa$.
For any multiindex $r=(i_1,\dots,i_{\dim(\kaaa)})$ we set $X_r=\prod_{l=1}^{\dim(\kaaa)} X_l^{i_l}$, $|r|=\sum_{l=1}^{\dim(\kaaa)} i_l$,   and for $f\in C^{\infty}(K,V_{\sigma_\lambda})$
we define the seminorm
$$\|f\|_r=\sup_{k\in K} |f(X_rk)|\ .$$
It is well known that the system $\{\|.\|_r\}$, $r$ running over all multiindices,
defines the Fr\'echet topology of $C^{\infty}(K,V_{\sigma_\lambda})$ and by restriction the topology
of $C^{\infty}(\partial X,V(\sigma_\lambda))$.
 
We extend $f\in C^{\infty}(K,V_{\sigma_\lambda})$ to a function $f_\lambda$ on $G$ by setting
$f_\lambda(kan):=f(k)a^{\lambda-\rho}$.
Then we have 
$$\hat{J}^w_{\sigma,\lambda} (f)(k)=\int_{\bar{N}} f_\lambda (k w\bar{n}) d\bar{n}\ .$$
For any $\lambda_0\in \aca$ with $\Ree(\lambda)<0$ and $\delta>0$ we can find
an $\epsilon>0$ such that for $|\lambda-\lambda_0|<\epsilon $
$$ \int_{\bar{N}} |a(\bar{n})^{\lambda_0-\rho}-a(\bar{n})^{\lambda-\rho}| d\bar{n} <\delta\ .$$
We obtain
\begin{eqnarray*}
\|\hat{J}^w_{\sigma, \lambda_0}  f -\hat{J}^w_{\sigma,\lambda}  f \|_r&=&\sup_{k\in K} \int_{\bar{N}} (f_{\lambda_0} (X_rk w\bar{n}) - (f_{\lambda } (X_rk w\bar{n})  ) d\bar{n}\\
&=& \sup_{k\in K} \int_{\bar{N}} f(X_rkw\kappa(\bar{n})) (a(\bar{n})^{\lambda_0-\rho}-a(\bar{n})^{\lambda-\rho})d\bar{n}\\
&\le& \|f\|_r \int_{\bar{N}}   |a(\bar{n})^{\lambda_0-\rho}-a(\bar{n})^{\lambda-\rho}|d\bar{n}\\
&\le& \delta \|f\|_r
\end{eqnarray*}
This immediately implies that $\lambda\mapsto \hat{J}^w_{\sigma,\lambda}$
is a continuous family of operators on the space of smooth functions.
The fact that the family $\hat{J}^w_{\sigma,\lambda}$, $\Ree(\lambda)<0$, depends holomorphically on $\lambda$
is now easy to check (apply \cite{wallach92}, Lemma 10.1.3 and Fact \ref{holla}).
\hB 

\begin{lem}\label{jjj}
The family of intertwining operators 
$$\hat{J}^w_{\sigma,\lambda}:C^{\infty}(\partial X,V(\sigma_\lambda)) \rightarrow C^{\infty}(\partial X ,V(\sigma^w_{-\lambda}))$$
extends   meromorphically to all of $\aca$. 
\end{lem}
\proof
We employ \cite{wallach92}, Thm. 10.1.5, which states that there are polynomial
maps $b:\aca\rightarrow \C$ and $D:\aca\rightarrow \cU(\gaaa)^K$,
such that
\begin{equation}\label{shifty}
b(\lambda)\hat{J}^w_{\sigma,\lambda}=\hat{J}^w_{\sigma, \lambda-4\rho}\circ \pi^{\sigma,\lambda-4\rho}(D(\lambda))\ .
\end{equation}
This formula requires some explanation.
We  identify $$C^{\infty}(\partial X, V(\sigma_\lambda))\cong C^\infty(K,V_\sigma)^M$$
canonically. Then all operators act on the same space $C^\infty(K,V_\sigma)^M$.
 
If we know that $\hat{J}^w_{\sigma,\lambda}$ is meromorphic up to $\Ree(\lambda)<\mu$,
then we conclude that
$$\hat{J}^w_{\sigma,\lambda}= \frac{1}{b(\lambda) } \hat{J}^w_{\sigma,\lambda-4\rho}\circ \pi^{\sigma,\lambda-4\rho}(D(\lambda))$$
is meromorphic up to $\Ree(\lambda)<\mu+4\rho$.
Thus the lemma follows from Lemma \ref{iny1}.
\hB

\begin{lem}\label{off}
Let $\chi,\phi\in C^\infty(\partial X)$ such that $\supp(\phi)\cap \supp(\chi)=\emptyset$.
Then $\chi\hat{J}^w_{\sigma,\lambda}\phi$ is a holomorphic family of smoothing operators.
In particular, the residues of $\hat{J}^w_{\sigma,\lambda}$ are differential
operators.
\end{lem}
\proof 
Since $\supp(\phi)\cap\supp(\chi)=\emptyset$,
there exists a compact set $V\subset \bar{N}$ such that
$$\kappa(\supp(\chi)M w (\bar{N}\setminus V))M \subset (\partial X\setminus \supp(\phi))M\ .$$
For $\Ree(\lambda)<0$ and $f\in C^\infty(\partial X, V(\sigma_\lambda))$
we have (viewing $f$ as a function on $K$ with values in $V_{\sigma_\lambda})$
\begin{eqnarray*}
(\chi \hat{J}^w_{\sigma,\lambda} \phi f)(k) &=& \int_{\bar{N}} \chi(k)f(\kappa(kw\bar{n}))\phi(\kappa(kw\bar{n})) a(\bar{n})^{ \lambda-\rho}d\bar{n}\\
&=& \int_{ V }\chi(k) f(\kappa(kw\bar{n})) \phi(\kappa(kw\bar{n})) a(\bar{n})^{ \lambda-\rho}d\bar{n}\ .
\end{eqnarray*}
The right-hand side of this equation extends to all of $\aca$
and defines a holomorphic family of operators.
This proves the first part of the lemma.
It in particular implies that the residues of $\hat{J}^w_{\sigma,\lambda}$ are local operators.
Hence the second assertion follows.
\hB

We need the following consequence of Lemma \ref{off}.
Let $W\subset \partial X$ be a closed
subset and let
$$\cG_\lambda:=\{f\in C^{-\infty}(\partial X, V(\sigma_\lambda))|
f_{|\partial X\setminus W}\in  C^{\infty}(\partial X\setminus W, V(\sigma_\lambda))\}\ .
$$
We equip $\cG_\lambda$ with the weakest topology such that the maps
$\cG_\lambda\hookrightarrow C^{-\infty}(\partial X, V(\sigma_\lambda))$
and $\cG_\lambda\rightarrow  C^{\infty}(\partial X\setminus W, V(\sigma_\lambda))$
are continuous.
Let $U\subset \bar{U}\subset \partial X\setminus W$ be open.
\begin{kor}\label{dissmo1}
The composition 
$$res_U \circ \hat{J}^w_{\sigma,\lambda}:\cG_\lambda\rightarrow C^\infty(U,V(\sigma^w_{-\lambda}))$$
is well-defined and depends meromorphically on $\lambda$.
\hB
\end{kor}

Below we shall work with a slight modification of $\hat{J}^w_{\sigma,\lambda}$ not depending on the particular representative $w$
and having the intertwininig properties we wish.

If $\sigma$ is equivalent to $\sigma^w$, then we say that $\sigma$ is Weyl-invariant. 
Unless indicated otherwise $\sigma$ shall always denote a Weyl-invariant representation
of $M$ which is either irreducible or of the
form $\sigma^\prime\oplus \sigma^{\prime w}$ with $\sigma^\prime$ irreducible 
and not Weyl-invariant. In both cases the representation of $M$ on $V_\sigma$ can be extended to a representation of $N_K(\aaaa)$ which we also denote by $\sigma$. This extension is unique up to a character of the Weyl group, i.e.,
the two possible choices of $\sigma(w)$ can differ by a sign, only.
We fix such an extension and
define the $G$-intertwining operator $$\hat{J}_{\sigma,\lambda}:C^{\pm\infty}(\partial X,V(\sigma_\lambda)) \rightarrow C^{\pm\infty}(\partial X ,V(\sigma_{-\lambda}))$$
by $\hat{J}_{\sigma,\lambda}:=\sigma(w)\hat{J}^w_{\sigma,\lambda}$. Then the operator $\hat{J}_{\sigma,\lambda}$ does not depend on the choice of $w$.     

In order to define normalized intertwining operators we first
have to recall properties of $c$-functions and the functional equation of $\hat{J}_{\sigma,\lambda}$.
Let $\gamma$ be a finite-dimensional representation of $K$ and let $T\in \Hom_M(V_\sigma,V_\gamma)$.  
We define $T^w\in \Hom_M(V_\sigma,V_\gamma)$ by
$$T^w:=\gamma(w)T\sigma(w^{-1})\ .$$
$T^w$ does not depend on the choice of the representative $w$.
In a similar manner for $T\in \End_M(V_\gamma)$ we define
$T^w:=\gamma(w)T\gamma(w^{-1})$.
Let $c_\gamma: \aca\rightarrow \End_M(V_\gamma)$ be the meromorphic function given for $\Ree(\lambda)>0$ by
\begin{equation}\label{cgamma}
c_{\gamma}(\lambda):=\int_{\bar{N}} a(\bar{n})^{-(\lambda+\rho)}
\gamma(\kappa(\bar{n}))\: d\bar{n}\ .
\end{equation}
\begin{lem}\label{149u}
Let $T\in \Hom_M(V_\sigma,V_\gamma)$ and define
$T^\sharp\in \Hom_K(C^{-\infty}(\partial X,V(\sigma_\lambda)),V_\gamma)$ by
$T^\sharp(f):=(P^T_\lambda f)(1)$, $f\in C^{-\infty}(\partial X,V(\sigma_\lambda))$.
\begin{enumerate}
\item We have
$T^\sharp\circ \hat{J}_{\sigma,\lambda}=[(c_\gamma(\lambda)T)^w]^\sharp$.
\item The Poisson transform satisfies the following functional equation
 \begin{equation}\label{mi9}
P^T_\lambda \circ \hat J_{-\lambda} = P^{(c_{ \gamma}(\lambda)T)^w}_{-\lambda }  \ .\end{equation}
\item 
There is a meromorphic function $p_\sigma:\aca\rightarrow \C$ 
such that
\begin{equation}\label{spex1} 
\hat{J}_{\sigma,-\lambda}\circ \hat{J}_{\sigma,\lambda}=\frac{1}{p_\sigma(\lambda)} \id\ . 
\end{equation}
\item We have
\begin{equation}\label{cgcgw}
c_\gamma(-\lambda)^wc_\gamma(\lambda)T=\frac{1}{p_\sigma(\lambda)} T\ .
\end{equation}
\item The restriction of $\hat{J}_{\tilde{\sigma},\lambda}$ to
$C^\infty(\partial X,V(\tilde{\sigma}_\lambda))$ coincides with
the adjoint of $\hat{J}_{\sigma,-\lambda}$.
\item We have 
\begin{equation}\label{cgcgw1}
c_\gamma(\lambda)^*=c_\gamma(\bar{\lambda})^w\ .
\end{equation}
\end{enumerate}
\end{lem}
\proof
The identity 1. immediately follows from (\ref{mi9}).
The latter can be read off from the asymptotics of the Poisson transforms
(\ref{las1}) and (\ref{las2}). The functional equation
$$\hat{J}^{w^{-1}}_{\sigma^w,-\lambda}\circ \hat{J}^w_{\sigma,\lambda}=
\frac{1}{p_\sigma(\lambda)} \id$$
can be found in  \cite{knappstein71}.
We compute
\begin{eqnarray*}
\hat{J}_{\sigma,-\lambda}\circ \hat{J}_{\sigma,\lambda}&=&
\sigma(w) \circ \hat{J}^w_{\sigma,-\lambda}\circ  \sigma(w)\circ \hat{J}^w_{\sigma,\lambda}\\
&=& \sigma(w^{-1}) \circ \hat{J}^{w^{-1}}_{\sigma,-\lambda}\circ  \sigma(w)\circ \hat{J}^w_{\sigma,\lambda}\\
&=&\hat{J}^{w^{-1}}_{\sigma^{w},-\lambda}\circ  \sigma(w^{-1})\circ  \sigma(w)\circ \hat{J}^w_{\sigma,\lambda}\\
&=&\hat{J}^{w^{-1}}_{\sigma^{w},-\lambda}\circ\hat{J}^w_{\sigma,\lambda}\\
&=&\frac{1}{p_\sigma(\lambda)} \id\ .
\end{eqnarray*}
Equation
(\ref{cgcgw}) can be obtained combining (\ref{spex1}) with
(\ref{mi9}). The relation 5. follows from   \cite{knappstein71}, Lemma 24.
and (\ref{cgcgw1}) is a consequence of 5. and (\ref{cgcgw}).\hB

Explicit formulas for the Plancherel density $p_\sigma$ can be found
e.g. in \cite{knappstein71}, Ch. 12.

For any (irreducible) $\sigma\in\hat{M}$ we fix a minimal $K$-type
$\gamma_\sigma\in\hat{K}$ of $\pi^{\sigma,\lambda}$.
Note that $\Hom_M(V_\sigma,V_{\gamma_\sigma})$ is one-dimensional.
Therefore we can define a meromorphic function $c_\sigma:\aca\mapsto\C$
such that 
\begin{equation}\label{mao}
c_{\gamma_\sigma}(\lambda) T = c_\sigma(\lambda) \left\{\begin{array}{cc} T^w,& \sigma=\sigma^w\\
T,&\sigma\not=\sigma^w\end{array}\right.\ .\end{equation}
If $\sigma$ is of the form $\sigma^\prime\oplus \sigma^{\prime w}$ for
$\sigma^\prime$ not Weyl-invariant, then we define
$c_\sigma(\lambda)\in \End_M(V_\sigma)$
such that it acts on the $M$-isotypic components $V_\sigma(\sigma^\prime)$ ($V_\sigma(\sigma^{\prime w})$) as multiplication by $c_{\sigma^\prime}(\lambda)$ ($c_{\sigma^{\prime w}}(\lambda)$). 

Let now $\sigma$ be a Weyl-invariant representation of $M$
satisfying our general convention introduced above.
We define the normalized intertwining operator by
$$J_{\sigma,\lambda}:= \hat{J}_{\sigma,\lambda}c_{  \sigma}(-\lambda )^{-1}\ .$$
Combining (\ref{spex1}) and (\ref{cgcgw}) we obtain the following
functional equation: 
\begin{equation}\label{spex}
J_{\sigma,-\lambda}\circ J_{\sigma,\lambda} =\id \ .
\end{equation}

We shall often omit the subscript $\sigma$ in the notation of the intertwining operators. Tensoring with a finite-dimensional representation $(\varphi,V_\vp)$ of $\Gamma$ we obtain a
meromorphic-family of $\Gamma$-intertwining operators which we denote
by the same symbol
$$\hat J_\lambda:=\hat{J}_{\sigma,\lambda}\otimes\id: 
C^{\pm\infty}(\partial X,V(\sigma_\lambda,\varphi)) \rightarrow C^{\pm\infty}(\partial X ,V(\sigma_{-\lambda},\varphi))$$
and its normalized version $J_\lambda$. Sometimes we are forced to consider
also intertwining operators for irreducible, non-Weyl-invariant representations
$\sigma$. Then $\hat{J}_{\sigma,\lambda}$ denotes the restriction of
$\hat{J}_{\sigma\oplus\sigma^w,\lambda}$ to $C^{\pm\infty}(\partial X,V(\sigma_\lambda,\varphi))$.

We  now turn to the definition of the (normalized) scattering matrix
as a family of operators
\begin{eqnarray*}
\hat{S}_\lambda  &:&C^{ *}(B,V_B(\sigma_\lambda,\varphi))\rightarrow 
C^{*}(B,V_B(\sigma_{-\lambda },\varphi)),\quad *=\infty,-\infty \\
 S_\lambda  &:&C^{ *}(B,V_B(\sigma_\lambda,\varphi))\rightarrow 
C^{*}(B,V_B(\sigma_{-\lambda },\varphi))\ .
\end{eqnarray*}
\begin{ddd}\label{scatdef}
For $\Ree(\lambda)>\delta_\Gamma+\delta_\varphi$ we define
\begin{equation}\label{scatde}
\hat{S}_{ \lambda} :=res\circ \hat{J}_{\lambda}\circ ext\quad, \quad 
S_{ \lambda} :=res\circ J_{\lambda}\circ ext \ .\end{equation}
\end{ddd}

\begin{lem}\label{kkk}
For $\Ree(\lambda)>\delta_\Gamma+\delta_\varphi$ 
the scattering matrix forms a meromorphic family of operators
$$C^{\pm\infty}(B,V_B(\sigma_\lambda,\vp))\rightarrow C^{\pm\infty}(B,V_B(\sigma_{-\lambda },\vp))\ .$$ If $\hat{S}_\lambda$ is singular and  $\Ree(\lambda)>\delta_\Gamma+\delta_\varphi$, then the residue of $\hat{S}_\lambda$ is a differential operator.
\end{lem}
\proof
The assertion for the scattering matrix acting on distributions follows from
the holomorphy of $ext$,  
Lemma \ref{jjj}, Lemma \ref{off}, Lemma \ref{lopi}, and  Fact \ref{comp}.
By Corollary \ref{dissmo1} it restricts nicely to the space of smooth sections. The last assertion is a consequence of Lemma \ref{off}.
\hB

\begin{lem}\label{lok}
If  $\Ree(\lambda)>\delta_\Gamma+\delta_\varphi$, then the adjoint
$${}^tS_{ \lambda}: C^{\infty}(B,V_B(\tilde{\sigma}_\lambda,\tilde\varphi))\rightarrow C^{\infty}(B,V_B(\tilde{\sigma}_{-\lambda },\tilde\varphi))$$
of
$$S_{ \lambda}: C^{-\infty}(B,V_B(\sigma_\lambda,\varphi))\rightarrow C^{-\infty}(B,V_B(\sigma_{-\lambda },\varphi))$$
coincides with the restriction of 
$$S_\lambda: C^{-\infty}(B,V_B(\tilde{\sigma}_\lambda,\tilde\varphi))\rightarrow C^{-\infty}(B,V_B(\tilde{\sigma}_{-\lambda },\tilde\varphi))$$ 
to $C^{\infty}(B,V_B(\tilde{\sigma}_\lambda,\tilde\varphi))$.
\end{lem}
\proof
We employ the fact (Lemma \ref{149u}, 5.) that the corresponding relation holds for the intertwining operators (step (\ref{zx2}) below).
Recall the definition of $\pi^*$ from the proof of Lemma \ref{lopi}.
In the domain of convergence of $\pi_*$ we have
$$ \langle \phi,\pi_*(h)\rangle=\sum_{g\in\Gamma} \langle \pi^*(\phi),\pi(g)h\rangle\ , 
$$
where $\phi\in C^{\infty}(B,V_B(\sigma_\lambda,\varphi))$, $h\in C^{\infty}(\partial X,V(\tilde\sigma_{-\lambda},\tilde\varphi))$, and
$\pi(g)=\pi^{\tilde\sigma,-\lambda}(g)\otimes \tilde\vp(g)$. We will use this
formula in step (\ref{kairo}) below.

Let $\phi\in C^{\infty}(B,V_B(\sigma_\lambda,\varphi))$,
$f\in C^{\infty}(B,V_B(\tilde{\sigma}_\lambda,\tilde\varphi))$,
and consider $\phi$ as a distribution section.
Then
\begin{eqnarray}
\langle \phi ,{}^tS_\lambda f\rangle &=&\langle  S_\lambda \phi,  f\rangle \nonumber\\
&=& \langle  (res\circ J_\lambda\circ ext) \phi ,  f\rangle\nonumber\\
&=&\langle  (\tilde{res}\circ J_\lambda\circ ext) \phi ,  f\rangle\nonumber\\
&=&\langle (J_\lambda\circ ext) \phi , \pi^* f   \rangle\nonumber\\ 
&=&\langle  ext \phi , ({}^t J_\lambda \circ\pi^*)f  \rangle \nonumber\\
&=&\langle   ext \phi , (J_\lambda \circ\pi^*) f  \rangle\label{zx2}\\
&=&\langle   \phi , (\pi_*\circ J_\lambda \circ\pi^*) f  \rangle\nonumber\\
&=&\sum_{g\in\Gamma } \langle     \pi^*\phi ,   (\pi(g)\circ J_\lambda \circ\pi^*)   f\rangle \label{kairo}\\
&=&\sum_{g\in\Gamma } \langle     \pi^*\phi ,   ( J_\lambda\circ\pi(g) \circ\pi^*)   f\rangle \nonumber\\
&=&\sum_{g\in\Gamma } \langle     (\pi(g)\circ J_\lambda \circ\pi^*)\phi , \pi^*   f\rangle \label{zx3}\\
&=&   \langle \phi , S_\lambda f\rangle\ . \label{zx4}
\end{eqnarray}
In order to obtain (\ref{zx4}) from (\ref{zx3}) we do the transformations
backwards with the roles of $\phi$ and $f$ interchanged.
\hB

\begin{lem}\label{pfun}
If $|\Ree(\lambda)|<- (\delta_\Gamma+\delta_\varphi)$,
then the scattering matrix satisfies the functional equation
(viewed as an identity of meromorphic families of operators)
$$S_{-\lambda}\circ S_\lambda = \id\ .$$
\end{lem}
\proof
We employ Lemma \ref{mainkor}, Theorem \ref{pokm}, and (\ref{spex}) in order to compute 
\begin{eqnarray*}
S_{ -\lambda }\circ S_{ \lambda}&=&res\circ J_{ -\lambda }\circ ext\circ res \circ J_{ \lambda}\circ ext\\
&=&res\circ J_{ -\lambda }\circ J_{ \lambda}\circ ext \\
&=&res\circ ext\\
&=& \id\ .
\end{eqnarray*}
\hB

The main result of this section is 

\begin{theorem}\label{part1}
Let $X\not=\OO H^2$.
Then the scattering matrix
$${S}_\lambda :C^{\pm\infty}(B,V_B(\sigma_\lambda,\varphi))\rightarrow C^{\pm\infty}(B,V_B(\sigma_{-\lambda },\varphi))$$
and  the extension map 
$$ext: C^{-\infty}(B,V_B(\sigma_\lambda,\varphi))\rightarrow {}^\Gamma C^{-\infty}(\partial X,V(\sigma_\lambda,\varphi))\ ,$$
initially defined for $\Ree(\lambda)>\delta_\Gamma+\delta_\varphi$,
have  meromorphic continuations to 
all of $\aca$. 
In particular, we have    
\begin{eqnarray}\label{forme1}
ext & =& J_{-\lambda} \circ ext \circ S_{\lambda}\ , \\
S_{-\lambda}\circ S_\lambda &=& \id\ .\label{forme}
\end{eqnarray}  
Moreover,
$ext$ has at most finite-dimensional singularities. 
\end{theorem}

The proof of the theorem will occupy the remainder of this section. We first use the meromorphic Fredholm theory in order to understand a special case. We introduce the element $\beta\in\aaaa^*$ as follows\\
\centerline{\begin{tabular}{|c||c|c|c|c|}
\hline
$X$&$\R H^n$&$\C H^n$&$\HH H^n$&$\OO H^n$\\
\hline
$\beta$&$0$&$0$&$2\alpha$&$6\alpha$\\
\hline
\end{tabular}\ . }
 
\begin{lem}\label{kugel}
If $\delta_\Gamma+\delta_\varphi<0$ and $\hat M\ni\sigma=1$ is the trivial representation, then $S_\lambda$ and $ext$ have meromorphic continuations to
$W:=\{\lambda\in\aca\:|\:\Ree(\lambda)>-\rho+\beta\}$. On $W$ the remaining assertions of Theorem
\ref{part1} hold true. 
\end{lem}
\proof
We construct the meromorphic continuation of 
$$S_\lambda:C^{\infty}(B,V_B(1_\lambda,\vp))\rightarrow C^{\infty}(B,V_B(1_{-\lambda },\vp))\ ,$$ 
and then we extend this continuation to distributions
by duality using Lemma \ref{lok}. The idea is to set $S_\lambda:=S_{-\lambda}^{-1}$ for $\Ree(\lambda)<-(\delta_\Gamma+\delta_\varphi)$ and to show that $S_{-\lambda}^{-1}$ forms
a meromorphic family on $$U:=\{\lambda\in\aca\:|\:\max\{\delta_\Gamma+\delta_\varphi,-\rho+\beta\}<\Ree(\lambda)<\rho-\beta\}\ .$$ 

Let $\{U_\alpha\}$ be a finite open covering of $B$ and let $\tilde{U}_\alpha$ be
diffeomorphic lifts of $U_\alpha$.
Choose a subordinated partition of unity $\{\phi_\alpha\}$.
Let $\tilde\phi_\alpha$ be the corresponding compactly supported function on $\tilde{U}_\alpha$. 
Then we define
$\chi\in {}^\Gamma C^\infty(\Omega\times\Omega)$ by 
$$\chi(x,y):=\sum_\alpha\sum_{g\in\Gamma} \tilde\phi_\alpha(gx)\tilde\phi_\alpha(gy)\ .$$ 
Let 
$${J}^{diag}_{\lambda}:C^\infty(B, V_B(1_\lambda,\varphi))\rightarrow C^\infty(B, V_B(1_{-\lambda},\varphi))$$ be the meromorphic family
of operators obtained by multiplying
the distribution kernel of $\hat{J}_{ \lambda}$ by $\chi$.
If $f\in C^\infty(B, V_B(1_\lambda,\varphi))$, then
$$({J}^{diag}_{ \lambda})f =\sum_{\alpha} \tilde\phi_\alpha {J}_\lambda(\tilde\phi_\alpha f)$$
using the canonical identifications (see the proof of Lemma \ref{lopi}).

For $\lambda\in U$ define 
\begin{equation}\label{klio}
R(\lambda):= {J}^{diag}_{-\lambda} \circ {S}_\lambda - \id\ .
\end{equation}
The inverse of the normalized scattering matrix for $\lambda\in U$
should be  given by
\begin{equation}\label{finit}
{S}_{ \lambda}^{-1}=(\id+ R(\lambda))^{-1}\circ {J}^{diag}_{ -\lambda }\ .
\end{equation}
It exists as a meromorphic family if $(\id+ R(\lambda))^{-1}$ does.

We apply the  meromorphic Fredholm theory (Lemma \ref{merofred})
in order to invert $\id +R(\lambda)$ for $\lambda\in U$
and to conclude that $(\id+R(\lambda))^{-1}$
is meromorphic. 
 
We check the assumption of Lemma \ref{merofred}.
We choose a Hermitian metric on $V_B(1_0,\varphi)$ and a volume form on $B$.
The Hilbert space $\cH$ of Proposition \ref{merofred}
is $L^2(B,V_B(1_0,\varphi))$ defined using these choices.  
The Fr\'echet space $\cF$ is just $C^\infty(B, V_B(1_0,\varphi))$.
Implicitly, we identify the spaces $C^\infty(B, V_B(1_\lambda,\varphi))$ with $C^\infty(B, V_B(1_0,\varphi))$ using a trivialization of the holomorphic
family of bundles $\{V_B(1_\lambda,\varphi)\}_{\lambda\in\aca}$.
 
We need the following well-known fact.
\begin{lem}\label{rund}
If $|\Ree(\lambda)|<\rho-\beta$, then the spherical intertwining operator
$$J_{1,\lambda}: C^{\infty}(\partial X,V(1_\lambda)) \rightarrow C^{\infty}(\partial X ,V(1_{-\lambda}))$$ is regular and invertible. Moreover,
$J_{1,0}=\id$.
\end{lem}
\proof  
If $|\Ree(\lambda)|<\rho-\beta$, then the principal series representation $C^{\infty}(\partial X,V(1_\lambda))$ is irreducible (see e.g. \cite{helgason94}, Ch. VI, Thm. 3.6). By definition of the normalized intertwining operator its restriction to the minimal $K$-type is regular and bijective. Now the kernel of the leading term in the Laurent expansion of $J_{1,\mu}$ at $\mu=\lambda$ is a subrepresentation of $C^{\infty}(\partial X,V(1_\lambda))$.
It follows that all singular terms  in the Laurent expansion vanish and that $J_{1,\lambda}$ is bijective.
By Schur's Lemma there is $\mu\in\C$ such that  $J_{1,0}= \mu\, \id$.
Our  normalization of $J_{1,.}$ implies that $\mu=1$. 
\hB

It follows that ${J}^{diag}_{-\lambda}$ as well as
${S}_\lambda=res\circ {J}_\lambda\circ ext$ are regular on $U$.
The difference $J^{off}_{-\lambda}:=res\circ J_{-\lambda}-J^{diag}_{-\lambda}\circ res$
is a continuous map from ${}^\Gamma C^{-\infty}(\partial X, V(1_{-\lambda},\vp))$ to
$C^{\infty}(B, V_B(1_{\lambda},\vp))$ by Lemma \ref{off}. Since $R(\lambda)=-J^{off}_{-\lambda}\circ J_\lambda\circ ext$, the family  
$R(\lambda)$ is a holomorphic family of smoothing operators on $U$. In addition,
$R(0)=0$. Thus $\id +R(\lambda)$ is invertible at $\lambda=0$.
 
We now have verified the assumptions of Lemma \ref{merofred}.
We conclude that 
the family  ${S}_{\lambda}^{-1}$ is meromorphic family on $U$ with finite-dimensional singularities.

Since $\delta_\Gamma+\delta_\varphi<0$, we have  
$-U\cup \{\lambda\in\aca\:|\:\Ree(\lambda)>\delta_\Gamma+\delta_\varphi\}=W$.
Furthermore, by Lemma \ref{pfun} we have $S_\lambda=S_{-\lambda}^{-1}$ on $-U\cap \{\lambda\in\aca\:|\:\Ree(\lambda)>\delta_\Gamma+\delta_\varphi\}$. Thus, setting
$S_\lambda:=S_{-\lambda}^{-1}$ for $-\lambda\in U$ we obtain a well-defined continuation of $S_\lambda$ to all of $W$.
By duality  this  continuation extends to distributions still having finite-dimensional
singularities. Moreover, the functional equation (\ref{forme}) holds by definition.

It remains to consider the extension map. 
We  employ the scattering matrix in order to define  for 
$\lambda\in W$,
$\Ree(\lambda)<-(\delta_\Gamma+\delta_\varphi)$
$$ ext_1 := J_{-\lambda} \circ ext \circ S_{\lambda} \ .$$
We claim that $ext=ext_1$. In fact, since $res$ is injective on  $\{|\Ree(\lambda)| < -(\delta_\Gamma+\delta_\varphi)\}$ by Theorem \ref{pokm}, the computation 
$$
res\circ ext_1= res\circ J_{-\lambda} \circ ext \circ S_\lambda
=  S_{-\lambda}\circ S_\lambda   
=\id
$$
and Lemma \ref{iden} imply the claim.

We now have constructed a meromorphic continuation of $ext$ to all of 
$W$. The relation (\ref{forme1})
between the scattering matrix and $ext$ follows by meromorphic continuation. This equation also implies that $ext$  has at most finite-dimensional singularities.
Thus the proof of Lemma \ref{kugel} is complete.
\hB

In the next step we drop the assumption $\delta_\Gamma+\delta_\varphi<0$ employing the fact that for $X\not=\OO H^2$ the symmetric space $X$ belongs to a series.

\begin{lem}\label{bed}
If $X\not=\OO H^2$, then Lemma \ref{kugel} holds true without the assumption
$\delta_\Gamma+\delta_\varphi<0$.
\end{lem} 
\proof
$X$ belongs to a series of rank-one symmetric
spaces.
Let $\dots \subset G^n\subset G^{n+1}\subset \dots$ be the corresponding sequence 
of real, semisimple, linear Lie groups inducing embeddings of the corresponding
Iwasawa constituents $K^n\subset K^{n+1}$,
$N^n \subset  N^{n+1}$, $ M^n\subset M^{n+1}$
such that $A=A^n = A^{n+1}$.
Then we have totally geodesic embeddings of the symmetric spaces
$X^n\subset X^{n+1}$
inducing embeddings of their boundaries
$\partial X^n\subset \partial X^{n+1}$.
If $\Gamma\subset G^n$ is convex-cocompact then it is still convex-cocompact
viewed  as a subgroup of $G^{n+1}$.
We obtain  embeddings
$\Omega^n\subset  \Omega^{n+1}$ inducing
$B^n\subset  B^{n+1}$
while the limit set $ \Lambda^n$ is identified with $ \Lambda^{n+1}$.
Let $\rho^n(H)=\frac{1}{2}\tr(\ad(H)_{|\naaa^n})$, $H\in\aaaa$.

The exponent of $\Gamma$ now depends on $n$ and is denoted by
$\delta_\Gamma^n$.
We have the relation $\delta_\Gamma^{n+1}=\delta_\Gamma^n-\zeta$, where $\zeta:=\rho^{n+1}-\rho^n>0$.
Thus $\delta_\Gamma^{n+m}\to -\infty$ as $m\to \infty$. Hence,
taking $m$ large enough we obtain $\delta_\Gamma^{n+m}+\delta_\varphi<0$.
The aim of the following discussion is to show how the meromorphic continuation of $ext^{n+1}$ leads to the continuation of  $ext^n$.

Let $P^{n}:=M^nA^nN^n$, $V(1_\lambda,\varphi)^n:=G^n\times_{P^n} V_{1_\lambda}\otimes V_\vp$, and $V_{B^n}(1_\lambda,\varphi)=\Gamma\backslash V(1_\lambda,\varphi)^n_{|\Omega^n}$. Here as always $(\varphi,V_\vp)$ is a finite-dimensional representation representation of $\Gamma$.
The representation $V_{1_\lambda}$ of $P^{n+1}$ restricts to the representation
$V_{1_{\lambda-\zeta}}$ of $P^n$.
This induces isomorphisms of bundles
$$V(1_\lambda,\varphi)^{n+1}_{|\partial X^n}\cong V(1_{\lambda-\zeta},\varphi)^n,\quad  V_{B^{n+1}}(1_\lambda,\varphi)_{|B^n}\cong V_{B^n}(1_{\lambda-\zeta},\varphi)\ .$$

Let 
\begin{eqnarray*}
i^*:C^\infty(B^{n+1},V_{B^{n+1}}(1_\lambda,\tilde\varphi))&\rightarrow& C^\infty(B^{n },V_{B^{n}}(1_{\lambda-\zeta},\tilde\varphi))\ ,\\
j^*:C^\infty(\partial X^{n+1},V(1_\lambda,\tilde\varphi)^{n+1})&\rightarrow& C^\infty(\partial X^n,V(1_{\lambda-\zeta},\tilde\varphi)^n)
\end{eqnarray*}
denote the maps given by restriction of sections.
Note that $j^*$ is $G^n$-equivariant.
The adjoint maps define
the push-forward of distribution sections
\begin{eqnarray*}
i_*:C^{-\infty}(B^n,V_{B^n}(1_{\lambda},\varphi))&\rightarrow& C^{-\infty}(B^{n+1},V_{B^{n+1}}(1_{\lambda-\zeta},\varphi))\ ,\\
j_*:C^{-\infty}(\partial X^n,V(1_{\lambda},\varphi)^n)&\rightarrow& C^{-\infty}(\partial X^{n+1},V(1_{\lambda-\zeta},\varphi)^{n+1})\ .
\end{eqnarray*}

If $\phi\in C^{-\infty}(B^n,V_{B^n}(1_\lambda,\varphi))$, then the push forward
$i_\ast \phi $ has support in $B^n\subset B^{n+1}$.
Since $res^{n+1}\circ ext^{n+1}=\id$
we have 
$\supp(ext^{n+1}\circ i_*)(\phi)\subset 
\Lambda^{n+1}\cup \Omega^{n}=\partial X^n$.

Assume that $ext^{n+1}$ is meromorphic on $W^{n+1}:=\{\lambda\in\aca\:|\:\Ree(\lambda)>-\rho^{n+1}+\beta\}$.
We are now going to continue $ext^n$ using $i_*$, $ext^{n+1}$ and a left inverse of $j_*$. 
As in previous occasions we trivialize the family $\{V(1_\lambda)^{n+1}\}_\lambda$. 
We identify $C^{\infty}(\partial X^{n+1},V(1_\lambda,\tilde\varphi)^{n+1})$
with $C^{\infty}(\partial X^{n+1})\otimes\tilde V_\vp$ for all $\lambda\in\C$.
Let $U\subset \partial X^{n+1}$ be a tubular neighbourhood
of $\partial X$ and fix a diffeomorphism $T:(-1,1) \times \partial X^n\stackrel{\cong}{\rightarrow} U$. Let $\chi\in C_c^\infty((-1,1))$ be a cut-off function satisfying $\chi(0)=1$. Then we define
a continuous extension $t:C^{\infty}(\partial X^n)\otimes V_{\tilde\vp}\rightarrow C^\infty(\partial X^{n+1})\otimes V_{\tilde\vp}$
by 
$$tf(x):=\left\{
\begin{array}{llll}
\chi(r)f(x^\prime)&&\mbox{if }x=T(r,x^\prime)\in U,&r\in (-1,1),\: x^\prime\in\partial X^n\\
0&&\mbox{if }x\not\in U
\end{array}\right.\ .$$ 
Let
$t^\prime:\:C^{-\infty}(\partial X^{n+1},V(1_{\lambda-\zeta},\varphi)^{n+1})\rightarrow C^{-\infty}(\partial X^n,V(1_{\lambda},\varphi)^n)$ be the adjoint of $t$. Then $t^\prime\circ j_*=\id$. Now we can define
$$\widetilde{ext}^n \phi:= (t^\prime \circ ext^{n+1}\circ i_*)(\phi)\ .$$
Then
$$\widetilde{ext}^n : C^{-\infty}(B^n,V_{B^n}(1_\lambda,\varphi)) \rightarrow C^{-\infty}(\partial X^n,V(1_\lambda,\varphi)^n)$$
is a meromorphic family on $W^{n+1}+\zeta=W^n$ of continuous maps with at most finite-dimensional singularities. 

In order to prove that $\widetilde{ext}^n$ provides the desired
meromorphic continuation it remains to show that it coincides with ${ext}^n$ 
in the region $\Ree(\lambda)>\delta_\Gamma^{n}+\delta_\varphi$. If $\Ree(\lambda)>\delta_\Gamma^{n}+\delta_\varphi$, then $\Ree(\lambda)-\zeta>\delta_\Gamma^{n+1}+\delta_\varphi$,
and the push-down maps $\pi_{*,-\lambda}^n$, $\pi_{*,-\lambda+\zeta}^{n+1}$ are defined.
It is easy to see from the definition of the push-down that in the domain of convergence
$$i^*\circ \pi_{*}^{n+1}=\pi_{*}^n\circ j^*\ .
$$
Taking adjoints we obtain $ext^{n+1}\circ i_*
=j_*\circ {ext}^n$.
Therefore we have
$$\widetilde{ext}^n
=t^\prime \circ ext^{n+1}\circ i_*
=t^\prime \circ j_*\circ {ext}^n
={ext}^n\ .$$
It follows by meromorphy that $\im(\widetilde{ext}^n)$ consists of $\Gamma$-invariant sections for all $\lambda\in W^n$.
 
We define the meromorphic continuation of the scattering matrix by
(\ref{scatde}).
Then it is easy to see that the scattering matrix has the properties
as asserted.
This finishes the proof of Lemma \ref{bed}.
\hB

We now use tensoring with finite-dimensional $G$-representations in order to
complete the proof of Theorem \ref{part1}. For a moment let $\sigma\in\hat M$.
Then the theory of highest weights for $G$ implies that there are sequences $\mu_i\in\aaaa^*$, $\mu_i\to\infty$ and $\pi_{\sigma,\mu_i}$ of
finite-dimensional irreducible representations of $G$ such that $\pi_{\sigma,\mu_i}$ has highest $\aaaa$-weight $\mu_i$, and the representation of $M$ on the corresponding highest weight space is equivalent to $\sigma$. 
More details on this can be found e.g. in \cite{bunkeolbrich955}, pp. 39-41. 
If $\sigma=\sigma^\prime\oplus\sigma^{\prime w}$ for some non-Weyl-invariant $\sigma^\prime\in\hat M$, then we set $\pi_{\sigma,\mu_i}
:=\pi_{\sigma^{\prime},\mu_i}\oplus\pi_{\sigma^{\prime w},\mu_i}$. 
As in the proof of Theorem \ref{pokm}, Equation (\ref{wippe}), we obtain an embedding of bundles
\begin{equation}\label{sub}
V(\sigma_{\lambda},\varphi)\hookrightarrow V(1_{\lambda+\mu_i},\varphi\otimes \pi_{\sigma,\mu_i})
\end{equation}
as well as
an injective, local, $\Gamma$-intertwining operator
\begin{equation}\label{romski}
i_{\sigma,\mu_i}: C^{-\infty}(\partial X,V(\sigma_{\lambda},\varphi)) \hookrightarrow
C^{-\infty}(\partial X,V(1_{\lambda+\mu_i},\varphi\otimes \pi_{\sigma,\mu_i}))\ .
\end{equation}
It induces the embedding 
$$i_{\sigma,\mu_i}^B: C^{-\infty}(B,V_B(\sigma_{\lambda},\varphi)) \hookrightarrow
C^{-\infty}(B,V_B(1_{\lambda+\mu_i},\varphi\otimes \pi_{\sigma,\mu_i}))\ .$$
Let ${i_{\sigma,\mu_i}}^\prime$ and ${i_{\sigma,\mu_i}^B}^\prime$ be the adjoint operators, i.e., the projections onto the spaces smooth sections of the corresponding dual bundles. Then in the domain of convergence of the push-down map $\Ree(\lambda)<-(\delta_\Gamma+\delta_\varphi)$ we have
$$ {i_{\sigma,\mu_i}^B}^\prime\circ\pi_*=\pi_*\circ {i_{\sigma,\mu_i}}^\prime\ .$$
Note that the domains of convergence of both sides coincide.
It follows that for $\Ree(\lambda)>\delta_\Gamma+\delta_\varphi$, 
$\phi\in  C^{-\infty}(B,V_B(\sigma_{\lambda},\varphi))$, \begin{equation}\label{hope}
ext\circ i_{\sigma,\mu_i}^B(\phi)=i_{\sigma,\mu_i}\circ ext(\phi)\ .
\end{equation}

Now let $\nu\in\aaaa*$ be arbitrary, and let $W_\nu:=\{\lambda\in\aca\:|\:\Ree(\lambda)>\nu\}$ be the corresponding half-plane. Choose $\mu_i$ large enough such that $\Ree(\lambda)+\mu_i>-\rho+\beta$ for all $\lambda\in W_\nu$. By Lemma \ref{bed}
the extension map on the left hand side of (\ref{hope}) has a meromorphic continuation to $W_\nu$ with finite-dimensional singularities. Moreover, 
$$ext\circ i_{\sigma,\mu_i}^B\left(C^{-\infty}(B,V_B(\sigma_{\lambda},\varphi))\right)\subset
i_{\sigma,\mu_i}\left(C^{-\infty}(\partial X,V(1_{\lambda+\mu_i},\varphi\otimes \pi_{\sigma,\mu_i}))\right)\ .$$ 
In fact, this is true for $\lambda$ in the domain
of convergence by (\ref{hope}), hence on all of $W_\nu$ by meromorphy. Therefore
we can define for $\phi\in  C^{-\infty}(B,V_B(\sigma_{\lambda},\varphi))$, $\lambda\in W_\nu$ 
$$ext(\phi):= (i_{\sigma,\mu_i})^{-1}\circ ext \circ i_{\sigma,\mu_i}^B(\phi)\ .$$
This gives the desired meromorphic continuation of $ext$ to $W_\nu$. Now
(\ref{hope}) holds on all of $W_\nu$. It follows that the singularities of $ext$ are at most finite-dimensional. Since $\nu$ was arbitrary $ext$ is meromorphic on all of $\aca$.

We define the meromorphic continuation of the scattering matrix by
(\ref{scatde}).
Using Lemma \ref{mainkor} and Theorem \ref{pokm} it is easy to verify the functional equations (\ref{forme}) and (\ref{forme1}) for $\Ree(\lambda)<-(\delta_\Gamma+\delta_\varphi)$. By meromorphy they hold on all of $\aca$. This finishes the proof of Theorem \ref{part1}.
\hB

\section{Invariant distributions on the limit set}\label{invvv}

In the present section we study the spaces 
${}^\Gamma C^{-\infty}(\Lambda,V(\sigma_\lambda,\varphi))$ of
invariant distributions which are
supported on the limit set. The main result of this section is

\begin{theorem}\label{disfin}
Let $X\not=\OO H^2$. Fix finite-dimensional representations $\sigma$ of $M$ and $\vp$ of
$\Gamma$. Then
\begin{enumerate}
\item The set $\{\lambda\in \aca\:|\:{}^{\Gamma} C^{-\infty}(\Lambda,V(\sigma_\lambda,\varphi))\not=0\}$ is discrete.
\item For each $\lambda\in\aca$ the space ${}^{\Gamma} C^{-\infty}(\Lambda,V(\sigma_\lambda,\varphi))$ is finite-dimensional.
\end{enumerate}
\end{theorem}

The proof is based on the following observation: Assume that $ext$ is singular at
$\lambda\in\aca$. Since $res\circ ext(\phi_\mu)=\phi_\mu$ is regular for any holomorphic family $\mu\mapsto\phi_\mu\in C^{-\infty}(B,V_B(\sigma_\mu,\varphi))$ the leading singular part of the Laurent expansion of 
$ext(\phi_\mu)$ at $\mu=\lambda$ belongs to ${}^\Gamma C^{-\infty}(\Lambda,V(\sigma_\lambda,\varphi))$.
We show that for almost all $\lambda$  these spaces are generated by the singular
parts of $ext$. 

First we need detailed knowledge of the asymptotics of the Poisson tansform.
\begin{lem}\label{poisasm}
 Let $\gamma$ be a finite-dimensional representation of $K$, $T\in \Hom_M(V_\sigma,V_\gamma)$, and let $w\in N_K(\aaaa)$  represent the non-trivial element of the Weyl group
$W(\gaaa,\aaaa)$.  
\begin{enumerate}
\item Let
$f\in C^{\infty}(\partial X,V(\sigma_\lambda))$. If $\Ree(\lambda)>0$, then
there exists $\epsilon>0$ (depending on $\lambda$ but not on $f$) such that for $a\to\infty$  we have
$$(P^T_\lambda f)(ka)=a^{\lambda-\rho} c_\gamma(\lambda)T  f(k)+O(a^{\Ree(\lambda)-\rho-\epsilon})$$
uniformly in $k\in K$.
If $|\Ree(\lambda)|<\frac{1}{2}|\alpha|$ and $\lambda\not=0$, then \\
\begin{equation}\label{murphy}
(P^T_\lambda f)(ka) = a^{\lambda-\rho} c_\gamma(\lambda)  Tf(k)+ a^{-\lambda-\rho}  T^w(\hat{J}_\lambda f)(k) + O(a^{-\frac{\alpha}{2}-\rho-\epsilon})
\end{equation}
uniformly in $k\in K$.
The remainder depends jointly continuously on $\lambda$ and $f$.

If $\hat J_0$ is regular, then (\ref{murphy}) remains valid for $\lambda=0$. Otherwise we have
\begin{eqnarray*}
(P^T_0 f)(ka)&=&\log(a) a^{-\rho}
\left(2\res_{\lambda=0}(c_\gamma(\lambda))T f(k)\right)\\
&&\quad\quad+a^{-\rho}\left(c_\gamma^0Tf(k)+T^w(\hat J_0^0f)(k)\right)+ O(a^{ -\rho-\alpha})\ ,
\end{eqnarray*}
where $c_\gamma^0$, $\hat J_0^0$ denote the constant terms in the Laurent expansions
of $c_\gamma$ and $\hat J_\lambda$ at $\lambda=0$.\\
\item Let $\partial X=U\cup Q$, where $U$ is open and $Q:=\partial X\setminus U$.
Let 
$f\in C^{-\infty}(\partial X,V(\sigma_\lambda))$ with $\supp f\subset Q$. Then there exist smooth functions $\psi_n$, $n\in\nat$,
on $U$ such that
\begin{equation}\label{epan}
(P^T_\lambda f)(ka)=a^{-(\lambda+\rho)}T^w(\hat{J}_\lambda f)(k) 
 +\sum_{n\ge 1} a^{-(\lambda+\rho)-n\alpha}\psi_n(k)\ ,\quad k\in U\ .
\end{equation}
The series converges uniformly for $a\gg 0$ and $k$ in compact subsets of $U$. 
In particular, for $a\to\infty$ we have 
$$(P^T_\lambda f)(ka)=a^{-\lambda-\rho}T^w(\hat{J}_\lambda f)(k)+O(a^{-\lambda-\rho-\alpha})$$
uniformly as $kM$ varies in compact subsets of $U$.\\
\item  
Let $U,Q$ be as in 2. and 
$f\in C^{-\infty}(\partial X,V(\sigma_\lambda))$ such that
$res_U f\in C^\infty(U,V(\sigma_\lambda))$. If $\Ree(\lambda)>0$, then
there exists $\epsilon>0$ such that for $a\to\infty$  we have
$$(P^T_\lambda f)(ka)=a^{\lambda-\rho} c_\gamma(\lambda)T  f(k)+O(a^{\lambda-\rho-\epsilon})$$
uniformly as $kM$ varies in compact subsets of $U$.
If $|\Ree(\lambda)|<\frac{1}{2}|\alpha|$, $\lambda\not=0$, then
we have for $a\to\infty$  
$$(P^T_\lambda f)(ka) = a^{\lambda-\rho} c_\gamma(\lambda)T  f(k)+ a^{-\lambda-\rho}  T^w(\hat{J}_\lambda f)(k) + O(a^{-\frac{\alpha}{2} -\rho-\epsilon})$$ 
uniformly as $kM$ varies in compact subsets of $U$.
The remainder depends jointly continuously on $\lambda$ and $f$.
\end{enumerate}
The asymptotic expansions can be differentiated with respect to $a$.
\end{lem}
\proof
Assertion 1 is a consequence of the general results concerning
the asymptotics of matrix coefficients of admissible representations including their dependence on parameters 
(\cite{wallach88}, Thm. 4.4.3, \cite{wallach92}, 12.4., 12.5., 12.6., \cite{olbrichdiss}) combined with the limit
formulas for the Poisson transform (see \cite{vdV} or \cite{olbrichdiss}, also
\cite{wallach88}, Thm. 5.3.4)
\begin{eqnarray}
\lim_{a\to\infty} a^{\rho-\lambda}
(P^T_\lambda f)(ka)&=&c_\gamma(\lambda)Tf(k),\qquad \Ree(\lambda)>0\ ,\label{las1}\\
\lim_{a\to\infty}a^{\rho+\lambda} 
(P^T_\lambda f)(ka)&=&T^w(\hat{J}_\lambda f)(k),\qquad \Ree(\lambda)<0\ .\label{las2}
\end{eqnarray}

3. is a consequence of 1. and 2.. Indeed let
$W,W_1\subset U$ be a compact subsets such that $W\subset \inter(W_1)$. Let $\chi\in C_c^\infty(U)$ be such that
$\chi_{|W_1}\equiv 1$. Then we can write $f=\chi f + (1-\chi) f$,
where $\chi f$ is smooth and $\supp(1-\chi)f\subset\partial X\setminus \inter(W_1)$. We now apply 1. to $\chi f$ and 2. to $(1-\chi)f$ for $kM\in W$.

It remains to prove assertion 2. We imitate the argument of \cite{vandenbanschlichtkrull89}, Thm. 4.8. which proves the assertion for the case $\sigma=1$.

Let $w\in N_K(\aaaa)$ be a representative of the non-trivial element of $W(\gaaa,\aaaa)$.
In the following computation we write the pairing of a distribution
with a smooth function as an integral. 
\begin{eqnarray}
(P^T_\lambda f)(ka)&=&\int_K a(a^{-1}k^{-1}h)^{-(\lambda+\rho)}\gamma(\kappa(a^{-1}k^{-1}h)) T f(h) dh\nonumber\\
&=& \int_{\bar{N}} a(  a^{-1}w\kappa(\bar{n}))^{-(\lambda+\rho)}       \gamma(\kappa( a^{-1}w\kappa(\bar{n})))  Tf(kw\kappa(\bar{n}))
a(\bar{n})^{-2\rho} d\bar{n} \nonumber\\
&=&\int_{\bar{N}} a(a  \bar{n}  a^{-1})^{-(\lambda+\rho)} a^{-(\lambda+\rho)} a(\bar{n})^{\lambda+\rho} \gamma(w) \gamma(\kappa(a\bar{n}a^{-1}))  Tf(kw\kappa(\bar{n}))
a(\bar{n})^{-2\rho} d\bar{n} \nonumber\\
&=&a^{-(\lambda+\rho)}\gamma(w) \int_{\bar{N}} a(\bar{n})^{\lambda-\rho}  
a(a\bar{n}a^{-1})^{-(\lambda+\rho)} \gamma(\kappa(a\bar{n}a^{-1})) T f(kw\kappa(\bar{n})) d\bar{n}\label{tyh}\ .
\end{eqnarray}
For $z\in \R^+$ define $a_z\in A$ by $z=a_z^{-\alpha}$. We consider
the map
$\Phi:(0,\infty)\times \bar{N}\ni (z,\bar{n})\mapsto a_z\bar{n}a_z^{-1}\in \bar{N}$  which according to the decomposition of $\naaa$ into root spaces $\naaa=\naaa_\alpha\oplus\naaa_{2\alpha}$ can also be written as
$$\Phi(z,\exp(X+Y)):= \exp(zX+z^2Y), \quad X\in\naaa_\alpha, Y\in\naaa_{2\alpha}\ .$$
Thus $\Phi$ and hence $ (z,\bar{n})\mapsto a(a_z\bar{n}a_z^{-1})^{-(\lambda+\rho)} \gamma(\kappa(a_z\bar{n}a_z^{-1}))$ extend analytically to  $\R\times\bar{N}$.
Taking the Taylor expansion with respect to $z$ at $z=0$ we obtain  
   $$a(a_z\bar{n}a_z^{-1})^{-(\lambda+\rho)} \gamma(\kappa(a_z\bar{n}a_z^{-1}))=\id + \sum_{n\ge 1}A_n(\bar{n}) z^n\ .$$
Here $A_n:\bar{N}\rightarrow \End(V_\gamma)$ are analytic and the series converges
in the spaces of smooth functions on $\bar{N}$ with values in $\End(V_\gamma)$.

Inserting this expansion into (\ref{tyh}) we obtain
$$(P^T_\lambda f)(ka)=a^{-(\lambda+\rho)} T^w (\hat{J}_\lambda f)(k) + \sum_{n\ge 1} a^{-(\lambda+\rho)-n\alpha} \psi_n(k)\ ,$$
where 
$$\psi_n(k):=\gamma(w)\int_{\bar{N}} A_n(\bar{n})  T a(\bar{n})^{\lambda-\rho }f(kw\kappa(\bar{n})) d\bar{n}\ .$$
Note that $k\mapsto f(kw\kappa(.))$ is a smooth family of distributions
with compact support in $\bar{N}$.
Thus $\psi_n$ is smooth. This finishes the proof of the lemma.
\hB

\begin{kor}\label{lim}
Let $\phi\in C^{-\infty}(\partial X,V(\sigma_\lambda))$ and $f \in C^{\infty}(\partial X,V((\tilde\gamma_{|M})_{-\lambda}))$. If $\Ree(\lambda)>0$, then there exists $\epsilon>0$ such that for $a\to\infty$ we have 
$$\int_K \langle P^T_\lambda \phi(ka),f(k)\rangle\ dk=a^{\lambda-\rho} \langle \phi,{}^t(c_\gamma(\lambda)T)f\rangle+O(a^{\lambda-\rho-\epsilon})\ .$$
There are corresponding formulas for $\Ree(\lambda)=0$.
\end{kor}
\proof
The argument is adapted from \cite{schlichtkrull84}, 5.1. Let $v_i$ be a basis of $V_{\tilde\gamma}$, and write $f(k)=\sum f_i(k) v_i$.
Define $\check{}: C^{\infty}(K)\rightarrow C^\infty(K)$ by $\check h(k):=h(k^{-1})$. Let $*$ be the convolution on the group $K$. Then we have
\begin{eqnarray*} 
\int_K \langle P^T_\lambda \phi(ka),f(k)\rangle\ dk&=&
\sum_i \langle v_i,\int_K P^T_\lambda \phi(ka)f_i(k)\ dk\rangle\\
&=& \sum_i\langle v_i,(\check f_i*P^T_\lambda \phi(.a))(1)\rangle\\
&=& \sum_i\langle v_i, P^T_\lambda (\check f_i*\phi)(a)\rangle\ .
\end{eqnarray*}
Now observe that $\check f_i*\phi\in C^{\infty}(\partial X,V(\sigma_\lambda))$ and
apply Lemma \ref{poisasm}.
\hB

Now we show a variant of Green's formula.
We need nice cut-off functions which exist by the following lemma.
Let $\Delta_X$ be the Laplace-Beltrami operator of $X$. 
\begin{lem}\label{lll}
There exists a cut-off function $\chi\in C_c^\infty(X\cup \Omega)$ such that
\begin{enumerate}
\item $\sum_{g\in\Gamma}g^*\chi = 1$,
\item $\sup_{gK\in X} a_g^{\alpha}\:|d\chi(gK)|<\infty$, 
\item $\sup_{gK\in X} a_g^\alpha\:|\Delta_X\chi(gK)|<\infty$, 
\item $\sup_{gK\in X} |\chi(DgK)|<\infty$, $\forall D\in\cU(\gaaa)$.
\end{enumerate}
We denote the restriction of $\chi$ to $\Omega$ by $\chi_\infty$.
\end{lem}
\proof
Let $W\subset X\cup \Omega$ be a compact subset such that $\bigcup_{g\in\Gamma} gW=X\cup \Omega$. Then we choose a cut-off  function $\psi\in C_c^\infty(X\cup \Omega)$ such that
$\psi_{|W}=1$.  
 We define
$$\chi:=\frac{\psi}{\sum_{g\in\Gamma}g^*\psi}\ .$$
Note that $\chi$ is well-defined since
$\sum_{g\in\Gamma}g^*\psi$ never vanishes on $X\cup\Omega$.
Property 1 is obvious by the definition of $\chi$.

In the following we consider $\chi$ as an element of $C^\infty(\bar{X})$.
Since $G$ acts smoothly on the compact manifold $\bar{X}$
we have $L^\sharp\in C^\infty(\bar{X},T\bar{X})$, where
$L^\sharp$ denotes the fundamental vector field corresponding to $L\in\gaaa$. 
If $Y\in C^\infty(\bar{X},T\bar{X})$, then $Y(C^\infty(\bar{X}))\subset C^\infty(\bar{X})$. Thus the left action of $\gaaa$ and hence of
$\cU(\gaaa)$ preserves $C^\infty(\bar{X})$. This shows assertion 4.

We have $d\chi\in C^\infty(\bar{X},T^*\bar{X})$.
The tensor field in $C^\infty(X,S^2TX)$ which is dual to the Riemannian metric of $X$
has a continuous extension to $\bar{X}$ vanishing of second order
at $\partial X$. We conclude that $|d\chi|$ vanishes
of first order at the boundary of $\bar{X}$. This shows assertion
2.

Since the coefficients of $\Delta_X$ vanish at $\partial X$ of at least first
order assertion 3 follows. 
\hB

Let $\phi\in {}^\Gamma C^{-\infty}(\Lambda,V(\sigma_\lambda,\varphi))$. Then
$res\circ \hat{J}_\lambda(\phi)\in C^{\infty}(B,V_B(\sigma_{-\lambda},\varphi))$ is well-defined
even if $\hat{J}_\mu$ has a pole at $\mu=\lambda$. In the latter case the residue
of $\hat{J}_\mu$ at $\mu=\lambda$ is a differential operator $D_\lambda$ (see Lemma \ref{off})
and $res\circ D_\lambda(\phi)=0$.

\begin{prop}\label{green}
If
$\phi\in {}^\Gamma C^{-\infty}(\Lambda,V(\sigma_\lambda,\varphi))$ and 
$f\in {}^\Gamma C^{-\infty}(\partial X,V(\tilde{\sigma}_\lambda,\tilde\varphi))$, 
then $$\langle res\circ \hat J_\lambda (\phi),res(f)\rangle =0\ .$$
\end{prop}
\proof
At first we need
\begin{lem}\label{ddeenn}
The space 
$${}^\Gamma C^{-\infty}_\Omega(\partial X,V(\sigma_\lambda,\varphi)):=\{f\in {}^\Gamma C^{-\infty}(\partial X,V(\sigma_\lambda,\varphi))\:|\: f_{|\Omega}\in C^\infty(\Omega,V(\sigma_\lambda,\varphi))\}$$ 
is dense in 
${}^\Gamma C^{-\infty}(\partial X,V(\sigma_\lambda,\varphi))$.
\end{lem}
\proof
By Theorem \ref{part1} $ext$ has an at most finite-dimensional singularity at $\lambda$.  Thus there is a finite-dimensional
subspace $W\subset C^\infty(B,V_B(\tilde\sigma_{-\lambda},\tilde\varphi))$ such that 
$$ext_{|W^\perp}:
W^\perp\rightarrow {}^\Gamma C^{-\infty}(\partial X,V(\sigma_\lambda,\varphi))$$ is a well-defined continuous map, where $W^\perp:=\{\phi\in C^{-\infty}(B,V_B(\sigma_\lambda,\varphi))\:|\: \langle\phi,W\rangle=\{0\}\}$.
Since $C^\infty(B,V_B(\sigma_\lambda,\varphi))\subset C^{-\infty}(B,V_B(\sigma_\lambda,\varphi))$
is dense we can choose a complement $\tilde{W}\subset C^\infty(B,V_B(\sigma_\lambda,\varphi))$
such that $C^{-\infty}(B,V_B(\sigma_\lambda,\varphi))=W^\perp\oplus \tilde{W}$.

Let $f\in  {}^\Gamma C^{-\infty}(\partial X,V(\sigma_\lambda,\varphi))$.
Then we can write $res (f)=g=g^\perp\oplus \tilde{g}$,
$g^\perp\in W^\perp$, $\tilde{g}\in \tilde{W}$.
Now $res(f-ext(g^\perp))=\tilde{g}$.
It follows that $f-ext(g^\perp)\in {}^\Gamma C^{-\infty}_\Omega(\partial X,V(\sigma_\lambda,\varphi))$. Let now $g_i\in C^\infty(B,V_B(\sigma_\lambda,\varphi))$ be a sequence
such that $\lim_{i\to\infty} g_i=g$. Then we can decompose
$g_i=g_i^\perp+\tilde{g}_i$. The sections $g_i^\perp$ are smooth
since $g_i$ and $\tilde{g}_i\in\tilde{W}$ are so.
It follows that $ext (g^\perp_i)\in {}^\Gamma C^{-\infty}_\Omega(\partial X,V(\sigma_\lambda,\varphi))$.
By continuity of $ext_{|W^\perp}$
we have $ext(g^\perp)=\lim_{i\to\infty} ext(g^\perp_i)$.
The assertion of the lemma now follows from
$f=f-ext(g^\perp)+\lim_{i\to\infty} ext(g^\perp_i)$.
\hB

We now prove Proposition \ref{green} for the case $\sigma=1$, $\Ree(\lambda)>0$.
We consider the Poisson transforms, both denoted by $P$,
$$ P=P_\lambda\otimes\id: C^{-\infty}(\Lambda,V(1_\lambda,\varphi))\rightarrow C^\infty(X,V(1,\varphi)) $$
and
$$ P=P_\lambda\otimes\id: C^{-\infty}(\Lambda,V(1_\lambda,\tilde\varphi))\rightarrow C^\infty(X,V(1,\tilde\varphi)) \ . $$
Let $D:=\Delta_X \otimes\id-\rho^2+\lambda^2$ be the shifted Laplace operator acting on $C^\infty(X,V(1,\varphi))$ and $C^\infty(X,V(1,\tilde\varphi))$,
respectively. Then we have $D\circ P=0$.
Let $\chi$ be a cut-off function as in
Lemma \ref{lll}.
By $B_R$ we denote the metric $R$-ball centered at the origin of $X$.
The following is an application of Green's formula: 
\begin{eqnarray} 
0&=&\int_{B_R}\chi(x)\left( \langle  D P\phi(x),Pf(x)\rangle- \langle P\phi(x),D Pf(x)  \rangle\right)\: dx\label{limiz}\nonumber\\
&=& \int_{B_R}\left(\langle D \chi(x)  P\phi(x),Pf(x)\rangle- \langle \chi(x)  P\phi(x),D Pf(x)  \rangle - \langle [ D, \chi]  P\phi(x),Pf(x)\rangle\right)\: dx\nonumber                           \\
&=& - \int_{\partial B_R} \left(\langle \nabla_n \chi(y)  P\phi(y),Pf(y)\rangle-\langle \chi(y)  P\phi(y),\nabla_n Pf(y)  \rangle\right)\:dy\nonumber\\
&& - \int_{B_R}\langle [ D, \chi]  P\phi(x),Pf(x)\rangle\: dx\ , 
\end{eqnarray}
where   $n$ is the exterior unit normal vector field at $\partial B_R$.

By Lemma \ref{ddeenn} we can assume that $f_{|\Omega}$ is smooth. Then Lemma
\ref{poisasm}, 2. and 3. combined with
properties  2 and 3 of $\chi$ implies that
$|\langle [ D, \chi]  P\phi,Pf\rangle|$
is integrable over all of $X$. From  Lemma \ref{lll}, property 1, and the $\Gamma$-invariance
of $\langle P\phi,Pf\rangle$ it follows that
$\int_X\langle [ D, \chi]  P\phi,Pf\rangle\:dx=0$.
Taking the limit $R\to\infty$ in (\ref{limiz}) we obtain by Lemma
\ref{poisasm}, 2. and 3.
\begin{eqnarray*}
0&=&  (\lambda+\rho)\int_{\partial X} \chi_\infty(k) \langle (\hat{J}_\lambda \phi )(k), c_1(\lambda)f(k)\rangle\:dk  \\
&&+(\lambda-\rho)\int_{\partial X} \chi_\infty(k) \langle (\hat{J}_\lambda \phi )(k),c_1(\lambda)f(k)\rangle\:dk  \\
&=&2\lambda c_1(\lambda) \int_{\partial X} \chi_\infty(k) \langle  (\hat{J}_\lambda \phi )(k),f(k)\rangle\:dk \\
&=&2\lambda c_1(\lambda)\langle res\circ\hat{J}_\lambda (\phi),res( f)\rangle  \ .
\end{eqnarray*}
This is the assertion of the proposition in our special case since $c_1(\lambda)\not=0$ for $\Ree(\lambda)>0$.

Note that almost the same proof would also work for general $\sigma$  and $\Ree(\lambda)\ge 0$. But for $\Ree(\lambda)<0$ we have to employ tensoring
with finite-dimensional representations, anyway. This method will reduce matters to
the case $\sigma=1$, $\Ree(\lambda)>0$ treated above. As in the proof of Theorem \ref{part1}
we consider the finite-dimensional representation $\pi_{\sigma,\mu}$ of $G$
and the embedding
$$i_{\sigma,\mu}: C^{-\infty}(\partial X,V(\sigma_{\lambda},\varphi)) \hookrightarrow
C^{-\infty}(\partial X,V(1_{\lambda+\mu},\varphi\otimes \pi_{\sigma,\mu}))\ .$$
In addition, the projection onto the lowest weight space of $\pi_{\sigma,\mu}$ induces a surjection
$$ p_{\sigma,\mu}: C^{-\infty}(\partial X,V(1_{-\lambda-\mu},\varphi\otimes \pi_{\sigma,\mu})) \hookrightarrow
C^{-\infty}(\partial X,V(\sigma_{-\lambda},\varphi))\ .$$
Then we have the following identity of meromorphic families of operators (see
\cite{wallach92}, 10.2.6)
\begin{equation}
\hat J_{\sigma,\lambda}\otimes\id_{V_\vp}=p_{\sigma,\mu}\circ \left(\hat J_{1,\lambda+\mu}\otimes\id_{V_\vp\otimes V_{\pi_{\sigma,\mu}}}\right) \circ i_{\sigma,\mu}\ .
\end{equation}
We also consider the induced operators on $B$ denoted by $i_{\sigma,\mu}^B$ and $p_{\sigma,\mu}^B$. Now let $\phi$ and $f$ be as in the statement of the proposition. Choose $\mu\in\aaaa^*$ large enough such that $\pi_{\sigma,\mu}$ exists and $\Ree(\lambda+\mu)>0$. Then we have by the first part of the proof
\begin{eqnarray*}
\langle res\circ \hat J_\lambda (\phi),res(f)\rangle&=&\langle res\circ p_{\sigma,\mu}\circ \hat J_{\lambda+\mu} \circ i_{\sigma,\mu} (\phi),res(f)\rangle\\
&=&\langle p_{\sigma,\mu}^B\circ res\circ  \hat J_{\lambda+\mu}\circ i_{\sigma,\mu} (\phi),res(f)\rangle\\
&=&\langle res\circ  \hat J_{\lambda+\mu}\circ i_{\sigma,\mu} (\phi), i_{\tilde\sigma,\mu}^B\circ res(f)\rangle\\
&=&\langle res\circ  \hat J_{\lambda+\mu}\circ i_{\sigma,\mu} (\phi),  res\circ i_{\tilde\sigma,\mu}(f)\rangle\\
&=&0\ .
\end{eqnarray*}
This proves the proposition.
\hB

Let $\maaa$ denote the Lie algebra of $M$. 
We choose a Cartan subalgebra $\taaa$ of $\maaa$. Then $\taaa\oplus\aaaa=:\haaa$ is a Cartan algebra of $\gaaa$. 
Via the Harish-Chandra isomorphism characters of $\cZ$ are
parametrized by elements of $\haaa_\C^*/W(\gaaa_\C,\haaa_\C)$, where $W(\gaaa_\C,\haaa_\C)$ is the Weyl group of $(\gaaa_\C,\haaa_\C)$. A character $\chi_\nu$, $\nu\in \haaa_\C^*$, is called integral,
if 
\begin{equation}\label{jemi}
2\frac{\langle \nu,\varepsilon\rangle}{\langle\varepsilon,\varepsilon\rangle}\in\Z
\end{equation}
for all roots $\varepsilon$ of $(\gaaa_\C,\haaa_\C)$. 

We further choose a positive root system of $(\maaa_\C,\taaa_C)$.
Let $\rho_m$ denote half of the sum of these positive roots.
For $\sigma\in \hat{M}$ let $\mu_\sigma\in\ii\taaa^*$ be its highest weight. The infinitesimal character of the principal series representation $\pi^{\sigma,\lambda}$ of $G$, $\sigma\in \hat{M}$, $\lambda\in\aca$, is now given
by $\chi_{\mu_\sigma+\rho_m-\lambda}$.
If $X\not= \R H^2$, then for $\sigma\in\hat M$
we define the lattice $$I_\sigma:=\{\lambda\in\aaaa^*\:|\:\chi_{\mu_\sigma+\rho_m-\lambda}\mbox{ is integral }\}\ .$$
If $X=\R H^2$ and $G=SL(2,\R)$, then $M\cong \Z_2$. Let $\pm 1$ denote the trivial (+), resp. non-trivial (-) irreducible representation of $M$. We define 
$$I_{1}:=(\frac{1}{2}+\Z)\alpha,\quad I_{-1}:=\Z\alpha\ .$$
If $G=PSL(2,\R)$, then $M=\{1\}$, and we define $I_{1}:=(\frac{1}{2}+\Z)\alpha$.
 
If $\lambda\not\in I_\sigma$, then the principal series representation $\pi^{\sigma,\lambda}$ is irreducible (see \cite{collingwood85}, 4.3.3). Note that $I_{\sigma^w}=I_\sigma$, so the definition is compatible with 
our previous convention concerning the Weyl-invariance of $\sigma$. 
Let $I_\aaaa\subset\aaaa^*$ be the $\Z$-module generated by the short root $\alpha$,
if $2\alpha$ is a root, or by $\alpha/2$, if not. Then all poles of $\hat J_\lambda$ are located in $I_\aaaa$ (\cite{knappstein71}, Thm. 3 and Prop. 43). Note that for any $\sigma$ we have that $I_\sigma\subset I_\aaaa$ is a sublattice of index $2$.
 
\begin{lem}\label{th43}
Assume that $\sigma\in\hat M$, $\lambda\in \aca$ satisfy one of the following conditions 
\begin{enumerate}
\item $\Ree(\lambda)\ge 0$ and $\sigma=1$.
\item $\Ree(\lambda)\ge 0$, $G=SL(2,\R)$ and $\sigma=-1$. 
\item $\Ree(\lambda)\ge 0$ and $\lambda\not\in I_\sigma$.
\item $\Ree(\lambda)<0$ and $\lambda\not\in I_\aaaa$. 
\end{enumerate}
Let $U\subset \partial X$ be a non-empty open subset,
and  let $\phi\in C^{-\infty}(\partial X, V(\sigma_\lambda))$ be such that $\phi_{|U}= 0$ and
$(\hat{J}_{\sigma,\lambda}\phi)_{|U}= 0$. Then $\phi=0$.
\end{lem}
Before turning to the proof we recall that Lemma \ref{off} implies that  $(\hat{J}_{\sigma,\lambda}\phi)_{|U}$
is well-defined even if $\hat{J}_{\sigma,\lambda}$ has a pole.\\[0.5cm]\noindent
\proof 
We modify an argument given by van den Ban-Schlichtkrull \cite{vandenbanschlichtkrull89}
for the case $\sigma=1$.

Assume that $(\sigma,\lambda)$ satisfies condition 4. Since $\Ree(\lambda)<0$, the operator $\hat J_{\sigma,\lambda}$ is regular and non-vanishing. 
Since $\lambda\not\in I_\aaaa\supset I_\sigma$ the principal series representation $\pi^{\sigma,\lambda}$ is
irreducible. This
implies that $\hat J_{\sigma,\lambda}$ is bijective. Moreover, $\hat J_{\sigma,-\lambda}$ is regular, since $-\lambda\not\in I_\aaaa$. Thus, by the functional equation (\ref{spex}) we can reduce the proof to case 1, 2 or 3 replacing $\phi$ by $\hat{J}_{\sigma,\lambda}(\phi)$. 

We now assume that $\Ree(\lambda)\ge 0$, $\sigma=1$ and $\phi\not\equiv 0$.
Then the Poisson transform $P_\lambda\phi\in C^\infty(X)$ does not vanish and is annihilated
by $D:=\Delta_X-\rho^2+\lambda^2$.
Without loss of generality we can assume that $M\in U$.
Since $P_\lambda\phi$ is real analytic and not identically zero the 
expansion (\ref{epan}) has non-trivial terms.
Let $m$ be the smallest integer such that $\psi_m\not\equiv 0$ near $M$,
where we set $\psi_0:= \hat{J}_{1,\lambda} \phi$.
We have to show that $m=0$.

With respect to the coordinates $k,a$ the operator $D$ has the form
$D=D_0+a^{-2\alpha}R(a,k)$, where $D_0$ is a constant coefficient
operator on $A$ and $R$ is a differential operator with coefficients which remain bounded if $a\to\infty$ (see \cite{helgason84}, Ch. IV, \S 5, (8)).
Moreover, it is known that $D_0$ coincides with the $\bar{N}$-radial part of $D$.

We consider the $\bar{N}$-invariant function $f\in C^{\infty}(X)$
defined by
$f(\bar{n}a):=a^{-(\lambda+\rho+m\alpha)}$.
Since $D$ annihilates the asymptotic expansion (\ref{epan}) we have
$Df=D_0f=0$. On the other hand, $f$ satisfies $(\Delta_X-\rho^2+(\lambda+m\alpha)^2)f=0$. Hence $(\lambda+m\alpha)^2=\lambda^2$.
Since $\Ree(\lambda)\ge 0$ we conclude that $m=0$.

The proof for general $\sigma$ proceeds similarly. Let $\Ree(\lambda)\ge 0$ and
$\lambda\not\in I_\sigma$. In particular, the principal series representation
$\pi^{\sigma,\lambda}$ is irreducible.
We choose $0\not=T\in \Hom_M(V_\sigma,V_\gamma)$ for a suitable $\gamma\in\hat K$. Then $P:=P^T_\lambda$ is injective, and
the range of $P$ can be identified with the kernel of a certain invariant
differential operator $D:C^{\infty}(X,V(\gamma))\rightarrow C^{\infty}(X,V(\gamma^\prime))$ 
for some representation $\gamma^\prime$ of $K$ (see \cite{bunkeolbrich947}, Sec.3).

We now assume $\phi\not\equiv 0$.
Moreover, without loss of generality we can assume that $ M\in U$.
Since $P\phi$ is real analytic, the expansion (\ref{epan}) does not vanish.
Let $m$ be the smallest integer such that $\psi_m\not\equiv 0$ near $M$
(where $\psi_0:=T^w\hat{J}_{\sigma,\lambda} \phi)$.
Again, $D=D_0+a^{-\alpha}R(a,k)$, where $D_0$ is the constant coefficient
operator on $A$ given by the $\bar{N}$-radial part of $D$, and $R$ remains bounded if $a\to\infty$ (see \cite{warner721}, Thm. 9.1.2.4). 

Choose  $k\in K$  near $1$ and $\sigma^\prime\subset \gamma_{|M}$ such that
that there exists an orthogonal projection $S\in \Hom_M(V_\gamma,V_{\sigma^\prime})$ with 
 $S \gamma(k)\psi_m(k) =:v \not= 0$.
Consider the $\bar{N}$-invariant section $f\in C^{\infty}(X,V(\gamma))$
defined by
$$f(\bar{n}a):=a^{-(\lambda+\rho+m\alpha_1)} S^* v\ .$$
Since $D$ annihilates the asymptotic expansion (\ref{epan}), one can check
that $Df=D_0f=0$ and thus $f=P\phi_1$ for some $\bar{N}$-invariant $\phi_1\in C^{-\infty}(\partial X,V(\sigma_\lambda))$.

Now $f=P^S_{\lambda+m\alpha_1}\delta v$, where $\delta v\in C^{-\infty}(\partial X,V(\sigma^\prime_{\lambda+m\alpha_1}))$
is the delta distribution at $1$ with vector part $v$. 
Since $D$ and $P^S_{\lambda+m\alpha_1}$ are $G$-equivariant and $\delta v$ generates the $G$-module $C^{-\infty}(\partial X,V(\sigma^\prime_{\lambda+m\alpha}))$, we obtain
a non-trivial
intertwining operator $I$ from $C^{-\infty}(\partial X,V(\sigma^\prime_{\lambda+m\alpha}))$ 
to the kernel of $D$, hence to $C^{-\infty}(\partial X,V(\sigma_\lambda))$.
This implies 
\begin{equation}\label{hugo}
\chi_{\mu_\sigma+\rho_m-\lambda}
=\chi_{\mu_{\sigma^\prime}+\rho_m-\lambda-m\alpha}\ .
\end{equation}
We conclude that $\chi_{\mu_{\sigma^\prime}+\rho_m-\lambda-m\alpha}$ is not integral and thus $\pi^{\sigma^\prime,\lambda+m\alpha}$ is irreducible. Hence $I$ is an isomorphism. Counting
$K$-types one finds that $\sigma=\sigma^\prime$ or $\sigma=w\sigma^\prime$. It follows
that $|\mu_\sigma+\rho_m|^2 = |\mu_{\sigma^\prime}+\rho_m|^2$, and hence $\lambda^2=(\lambda+m\alpha)^2$. The condition $\Ree(\lambda)\ge 0$ implies that $m=0$. 
Hence $(\hat{J}_{\sigma,\lambda}\phi)_{|U}\not= 0$.
 
In case $G=SL(2,\R)$ and $\sigma=-1$ we argue as above in order to obtain (\ref{hugo}) which is in this case equivalent to $\lambda^2=(\lambda+m\alpha)^2$.
This completes the proof of the lemma.
\hB

The above argument can be extended to cover also some cases
of $\sigma\not=1$, $\Ree(\lambda)\ge 0$ with $\lambda\in I_\sigma$. 
This would lead to stronger vanishing results 
for ${}^\Gamma C^{-\infty}(\Lambda,V(\sigma_\lambda,\varphi))$ below.
But there exist examples of $\sigma\in\hat{M}$
and $\lambda\in\aaaa^*$ with $\Ree(\lambda)\ge 0$ and 
$\lambda\in I_\sigma$, where
the assertion of Lemma \ref{th43} is false.
This is connected with the non-triviality of the spaces $U_\Lambda(\sigma_\lambda,\vp)$ introduced
in Definition \ref{lenin} below.

\begin{kor}\label{ujn}
Under the assumptions of Lemma \ref{th43}
$$res_\Omega\circ \hat{J}_{\sigma,\lambda}:C^{-\infty}(\Lambda,V(\sigma_\lambda))\rightarrow C^{\infty}(\Omega,V(\sigma_{-\lambda}))$$
is injective.
\hB
\end{kor}

For the remainder of this section we assume $X\not=\OO H^2$. The following
corollary now gives the first part of Theorem \ref{disfin}.
\begin{kor}\label{isol}
For fixed $\sigma$ and $\varphi$ the set $\{\lambda\in \aca\:|\:{}^{\Gamma} C^{-\infty}(\Lambda,V(\sigma_\lambda,\varphi))\not=0\}$ is discrete.
\end{kor}
\proof
Assume that $ext: C^{-\infty}(B,V_B(\tilde\sigma_\lambda,\tilde\varphi))\rightarrow
{}^{\Gamma} C^{-\infty}(\partial X,V(\tilde\sigma_\lambda,\tilde\varphi))$
is regular at $\lambda$. Then $res: {}^{\Gamma} C^{-\infty}(\partial X,V(\tilde\sigma_\lambda,\tilde\varphi))\rightarrow C^{-\infty}(B,V_B(\tilde\sigma_\lambda,\tilde\varphi))$ is surjective. It follows from Proposition \ref{green} that $res\circ\hat J_\lambda\left({}^{\Gamma} C^{-\infty}(\Lambda,V(\sigma_\lambda,\varphi))\right)=0$. Assume in
addition that $\lambda\not\in I_\aaaa$. Then ${}^{\Gamma} C^{-\infty}(\Lambda,V(\sigma_\lambda,\varphi))=0$ by Corollary \ref{ujn}. We conclude that ${}^{\Gamma} C^{-\infty}(\Lambda,V(\sigma_\lambda,\varphi))\not=0$ implies that $ext$ is singular at $\lambda$ or $\lambda\in I_\aaaa$.
\hB
  
We denote by $\cO_\lambda C^{-\infty}(\partial X,V(\sigma_.,\varphi))$ and
$\cO_\lambda C^{-\infty}(B,V_B(\sigma_.,\varphi))$ the spaces of germs at $\lambda$ of holomorphic families $\mu\mapsto f_\mu\in C^{-\infty}(\partial X,V(\sigma_\mu,\varphi))$ and $\mu\mapsto f_\mu\in C^{-\infty}(B,V_B(\sigma_\mu,\varphi))$, respectively.
\begin{ddd}\label{horror}
Let 
$$\cO_\lambda^0C^{-\infty}(B,V_B(\sigma_.,\varphi)):=\{f_\mu\in \cO_\lambda C^{-\infty}(B,V_B(\sigma_.,\varphi))\:|\: (\mu-\lambda)ext f_\mu \in \cO_\lambda C^{-\infty}(\partial X,V(\sigma_.,\varphi))\}\ .$$
We define
the space of invariant distributions on the limit set which are generated by the singular parts of $ext$ by
$$
E_\Lambda(\sigma_\lambda,\varphi):=\{\res_{\mu=\lambda}ext(f_\mu)\:| f_\mu\in \cO_\lambda^0C^{-\infty}(B,V_B(\sigma_.,\varphi))\}\ .$$
\end{ddd}

\begin{prop}\label{ghu}
The space $E_\Lambda(\sigma_\lambda,\varphi)$ is finite-dimensional and
\begin{equation}\label{aal}
E_\Lambda(\sigma_\lambda,\varphi)\subset {}^{\Gamma} C^{-\infty}(\Lambda,V(\sigma_\lambda,\varphi))\ .
\end{equation}
Assume that $(\sigma,\lambda)$ satisfies the assumptions of Lemma \ref{th43}. Then equality holds in (\ref{aal}) and $ext$ is regular at $\lambda$ if and only if $\:{}^\Gamma C^{-\infty}(\Lambda,V(\sigma_\lambda,\varphi))=0$.
\end{prop} 
\proof
Since $ext$ has at most finite-dimensional singularities $E_\Lambda(\sigma_\lambda,\varphi)$ is finite-dimensional. (\ref{aal}) follows 
from the meromorphy of $ext$ and the equation $res\circ ext=\id$. It remains to show that under the assumptions of Lemma \ref{th43}
$$\dim E_\Lambda(\sigma_\lambda,\varphi)=\dim {}^{\Gamma} C^{-\infty}(\Lambda,V(\sigma_\lambda,\varphi))\ .$$
The main step in the proof is to show that $\dim E_\Lambda(\sigma_\lambda,\varphi)\ge \dim\coker(res)$ (see Formula (\ref{surprise})
below). This is more or less obvious if $ext$ has a pole of first order at $\lambda$. The general case is more involved. 

Let $L_\lambda$ be the multiplication operator on $\cO_\lambda C^{-\infty}(\partial X,V(\sigma_.,\varphi))$ and $\cO_\lambda C^{-\infty}(B,V_B(\sigma_.,\varphi))$ given by $(L_\lambda f)_\mu:=(\mu-\lambda)f_\mu$. $\cO_\lambda C^{-\infty}(\partial X,V(\sigma_.,\varphi))$ becomes a
$\Gamma$-module by $(\pi(g)f)_\mu:=(\pi^{\sigma,\mu}\otimes\varphi)(g)f_\mu$, $g\in\Gamma$. The $\Gamma$-action commutes with $L_\lambda$. The restriction
map induces a morphism of $\C[L_\lambda]$-modules
$${}_{\cO_\lambda}res: \:{}^\Gamma\cO_\lambda C^{-\infty}(\partial X,V(\sigma_.,\varphi))
\rightarrow
\cO_\lambda C^{-\infty}(B,V_B(\sigma_.,\varphi))\ .$$
By Corollary \ref{isol} and Lemma \ref{mainkor} we have for $f\in
{}^\Gamma\cO_\lambda C^{-\infty}(\partial X,V(\sigma_.,\varphi))$
\begin{equation}\label{babel}
ext\circ{}_{\cO_\lambda}res (f)=f\ .
\end{equation}
Let $k$ be the order of the pole of $ext$ at $\lambda$. Set
$$ \cO_{\lambda,k}C^{-\infty}(\partial X,V(\sigma_.,\varphi)):=
\coker\left(L_\lambda^k: \cO_\lambda C^{-\infty}(\partial X,V(\sigma_.,\varphi))
\rightarrow \cO_\lambda C^{-\infty}(\partial X,V(\sigma_.,\varphi))\right)\ .$$
Then we consider the singular part of $ext$ as an operator between $\C[L_\lambda]$-modules
$$ ext^{<0}: \coker{}_{\cO_\lambda}res\rightarrow {}^\Gamma \cO_{\lambda,k}C^{-\infty}(\partial X,V(\sigma_.,\varphi))$$
given by
$$  ext^{<0}([f]):= ext(L_\lambda^k f)\ \mod\  \im L^k_\lambda\ ,\quad f\in \cO_\lambda C^{-\infty}(B,V_B(\sigma_.,\varphi))\ ,$$
which is well-defined by (\ref{babel}). Assume that $ext^{<0}([f])=0$. Then $ext(L_\lambda^kf)=L_\lambda^kg$ for some
$g\in \cO_\lambda C^{-\infty}(\partial X,V(\sigma_.,\varphi))$. It follows that
$f={}_{\cO_\lambda}res(g)$. We conclude that $ext^{<0}$ is injective. In particular, since $ext^{<0}$ is a finite-dimensional operator, the space $\coker{}_{\cO_\lambda}res$ is finite-dimensional. 
The map
$$ E_\Lambda(\sigma_\lambda,\varphi)\ni \res_{\mu=\lambda}ext(f_\mu)\mapsto
ext^{<0}([f])\in {}^\Gamma \cO_{\lambda,k}C^{-\infty}(\partial X,V(\sigma_.,\varphi)),\ 
f\in \cO_\lambda^0C^{-\infty}(B,V_B(\sigma_.,\varphi)),$$
is well-defined and identifies $E_\Lambda(\sigma_\lambda,\varphi)$ with
$\ker\left(L_\lambda:\im\: ext^{<0}\rightarrow\im\: ext^{<0}\right)$.
We also consider the usual point-wise restriction map
$$ res: {}^\Gamma C^{-\infty}(\partial X,V(\sigma_\lambda,\varphi))
\rightarrow
C^{-\infty}(B,V_B(\sigma_\lambda,\varphi))$$
and the surjection
$$\coker\left(L_\lambda:\coker{}_{\cO_\lambda}res\rightarrow
\coker{}_{\cO_\lambda}res\right)\rightarrow\coker(res)$$
induced by the point evaluation at $\lambda$.
Summarizing the above discussion we obtain
\begin{eqnarray}\label{surprise}
d:=\dim E_\Lambda(\sigma_\lambda,\varphi)&=&\dim 
\ker\left(L_\lambda:\im\: ext^{<0}\rightarrow\im\: ext^{<0}\right)\nonumber\\
&=&\dim\coker\left(L_\lambda:\im\: ext^{<0}\rightarrow\im\: ext^{<0}\right)\nonumber\\
&=&\dim\coker\left(L_\lambda:\coker{}_{\cO_\lambda}res\rightarrow\coker{}_{\cO_\lambda}res\right)\nonumber\\
&\ge&\dim\coker(res)\ .
\end{eqnarray}

Set $\tilde d:= \dim E_\Lambda(\tilde\sigma_\lambda,\tilde\varphi)$. 
If $(\sigma,\lambda)$ satisfies the assumptions of Lemma \ref{th43}, then so does 
$(\tilde\sigma,\lambda)$. Thus Proposition \ref{green} combined with Corollary \ref{ujn} implies 
$$\dim \coker(res)
\ge\dim {}^{\Gamma} C^{-\infty}(\Lambda,V(\tilde\sigma_\lambda,\tilde\varphi))\ .$$ 
It
eventually follows from (\ref{aal}) that
$$d\ge\dim\:\coker(res)
\ge \dim {}^{\Gamma} C^{-\infty}(\Lambda,V(\tilde\sigma_\lambda,\tilde\varphi))\ge\tilde d\ .$$
Changing the roles of $(\sigma,\varphi)$ and $(\tilde\sigma,\tilde\varphi)$
we obtain $d=\tilde d=\dim {}^{\Gamma} C^{-\infty}(\Lambda,V(\sigma_\lambda,\varphi))$.
This finishes the proof of the proposition.
\hB

\begin{kor}\label{speck}
For every datum $(\sigma_\lambda,\varphi)$ the space ${}^{\Gamma} C^{-\infty}(\Lambda,V(\sigma_\lambda,\varphi))$ is finite-dimensional.
\end{kor}
\proof
Choose $\mu\in\aaaa^*$ large enough such that
the finite-dimensional representation $\pi_{\sigma,\mu}$ exists and $\Ree(\lambda+\mu)\ge 0$. Recall the definition (\ref{romski}) of the embedding $i_{\sigma,\mu}$. Since $(1,\lambda+\mu)$ satisfies the assumptions of Lemma \ref{th43} we obtain 
by Proposition \ref{ghu} that
$$ i_{\sigma,\mu}\left({}^{\Gamma} C^{-\infty}(\Lambda,V(\sigma_\lambda,\varphi))\right)\subset {}^{\Gamma} C^{-\infty}(\Lambda,V(1_{\lambda+\mu},\varphi\otimes\pi_{\sigma,\mu}))=E_\Lambda(1_{\lambda+\mu},\varphi\otimes\pi_{\sigma,\mu})\ .$$
The corollary now follows since the space on the right hand side is finite-dimensional.
\hB

Thus we have completed the proof of Theorem \ref{disfin}. It is perhaps worth noting that the finite-dimensionality of ${}^{\Gamma} C^{-\infty}(\Lambda,V(\sigma_\lambda,\varphi))$ can also be proved without refering
to the meromorphy of $ext$. In fact, using the asymptotic expansion (\ref{epan}) it
is not difficult to show that there exists a number $k\in\nat_0$ (depending on $\lambda$) and a continuous embedding 
$$i: {}^{\Gamma} C^{-\infty}(\Lambda,V(\sigma_\lambda,\varphi))\hookrightarrow 
\left(C^{k}(\partial X,V(\tilde\sigma_{-\lambda},\tilde\varphi))
\right)^\prime\ .$$ 
Furthermore, the Banach space topology on ${}^{\Gamma} C^{-\infty}(\Lambda,V(\sigma_\lambda,\varphi))$ induced by $i$ coincides with the topology
as a closed subspace of the Montel space $C^{-\infty}(\partial X,V(\sigma_\lambda,\varphi))$. It follows that the unit ball in the Banach space $ {}^{\Gamma} C^{-\infty}(\Lambda,V(\sigma_\lambda,\varphi))$ is compact.  
This implies Corollary \ref{speck}. However, arguments of this type
seem not to be sufficient in order to establish results like Corollary \ref{isol}.

\section{Consequences of unitarity}\label{uni}

From now on we assume that the finite-dimensional representation $(\vp,V_\vp)$ of $\Gamma$ is unitary. Then $\delta_\vp=0$. If $X=\OO H^2$, then in addition 
we assume that $\delta_\Gamma<0$.
Hence in any case $ext$ is meromorphic on a half-plane 
$\{\lambda\in\aca\:|\:\Ree(\lambda)>-\epsilon\}$ for some $\epsilon>0$. The aim of this section is to work out the consequences of the unitarity of $\vp$ for the singularities of $ext$ in that region.
In particular, for fixed $\sigma$ and $\vp$  we show that the set of parameters $\lambda$ with  $\Ree(\lambda)\ge 0$ and allowing non-trivial invariant distributions
with support on the limit set is finite and real. This result will be the main ingredient in the proof of 
the finiteness of the discrete spectrum of $\cZ$ on
$L^2(Y,V_Y(\gamma,\vp))$, where $V_Y(\gamma,\vp)$ is the hermitian vector 
bundle over $Y$ defined by $V_Y(\gamma,\vp):=\Gamma\backslash V(\gamma,\vp)$.

In contrast to bilinear pairings $\langle .,.\rangle$ sesquilinear pairings
will always be written as $(.,.)$. By convention these pairings are $\C$-linear
with respect to the first variable. 

Since $\sigma$ and $\vp$ are unitary representations,
for $\lambda\in \imath\aaaa^*$ we have a positive definite conjugate linear
pairing $(V_{\sigma_\lambda}\otimes V_\vp)\otimes (V_{\sigma_\lambda}\otimes V_\vp) \rightarrow V_{1_{-\rho}}$. Since $V_B(1_{-\rho})\cong \Lambda^{max}T^*B$ integration
gives
a natural $L^2$-scalar product on $C^\infty(B,V_B(\sigma_\lambda,\varphi))$.
Let $L^2(B,V_B(\sigma_\lambda,\varphi))$ be associated Hilbert space.
Using Lemma \ref{lok} we see that the adjoint $S^*_\lambda$
with respect to this Hilbert space structure is just $S_{\bar\lambda}=S_{-\lambda}$.

\begin{lem}\label{unitary}
If $\Ree(\lambda)=0$, then $S_\lambda$ is regular and unitary.
\end{lem}
\proof
Let $\lambda$ be imaginary such that $S_{\pm\lambda}$ are regular. Let  
$f\in C^\infty(B,V_B(\sigma_\lambda,\vp))$. Then by the functional equation (\ref{forme})
$$\|S_\lambda f\|_{L^2(B,V_B(\sigma_\lambda,\varphi))}^2=( S_{-\lambda}\circ S_\lambda f, f )_{L^2(B,V_B(\sigma_\lambda,\varphi))} = \|f\|^2_{L^2(B,V_B(\sigma_\lambda,\varphi))}\ .$$
By meromorphy of $S_\lambda$ this equation implies that $S_\lambda$ is regular
and unitary on $\imath\aaaa^*$.
\hB

\begin{lem}\label{extregatim}
If $\Ree(\lambda)= 0$ and $\lambda\not=0$, then 
$$ext:C^{-\infty}(B,V_B(\sigma_\lambda,\varphi))\rightarrow {}^\Gamma C^{-\infty}(\partial X,V(\sigma_\lambda,\varphi))$$
is regular and ${}^\Gamma C^{-\infty}(\Lambda,V(\sigma_\lambda,\varphi))=0$.
\end{lem}
\proof
Note that $(\sigma,\lambda)$ satisfies the assumptions of Lemma \ref{th43}.
Assume that $ext$ is singular at a non-zero imaginary $\lambda$. The leading singular part of $ext$ maps
to distributions which are supported on the limit set $\Lambda$.
Thus by Corollary
\ref{ujn} the scattering matrix 
$\hat S_\lambda=res\circ \hat J_\lambda\circ ext$ is singular at $\lambda$, too.
Since $\hat S_\lambda= S_\lambda c_\sigma(-\lambda)$, and $c_\sigma$ is regular
on $\ii\aaaa^*\setminus\{0\}$, this contradicts Lemma \ref{unitary}.
Thus $ext$ is regular at $\lambda$. It follows by Proposition \ref{ghu} that
${}^\Gamma C^{-\infty}(\Lambda,V(\sigma_\lambda,\varphi))=0$.
\hB

\begin{prop}\label{firtu}
\mbox{}
\begin{enumerate}
\item 
If $\Ree(\lambda)>0$, then for any finite-dimensional representation $\gamma$ of $K$ and $T\in\Hom_M(V_\sigma,V_\gamma)$ the Poisson transform $P^T_\lambda$ maps ${}^\Gamma C^{-\infty}(\Lambda,V(\sigma_\lambda,\varphi))$ to $L^2(Y,V_Y(\gamma,\vp))$.\\ If ${}^\Gamma C^{-\infty}(\Lambda,V(\sigma_\lambda,\varphi))\not=0$, then $\lambda$ is real. 
\item 
If $ext:C^{-\infty}(B,V_B(\sigma_\mu, \varphi))\rightarrow {}^\Gamma C^{-\infty}(\partial X,V(\sigma_\mu,\varphi))$ is singular at $\mu=\lambda$, $\Ree(\lambda)>0$, then $\lambda$ is real and the order of the pole is $1$.
\item 
If $p_\sigma(0)\not=0$, i.e., the intertwining operator $\hat J_0$ is regular
(see (\ref{spex})),\\ then ${}^\Gamma C^{-\infty}(\Lambda,V(\sigma_\lambda,\varphi))=0$ for $\Ree(\lambda)>0$, and
$ext$ is regular in that region.
\end{enumerate} 
\end{prop}
\proof
Let $\gamma$ be a finite-dimensional representation of $K$ and $T\in\Hom_M(V_\sigma,V_\gamma)$.
Let $\phi \in {}^\Gamma C^{-\infty}(\Lambda,V(\sigma_\lambda,\vp))$.
By Lemma \ref{poisasm} for any compact subset $F\subset \Omega$ there exists 
a constant $C$ such that for $a\in A_+$, $k\in FM$
\begin{equation}\label{new} 
|P_\lambda^T\phi(ka)|\le C a^{-\lambda-\rho} \ .
\end{equation}
In particular, $P_\lambda^T\phi\in L^2(Y, V_Y(\gamma,\varphi))$.
This shows the first part of assertion 1.

Now let $\gamma\in\hat K$ be a minimal $K$-type of $\pi^{\sigma,\lambda}$ and $T$ be injective. 
Then $P_\lambda^T$ is injective. There is a real constant $c(\sigma)$ 
(see \cite{bunkeolbrich955})
such that $(-\Omega_G+c(\sigma)+\lambda^2)\circ P_\lambda^T =0$.

If $\phi \in {}^\Gamma C^{-\infty}(\Lambda,V(\sigma_\lambda,\vp))$,
then $P_\lambda^T\phi$ is a square-integrable eigenfunction 
of $-\Omega_G+c(\sigma)+\lambda^2$ on $Y$. 
Since $Y$ is complete $\Omega_G$ is essentially selfadjoint on the domain $C_c^\infty(Y,V_Y(\gamma,\varphi))$.
Its selfadjoint closure
has the domain of definition
$\{f\in L^2(Y,V_Y(\gamma,\varphi))\:|\: \Omega_Gf\in L^2(Y,V_Y(\gamma,\varphi))\}$.  In particular, $\Omega_G$ can not have non-trivial eigenvectors in $L^2(Y,V_Y(\gamma,\varphi))$
to eigenvalues with non-trivial imaginary part. Since $\Ree(\lambda)>0$, $\lambda$ has to be real. This proves the second part of assertion 1.

Assume in addition that $p_\sigma(0)\not=0$. Then
there is a representation $\gamma^s$ of $K$, a $K$-equivariant embedding 
$i:\gamma\hookrightarrow \gamma^s$ and a locally invariant Dirac operator
$D(\sigma)$ acting on $C_c^\infty(Y,V_Y(\gamma^s,\varphi))$ such that \begin{equation}\label{huhu8}D(\sigma)^2=-\Omega_G+c(\sigma)\ .\end{equation}
 (If $X$ is e.g. an odd-dimensional hyperbolic 
space , then $\sigma=\sigma^\prime\oplus\sigma^{\prime w}$ for some  non Weyl-invariant
$\sigma^{\prime}\in \hat{M}$ and $\gamma_s$ is constructed in \cite{bunkeolbrich955},
see  pp. 28-29 for (\ref{huhu8}).) Since $D(\sigma)$ is selfadjoint
$-\Omega_G+c(\sigma)$ is non-negative. It follows that $\Ree(\lambda)=0$ which contradicts our assumption. 
Alternatively, one
could use the fact 
(\cite{baldonisilva81}, Thm. 6.1) that in case
$p_\sigma(0)\not=0$ for $\Ree(\lambda)>0$ the principal series 
$\pi^{\sigma,\lambda}$ does not contain unitary subrepresentations.
Compare Proposition \ref{nonte} below. This proves assertion 3.

Now let $f_\mu\in C^{-\infty}(B,V_B(\sigma_\mu,\varphi))$, $\mu\in \aca$,
be a holomorphic family such that $ext(f_\mu)$ has a pole of order
$k\ge 1$ at $\mu=\lambda$, $\Ree(\lambda)>0$. Let $0\not=\phi \in {}^\Gamma C^{-\infty}(\partial X,V(\sigma_\lambda,\varphi))$
be the leading singular part of $ext(f_\mu)$ at $\mu=\lambda$. In particular, $\lambda$ is real and $P_\lambda^T\phi\in L^2(Y, V_Y(\gamma,\varphi))$ by assertion 1.

Since the leading singular part of $ext$ is a finite-dimensional operator we
may assume that $f_\mu\in C^\infty(B,V_B(\sigma_\mu,\varphi))$ without changing
$\phi$.
Then by Lemma \ref{poisasm} for any compact subset $F\subset \Omega$ and $\mu_0>0$ there 
exist constants $C_1,C_2,C_3$ such that for $a\in A_+$, $k\in FM$, $\mu_0\le\mu<\lambda$,
\begin{eqnarray}\label{york}
|(\mu-\lambda)^kP_\mu^T ext(f_\mu)(ka)|&\le& C_1 |\mu-\lambda|^k a^{\mu-\rho}+C_2 a^{-\mu-\rho}\\
&\le& C_3(1+\log a)^{-k} a^{\lambda-\rho}\ .\nonumber
\end{eqnarray}
There is a constant $C_4$ such that for $\mu_0\le\mu<\lambda$
$$ |\mu-\lambda|^k|P_\mu^T ext(f_\mu)(x)|\:  |P_{\lambda}^T\phi(x)| \le C_4 (1+\log(a_x))^{-k}a_x^{-2\rho},\quad x\in F\ .$$
We now assume that $k\ge 2$. Then 
$$\int_{FK} (1+\log(a_x))^{-k}a_x^{-2\rho} dx<\infty\ ,$$
and we obtain by Lebesgue's
Theorem of dominated convergence  
\begin{eqnarray*}
\|P_\lambda^T \phi \|^2_{L^2(Y, V_Y(\gamma,\varphi))}&=& (\lim_{\mu\to\lambda \atop \mu<\lambda} (\mu-\lambda)^{k}  P_\mu^T ext(f_\mu), P_{\lambda}^T\phi )_{L^2(Y, V_Y(\gamma,\varphi))}\\ 
&=& \lim_{\mu\to\lambda \atop \mu<\lambda} (\mu-\lambda)^{k}  ( P_\mu^T ext(f_\mu), P_{\lambda}^T\phi )\ .
\end{eqnarray*}
On the other hand the estimates (\ref{new}) and (\ref{york}) allow  for partial integration,
and we obtain for $\mu<\lambda$
\begin{eqnarray*}
( P_\mu^T ext(f_\mu), P_{\lambda}^T\phi )&=&
\frac{1}{\lambda^2-\mu^2}( (-\Omega_G+c(\sigma)+\lambda^2)
P_\mu^T ext(f_\mu), P_{\lambda}^T\phi )\\
&=&
\frac{1}{\lambda^2-\mu^2}( 
P_\mu^T ext(f_\mu),(-\Omega_G+c(\sigma)+\lambda^2) P_{\lambda}^T\phi )\\
&=&0\ .
\end{eqnarray*}
Hence $\|P_\lambda^T \phi \|_{L^2(Y,V_Y(\gamma,\varphi))}=0$. Since $P_\lambda^T$
is injective this contradicts $\phi\not=0$. We conclude that $k=1$.
This proves assertion 2.
\hB

\begin{prop}\label{firtu0}
\mbox{}
\begin{enumerate}
\item 
The order of the pole of $ext$ at $0$ is at most $1$. 
\item
If $p_\sigma(0)\not=0$, then $ext$ is regular at $0$.
\item
If $p_\sigma(0)=0$, then the residue $\res_{\mu=0}ext$ 
vanishes on the $1$-eigenspace of the
involution $S_0$ and identifies the $-1$-eigenspace
with $E_\Lambda(\sigma_0,\vp)$. 
\end{enumerate}
\end{prop}
\proof 
The proof is analogous to the proof of Proposition \ref{firtu}.
Let $\gamma$ be the representation of $K$ on the sum of the $K$-isotypic
components of $C^{\infty}(\partial X,V(\sigma_0))$ that correspond to all minimal $K$-types of $C^{\infty}(\partial X,V(\sigma_0))$. ($\gamma$ is irreducible iff $p_\sigma(0)=0$.) Then one can find $T\in\Hom_M(V_\sigma,V_\gamma)$ such that the Poisson transform $$P_\mu:=P^T_\mu\otimes\id: C^{-\infty}(\partial X,V(\sigma_\mu,\varphi))\rightarrow C^\infty(X,V(\gamma,\varphi))$$
is injective for $\mu$ in a neighbourhood of $0$. 

Let $\mu\mapsto f_\mu\in C^\infty(B,V_B(\sigma_\mu,\varphi))$ be a holomorphic family defined on such
a neighbourhood. Assume that $ext f_\mu$ has a pole of order $k$ at $0$. We want to study the leading singular part $\phi\in {}^\Gamma C^{-\infty}(\partial X,V(\sigma_0,\varphi))$ of $ext f_\mu$
at $\mu=0$ via the leading singular part $P_0\phi$ of $P_\mu(ext f_\mu)$. For $\mu$ on a sufficiently small pointed disc $D$ around $0$ we obtain by Lemma \ref{poisasm}, 3.,
\begin{equation}\label{wurm}
(P_\mu (extf_\mu))(ka) = a^{\mu-\rho} c_\gamma(\mu)T  f_\mu(k)+ a^{-\mu-\rho}  T^w(\hat{S}_\mu f_\mu)(k) + O(a^{-\rho-\epsilon})
\end{equation}
for $a\to\infty$ uniformly in $\mu$ and $kM$ varying in compact subsets $D_0\subset D$, $F\subset\Omega$.

If $p_\sigma(0)\not=0$, then $c_\gamma(\mu)$ and $\hat{S}_\mu$ are regular 
at $\mu=0$. We multiply (\ref{wurm}) by $\mu^k$ and apply Cauchy's integral formula in order to conclude that for any compact subset $F\subset \Omega$ there exists 
a constant $C$ such that for $a\in A_+$, $k\in FM$
\begin{equation}\label{fortsatz}
|P_0\phi(ka)|\le C a^{-\rho-\epsilon} \ .
\end{equation}
If $p_\sigma(0)=0$,
then  $c_\gamma(\mu)$ and $\hat{S}_\mu$
have only first order poles at $\mu=0$.
If $k\ge 2$, then we can argue as above in order to obtain (\ref{fortsatz}).

In particular, $P_0\phi\in L^2(Y, V_Y(\gamma,\varphi))$. 
By (\ref{wurm}) and (\ref{fortsatz}) the pairing
$(P_\mu ext(f_\mu), P_{0}\phi)$ is defined for $0<|\mu|<\epsilon$. As
in the proof of Proposition \ref{firtu} we show by partial integration that
this pairing  vanishes and that $\|P_0 \phi \|_{L^2(Y,V_Y(\gamma,\varphi))}^2=0$. Hence $\phi=0$ unless $k=1$ and $p_\sigma(0)=0$.

It remains to show the last assertion of the proposition. Thus assume $k=1$ and $p_\sigma(0)=0$. Applying the residue theorem  and (\ref{mao}) to (\ref{wurm}) we obtain
$$(P_0 \phi)(ka) = a^{-\rho} (\res_{\mu=0}c_\sigma(\mu))T^w  (f_0(k)- {S}_0 f_0(k)) + O(a^{-\rho-\epsilon})\ .
$$
Thus the leading asymptotic coefficient of $P_0\phi$ vanishes iff $S_0 f_0 = f_0$.  
In this case $P_0\phi$ satisfies (\ref{fortsatz}) which implies that
$\phi=0$ as above. We conclude that $\ker\res_{\mu=0} ext = \ker (S_0-\id)$.
It follows from the definition of $E_\Lambda(\sigma_0,\vp)$ that
$$\res_{\mu=0} ext :\ker(S_0+\id)\rightarrow E_\Lambda(\sigma_0,\vp)$$
is an isomorphism.
The proof of the proposition is now complete.
\hB

For $0\not \in I_\sigma$ Proposition \ref{firtu0} could also have been proved by a refinement
of the proof of Lemma \ref{extregatim}.

Recall Definition \ref{horror} of $E_\Lambda(\sigma_\lambda,\varphi)$.
For $\Ree(\lambda)\ge 0$ we have just proved that this space coincides with the image of the residue of $ext$ at $\lambda$. 
\begin{ddd}\label{lenin}
For all $\lambda\in\aca$ we define the space
of "stable" invariant distributions supported on the limit set by
$$
U_\Lambda(\sigma_\lambda,\varphi):=\{\phi\in {}^\Gamma C^{-\infty}(\Lambda,V(\sigma_\lambda,\vp))\:|
\:res\circ\hat J_{\sigma,\lambda}(\phi)=0\}\ .$$
\end{ddd}

Corollary \ref{ujn} implies that if $\Ree(\lambda)\ge 0$ and $U_\Lambda(\sigma_\lambda,\varphi)$ is non-trivial, then $\lambda\in I_\sigma$. We can now refine Proposition \ref{ghu}.

\begin{prop}\label{engels}
For $\Ree(\lambda)\ge 0$ we have
$$ {}^\Gamma C^{-\infty}(\Lambda,V(\sigma_\lambda,\vp))=
E_\Lambda(\sigma_\lambda,\varphi)\oplus U_\Lambda(\sigma_\lambda,\varphi)\ .$$
\end{prop}
\proof
We assume that
$\lambda\in\aaaa^*$. Otherwise there is nothing to show. Let $\phi\in E_\Lambda(\sigma_\lambda,\varphi)\cap U_\Lambda(\sigma_\lambda,\varphi)$. Let $P$ be an injective Poisson transform as in the proofs of
Propositions \ref{firtu} and \ref{firtu0}. Using Lemma \ref{poisasm}, 1. in case $\lambda=0$ and 2. in case $\lambda>0$, we see that 
there exists 
a constant $\epsilon>0$ such that for any compact $F\subset\Omega$ there exists a constant $C$ such that for $a\in A_+$, $k\in FM$
$$
|P\phi(ka)|\le C a^{-\lambda-\rho-\epsilon}\ .$$
As in the proof of Proposition \ref{firtu0} this implies $\phi=0$. Thus  $res\circ\hat J_\lambda$ is injective on $E_\Lambda(\sigma_\lambda,\varphi)$. 
It remains to show that
\begin{equation}\label{trotzki}
\dim  E_\Lambda(\sigma_\lambda,\varphi)
=\dim res\circ\hat J_\lambda\left({}^\Gamma C^{-\infty}(\Lambda,V(\sigma_\lambda,\vp))\right)
\end{equation}
By (\ref{surprise}) we have $\dim  E_\Lambda(\sigma_\lambda,\varphi)\ge\dim\coker(res)$. On the other hand Proposition
\ref{green} implies that 
$$\dim\coker(res)\ge\dim res\circ\hat J_\lambda\left({}^{\Gamma} C^{-\infty}(\Lambda,V(\tilde\sigma_\lambda,\tilde\varphi))\right)\ .$$
Since $\lambda$ is real and $\vp$ is unitary the space on the
right hand side is conjugate linear isomorphic to $res\circ\hat J_\lambda\left({}^\Gamma C^{-\infty}(\Lambda,V(\sigma_\lambda,\vp))\right)$. 
This implies (\ref{trotzki}).
\hB

\begin{ddd}\label{kabel}    
Define 
\begin{eqnarray*}
PS(\sigma,\vp)&:=& \{\lambda\in\aca\:|\:\Ree(\lambda)\ge 0,\  {}^{\Gamma} C^{-\infty}(\Lambda,V(\sigma_\lambda,\varphi))\not=0\}\ ,\\
PS_\res(\sigma,\vp)&:=& \{\lambda\in\aca\:|\:\Ree(\lambda)\ge 0,\  E_\Lambda(\sigma_\lambda,\varphi)\not=0\}\ ,\\
PS_U(\sigma,\vp)&:=& \{\lambda\in\aca\:|\:\Ree(\lambda)\ge 0,\  U_\Lambda(\sigma_\lambda,\varphi)\not=0\}.
\end{eqnarray*}
\end{ddd}

\begin{prop}\label{upperbound}
$PS(\sigma,\vp)=PS_\res(\sigma,\vp)\cup PS_U(\sigma,\vp)$ is a finite subset of the interval $[0,\delta_\Gamma]\subset\aaaa^*$, and we have $PS_U(\sigma,\vp)\subset I_\sigma$.
The space
$$\bigoplus_{\lambda\in PS(\sigma,\vp)} {}^\Gamma C^{-\infty}(\Lambda,V(\sigma_\lambda,\varphi))$$
is finite-dimensional. If $\lambda\in PS(\sigma,\vp)\setminus I_\sigma$, then ${}^\Gamma C^{-\infty}(\Lambda,V(\sigma_\lambda,\varphi))$ is spanned by the
residue of $ext$ at $\lambda$.

If $p_\sigma(0)\not=0$, then $PS(\sigma,\vp)=PS_U(\sigma,\vp)\subset\{0\}$.

If $\delta_\Gamma>0$, then $\delta_\Gamma\in PS_\res(1,1)$ and $\dim{}^\Gamma C^{-\infty}(\Lambda,V(1_{\delta_\Gamma}))=1$.
\end{prop}
\proof
Combine the results of the present section with Theorem \ref{disfin} and Theorem \ref{pokm}. The last assertion follows from the construction of the
Patterson-Sullivan measure (\cite{patterson762}, \cite{sullivan79}, \cite{corlette90}) and the uniqueness of the eigen{\em function} $f\in L^2(Y)$ 
corresponding to the smallest eigenvalue of $\Delta$. 
\hB

\begin{kor}\label{treu}
If $X$ is an odd-dimensional hyperbolic space and $\sigma$ is a faithful representation of $M=Spin(n-1)$, then
$PS(\sigma,\vp)$ is empty.
\end{kor}
\proof
The condition on $\sigma$ ensures that $0\not\in I_\sigma$ and that 
$p_\sigma(0)\not=0$ (\cite{knappstein71}, Ch.12).
\hB 

\section{Abstract harmonic analysis on $\Gamma\backslash G$}\label{abs}

Let $(\vp,V_\vp)$ be a unitary representation of $\Gamma$. We consider the Hilbert space
$$ L^2(\Gamma\backslash G,\vp):=\{f:G\rightarrow V_\vp\:|\: f(gx)=\vp(g)f(x)\ \forall g\in\Gamma,\: x\in G,\:\ 
\int_{\Gamma\backslash G} |f(x)|^2\:dx<\infty\}\ .$$
The action of $G$ on $\Gamma\backslash G$ by right translations induces a unitary representation of $G$ on $L^2(\Gamma\backslash G,\vp)$. The abstract Plancherel theorem (see e.g. \cite{wallach92}, Thm. 14.10.5) provides a direct integral decomposition of
this representation into irreducibles
\begin{equation}\label{pl1} 
L^2(\Gamma\backslash G,\vp)\stackrel{\alpha}{\cong} \int_{\hat G}^\oplus M_\pi\hat\otimes V_\pi \: d\kappa(\pi)\ .
\end{equation}
Here $\hat G$ denotes the unitary dual of $G$, $V_\pi$ carries an irreducible unitary representation belonging to the class $\pi\in\hat G$, the Hilbert space $M_\pi$ is called
the multiplicity space of $\pi$ in $L^2(\Gamma\backslash G,\vp)$, and $\kappa$ is a
Borel measure on $\hat G$ called the Plancherel measure. The isomorphism (\ref{pl1}) is a unitary equivalence of representations, where $G$ acts on the
right hand side by $\id_{M_\pi}\otimes\pi$. Throughout the paper we shall not
distinguish between an element $\pi\in\hat G$ and a particular representative
$(\pi,V_\pi)$. Note that only the measure class of $\kappa$ is uniquely defined.
Strictly speaking, the spaces $M_\pi$ are defined on a complement of a set of
measure zero, only. 

If $\Gamma$ is the trivial group, then $M_\pi\cong V_\pi^\prime$, and $M_\pi\hat\otimes V_\pi$ can be identified with the space
${\cal L}^2(V_\pi)$ of Hilbert-Schmidt operators on $V_\pi$ which has a canonical scalar product. Then a canonical choice for $\alpha$ is  $\alpha(f)(\pi):=\pi(f)\in {\cal L}^2(V_\pi)$, $f\in  L^2(G)$, $\pi\in\hat G$.  
This determines the Plancherel measure.
In contrast, for non-trivial
$\Gamma$ there is no obvious normalization of the unitary equivalence $\alpha$,
of the Hilbert space structure on $M_\pi\hat\otimes V_\pi$ and of the Plancherel measure.   

The goal of this paper can now be phrased as follows: Make the equivalence (\ref{pl1}) and all its ingredients as explicit as possible. In the present section we specify this task. In particular, following the approach of Bernstein \cite{bernstein88} we identify $M_\pi$ with a subspace
of the "tempered" functionals on the space of smooth vectors $V_{\pi,\infty}$  
of $V_\pi$. This provides the bridge to the results of Section \ref{uni}.

We introduce the Harish-Chandra Schwartz space $\cC(G,V_\vp)$ of
$G$ in a way suitable for our purposes. Fix a base $\{X_i\}$ of $\gaaa$ and let $\cI_N$, $N\in\nat_0$, denote
the set of all multiindices $I=(i_1,\dots,i_{\dim(\gaaa)})$, $|I|\le N$. A multiindex $I\in \cI_N$ defines an element $X_I=X_1^{i_1}\dots X_{\dim(\gaaa)}^{i_{\dim(\gaaa)}}\in \cU(\gaaa)$. For $N\in\nat_0$ and a $V_\vp$-valued smooth function $f$ on $G$ we set 
$$
q_{N}(f)^2:=\sum_{(I,J)\in \cI_N\times\cI_N}\int_G |\log(a_g)^Nf(X_IgX_J)|^2\:dg\ .
$$
Then $\cC(G,V_\vp)=\{f\:|\:q_{N}(f)<\infty \ \forall N\in\nat_0\}$. The
seminorms $q_{N}$, $N\in\nat_0$ induce  the structure of
a Fr\'echet space on $\cC(G,V_\vp)$. For fixed $N\in\nat_0$ we define the Hilbert space $\cC^N(G,V_\vp)$ as the 
completion of $\cC(G,V_\vp)$ with respect to $q_{N}$. Then we have continuous inclusions $\cC(G,V_\vp)\subset \cC^N(G,V_\vp)\subset \cC^0(G,V_\vp)=L^2(G,V_\vp)$. 

Let $\chi$ be  the cut-off function $\chi$ constructed in Lemma \ref{lll}.
We consider $\chi$ as a right $K$-invariant function on $G$. 
\begin{ddd}\label{sscchh}
We define the Schwartz space on $\Gamma\backslash G$ by
$$ \cC(\Gamma\backslash G,\vp):=\{f\in L^2(\Gamma\backslash G,\vp)\:|\:
\chi f\in \cC(G,V_\vp)\}\ .$$
It inherits the structure of a Fr\'echet space
from $\cC(G,V_\vp)$. Similarly, we define intermediate Hilbert spaces $\cC^N(\Gamma\backslash G,\vp)$, $N\in\nat_0$, by
$$ \cC^N(\Gamma\backslash G,\vp):=\{f\in L^2(\Gamma\backslash G,\vp)\:|\:
\chi f\in \cC^N(G,V_\vp)\}\ .$$
Define $\cC^{-N}(\Gamma\backslash G,\vp)$ to be the conjugate linear dual of 
$\cC^N(\Gamma\backslash G,\vp)$. The space of tempered distributions $\cC^\prime(\Gamma\backslash G,\vp)$ is then
by definition the conjugate linear dual of $\cC(\Gamma\backslash G,\vp)$. 

Let $(\gamma,V_\gamma)$ be a finite-dimensional unitary representation of $K$.
Then $L^2(Y,V_Y(\gamma,\vp))\cong [L^2(\Gamma\backslash G,\vp)\otimes V_\gamma]^K$. Using this isomorphism we define for $*\in\{\emptyset,N,-N,\prime\}$
$$ \cC^*(Y,V_Y(\gamma,\vp)):=[\cC^*(\Gamma\backslash G,\vp)\otimes V_\gamma]^K \ .$$
\end{ddd}

\begin{lem}\label{komppp}
If $N$ is sufficiently large, then
the inclusion
$$i:\cC^N(\Gamma\backslash G,\vp)\hookrightarrow L^2(\Gamma\backslash G,\vp)$$
is Hilbert-Schmidt.
\end{lem}
\proof
Consider the inclusion $j:\cC^N(G,V_\vp)\hookrightarrow L^2(G,V_\vp)$. It is Hilbert-Schmidt for $N$ sufficiently large (see e.g. \cite{bernstein88}). Let
$\chi^\prime\in C_c^\infty(X\cup\Omega)$ be a second cut-off function
such that $\chi^\prime\chi=\chi$. Again, we consider $\chi^\prime$ as
a function on $G$. 
Then $p_{\chi^\prime}:L^2(G,V_\vp)\rightarrow L^2(\Gamma\backslash G,\vp)$ defined by 
$p_{\chi^\prime}(f)(x):=\sum_{g\in\Gamma} \vp(g) \chi^\prime(g^{-1}x)f(g^{-1}x)$ is continuous.
Let
$m_\chi:\cC^N(\Gamma\backslash G,\vp)\rightarrow \cC^N(G,V_\vp)$ be the operator
induced by multiplication with $\chi$. Now we can factorize $i$ over the 
Hilbert-Schmidt operator $j$: $i=p_{\chi^\prime}\circ j\circ m_\chi$. The lemma
follows.
\hB

In the following we choose $N$ suffiently large.
Fix an operator $\alpha$ providing the unitary equivalence (\ref{pl1}). 
By a theorem of Gelfand/Kostyuchenko (see \cite{bernstein88}) Lemma \ref{komppp} implies that
the composition
$$\cC^N(\Gamma\backslash G,\vp)\hookrightarrow L^2(\Gamma\backslash G,\vp)\stackrel{\alpha}{\longrightarrow} \int_{\hat G}^\oplus M_\pi\hat\otimes V_\pi \: d\kappa(\pi)$$
is pointwise defined, i.e., there exists a collection of continuous maps
$$\alpha_\pi:\cC^N(\Gamma\backslash G,\vp)\rightarrow M_\pi\hat\otimes V_\pi,\quad \pi\in\hat G$$ 
such that
for $f\in \cC^N(\Gamma\backslash G,\vp)$ we have $\alpha(f)(\pi)=\alpha_\pi(f)$.
By changing $\alpha_\pi$ on a set of $\pi$'s of measure
zero (w.r.t. $\kappa$) we
can assume that for all $\pi\in\hat G$ the map
$\alpha_\pi$ 
is an intertwining operator of $G$-representations.
Let
$$\beta_\pi: M_\pi\hat\otimes V_\pi\rightarrow \cC^{-N}(\Gamma\backslash G,\vp)$$ 
denote
the adjoint of $\alpha_\pi$. Note that $\cC^{-N}(\Gamma\backslash G,\vp)\subset\cC^{\prime}(\Gamma\backslash G,\vp)$. The composition of $\beta_\pi$ with this inclusion will also be denoted by $\beta_\pi$. 

Let $(\pi,V_\pi)$ be a representation of $G$ on a reflexive Banach space, and 
let $(\pi^\prime,V_{\pi^\prime})$ be its dual representation.
The space of distribution vectors $V_{\pi,-\infty}$ of $V_\pi$ is by definition $(V_{\pi^\prime,\infty})^\prime$, where the subscript $\infty$ indicates the transition to the subspace of smooth vectors, and the second dualization is with
respect to the canonical Fr\'echet topology on $V_{\pi^\prime,\infty}$. Then we have
the following inclusions of $G$-representations:
$V_{\pi,\infty}\subset V_\pi\subset V_{\pi,-\infty}$.

\begin{ddd}\label{laurent}
Let $(\pi,V_\pi)$ be a representation of $G$ on a reflexive Banach space, and
let $(\vp,V_\vp)$ be a finite-dimensional unitary representation of $\Gamma$. 
An invariant distribution vector $\phi\in {}^\Gamma (V_{\pi,-\infty}\otimes V_\vp)$ is called tempered (square integrable, resp.) if
for all $v\in V_{\pi^\prime,\infty}$ the function 
$$G\ni g\mapsto 
c_{\phi,v}(g):=\langle \phi, \pi^\prime(g) v \rangle \in V_\vp$$ belongs to $\cC^\prime(\Gamma\backslash G,\vp)$ ($L^2(\Gamma\backslash G,\vp)$). By
${}^\Gamma (V_{\pi,-\infty}\otimes V_\vp)_{d}\subset{}^\Gamma (V_{\pi,-\infty}\otimes V_\vp)_{temp}\subset{}^\Gamma (V_{\pi,-\infty}\otimes V_\vp)$
we denote the linear subspaces of square integrable and tempered invariant 
distribution vectors.
\end{ddd}

If $(\pi,V_\pi)$ is an admissible representation of finite length, 
we have the following characterizations and consequences of temperedness.

\begin{lem}\label{cassel}
Let $(\pi,V_\pi)$ be an admissible  $G$-representation of finite length on a reflexive Banach space, and let
$V_{\pi^\prime,K}\subset V_{\pi^\prime,\infty}$ be the underlying 
$(\gaaa,K)$-module of $K$-finite elements of the dual representation $V_{\pi^\prime}$. If $S\subset V_{\pi^\prime,K}$ is a generating set
and $\phi\in {}^\Gamma (V_{\pi,-\infty}\otimes V_\vp)$,
then the following conditions are equivalent:
\begin{enumerate}
\item $c_{\phi,v}\in \cC^\prime(\Gamma\backslash G,\vp)$ for all 
$v\in S$.
\item $\phi\in {}^\Gamma (V_{\pi,-\infty}\otimes V_\vp)_{temp}$\ .
\item $c_{\phi,v}\in \cC^\prime(\Gamma\backslash G,\vp)$ for all $ v\in V_{\pi^\prime,\infty}$,
 and the  map $c_{\phi,.}:V_{\pi^\prime,\infty} \rightarrow  \cC^\prime(\Gamma\backslash G,\vp)$ is continuous.
\end{enumerate}
The analogous assertions for ${}^\Gamma (V_{\pi,-\infty}\otimes V_\vp)_{d}$
and $L^2(\Gamma\backslash G,\vp)$ hold true, too.
\end{lem}
\proof 
The lemma is a consequence of the globalization theory of Casselman and Wallach
(\cite{casselman89}, \cite{wallach92}, Ch. 11). 
The conclusions 3. $\Rightarrow$ 2. $\Rightarrow$ 1. are obvious. We outline the argument for 1.$\Rightarrow$ 3.

Assume 1. 
Since the $(\gaaa,K)$-module 
$V_{\pi^\prime,K}$ is finitely generated and admissible, and $K$-finite matrix coefficients 
satisfy elliptic differential equations one can show that there exists $N\in\nat_0$ such that $c_{\phi,v}\in \cC^{-N}(\Gamma\backslash G,\vp)_K$ for all 
$v\in V_{\pi^\prime,K}$. Here $\cC^{-N}(\Gamma\backslash G,\vp)_K$ denotes the subspace of $K$-finite smooth vectors of $\cC^{-N}(\Gamma\backslash G,\vp)$. 
But $\cC^{-N}(\Gamma\backslash G,\vp)$ is a  Hilbert space on which $G$ acts continuously. The theorem of
Casselman and Wallach now states
that the $(\gaaa,K)$-module homomorphism $V_{\pi^\prime,K}\ni v \mapsto c_{\phi,v}\in \cC^{-N}(\Gamma\backslash G,\vp)_K$ extends to a continuous
homomorphism $V_{\pi^\prime,\infty}\rightarrow \cC^{-N}(\Gamma\backslash G,\vp)_\infty$. Condition 3 follows. The argument for ${}^\Gamma (V_{\pi,-\infty}\otimes V_\vp)_{d}$ is essentially the same.
\hB 

The notion of temperedness is compatible with the notion
of a tempered (irreducible) representation $(\pi,V)$ of $G$ as follows. For a moment let $\Gamma$ and $V_\vp$ be trivial, and define $V_{-\infty,temp}$ as above. Then  $\pi$ is tempered iff $V_{-\infty,temp}\subset V_{-\infty}$ is dense. Similarly, $\pi$ is square integrable, i.e., belongs to the discrete series, iff $V_{-\infty,d}\subset V_{-\infty}$ is dense. The above lemma in mind it is
not difficult to see that this characterization of temperedness and square-integrability is equivalent to the various definitions appearing in
the literature. In fact, more is true:

\begin{lem}\label{arthur}
If $(\pi,V)$ is an admissible $G$-representation of finite length on a reflexive Banach space, then $\pi$ is tempered iff $V_{-\infty,temp}=V_{-\infty}$. Moreover,
if $v\in V^\prime_\infty$ is fixed, then the matrix coefficient map
$$  V_{-\infty}\ni \phi\mapsto c_{\phi,v}\in \cC^\prime(G)$$
is continuos.
\end{lem}
\proof
The lemma follows from Lemma 10 in \cite{arthur70} which is based on subtle estimates of $K$-finite matrix coefficients. However, we would like
to indicate a different argument which is more in the spirit of the proof
of Lemma \ref{cassel}.

Fix $v\in V^\prime_\infty$. Then for all 
$\phi\in V_K$ the asymptotic expansions of 
matrix coefficients (\cite{wallach88}, 4.4.3.) give  
$c_{\phi,v}\in \cC^\prime(G)$. 
Thus by Lemma \ref{cassel} we have $v\in V^\prime_{-\infty,temp}$.
As in the proof of Lemma
\ref{cassel}
the map $\phi\mapsto c_{\phi,v}$ extends to a continuous map
$c_{.,v}:V_\infty\rightarrow\cC^{-N}(G)_\infty$ for some $N\in\nat_0$.
Since $\cC^{-N}(G)_\infty$ is a Hilbert space this map has a continuous right inverse,
the adjoint of which is the continuous extension of $c_{.,v}$ to $V_{-\infty}$
$$c_{.,v}:V_{-\infty}\rightarrow\cC^{N}(G)_{-\infty}\subset \cC^\prime(G)\ .$$
This proves the lemma. \hB

We now return to the discussion of the map $\beta_\pi$. 

\begin{lem}\label{lebal}
For any
$\pi\in\hat G$ there is an embedding
$$i_\pi : M_\pi\hookrightarrow {}^\Gamma (V_{\pi^\prime,-\infty}\otimes V_\vp)_{temp}$$
such that for $m\in M_\pi$ and $v\in V_{\pi,\infty}$
we have $\beta_\pi(m\otimes v)=c_{i_\pi(m),v}$.
\end{lem}
\proof
Let $N$ be sufficiently large.
Fix $m\in M_\pi$ and consider the $G$-intertwining operator $F_m : V_\pi\rightarrow \cC^{-N}(\Gamma\backslash G,\vp)$ given by $F_m(v):=\beta_\pi(m\otimes v)$. Then 
$$F_m(V_{\pi,\infty})\subset \cC^{-N}(\Gamma\backslash G,\vp)_\infty
\subset\cC^{-N}(\Gamma\backslash G,\vp)\cap C^{\infty}(\Gamma\backslash G,\vp)\ .$$
In particular, elements of $F_m(V_{\pi,\infty})$ can be evaluated at the identity
$e\in G$. We define $i_\pi$ by $\langle i_\pi(m),v\rangle:=F_m(v)(e)$,
$v\in V_{\pi,\infty}$. The assertion of the lemma is now obvious.
\hB 

Our first concretization of the abstract Plancherel decomposition (\ref{pl1})  is given by the following corollary.
\begin{kor}\label{abel}
There exists a collection of Hilbert spaces 
$N_\pi\subset
{}^\Gamma (V_{\pi^\prime,-\infty}\otimes V_\vp)_{temp}$, $\pi\in\hat G$,  and a direct integral $$ \int_{\hat G}^\oplus N_\pi\hat\otimes V_\pi \: d\kappa(\pi) $$ such that the
following holds:
\begin{enumerate}
\item The matrix coefficient map $c:{}^\Gamma (V_{\pi^\prime,-\infty}\otimes V_\vp)_{temp}\otimes V_{\pi,\infty}\rightarrow \cC^\prime(\Gamma\backslash G,\vp)$ gives rise to a map 
$$ c_\pi: N_\pi\hat\otimes V_\pi\rightarrow \cC^\prime(\Gamma\backslash G,\vp)\ .$$
\item Let $\cF_\pi: \cC(\Gamma\backslash G,\vp)\rightarrow N_\pi\hat\otimes V_\pi$ be the adjoint of $c_\pi$. Then the collection of maps $\cF_\pi$ extends to
a unitary equivalence
$$ \cF: L^2(\Gamma\backslash G,\vp)\stackrel{\cong}{\longrightarrow} \int_{\hat G}^\oplus N_\pi\hat\otimes V_\pi \: d\kappa(\pi)\ . $$
\item If ${}^\Gamma (V_{\pi^\prime,-\infty}\otimes V_\vp)_{d}\not=0$,
then $N_\pi={}^\Gamma (V_{\pi^\prime,-\infty}\otimes V_\vp)_{d}$.
\end{enumerate}
In particular, the Plancherel measure $\kappa$ is supported on the set
$$\{\pi\in\hat G\:|\: {}^\Gamma (V_{\pi^\prime,-\infty}\otimes V_\vp)_{temp}\not=0\}$$
and $\kappa(\{\pi\})\not=0$ iff ${}^\Gamma (V_{\pi^\prime,-\infty}\otimes V_\vp)_{d}\not=0$.

The Plancherel measure and the scalar product on the subspace ${}^\Gamma (V_{\pi^\prime,-\infty}\otimes V_\vp)_d\subset N_\pi$ can be chosen such
that $\kappa(\{\pi\})=1$, if $(V_{\pi^\prime,-\infty}\otimes V_\vp)_d\not=0$, and such that $c_\pi$ induces an isometric embedding of 
$\:{}^\Gamma (V_{\pi^\prime,-\infty}\hat\otimes V_\vp)_d\otimes V_\pi$ into 
$L^2(\Gamma\backslash G,\vp)$.
\end{kor}
\proof Set $N_\pi:=i_\pi(M_\pi)$.
\hB

The plan of the rest of the paper is now as follows. In the next section we
determine the spaces $(V_{\pi,-\infty}\otimes V_\vp)_{temp}$ and $(V_{\pi,-\infty}\otimes V_\vp)_d$ for all $\pi\in\hat G$. This is
based on the results of the Section \ref{uni}. 
In Section \ref{wxa} we study wave packets of Eisenstein series.
It turns out that
they span the orthogonal complement to the discrete subspace
$$L^2(\Gamma\backslash G,\vp)_d:=\bigoplus_{\{\pi\in\hat G\:|\: (V_{\pi,-\infty}\otimes V_\vp)_d\not=0\}}^{\mbox {\tiny Hilbert}}
\im\: c_\pi\subset L^2(\Gamma\backslash G,\vp)\ .$$ 
The proof of this fact heavily depends on our
a priori knowledge of the support of the Plancherel measure.
The last section contains the summary of our results, including the determination of the scalar products on $N_\pi$ and of the Placherel measure, as well
as the consequences for the spectral theory of the Casimir operator acting
on sections of the locally homogeneous vector bundle $V_Y(\gamma,\vp)$ over
the Kleinian manifold $Y$.

\section{Tempered invariant distribution vectors}\label{relsec}

In this section we determine the tempered and square integrable invariant distribution
vectors for all $\pi\in\hat G$. 

First we need a rough classification of the unitary dual $\hat G$. Recall
the notions of temperedness and square integrability of an irreducible
representation (see the discussion following Lemma \ref{cassel}). The
unitary dual is a disjoint union of the dicrete series, the unitary principal series, and the complementary series
$$\hat G=\hat G_d\cup\hat G_{u}\cup\hat G_{c}\ ,$$
where
\begin{eqnarray*}
\hat G_d     &:=& \{\pi\in\hat G\:|\: V_\pi \mbox{ is square integrable}\}\ ,\\
\hat G_{u} &:=& \{\pi\in\hat G\:|\: V_\pi \mbox{ is tempered}\}\setminus\hat G_d\ ,\\
\hat G_{c}  &:=& \{\pi\in\hat G\:|\: V_\pi \mbox{ is not tempered}\}\ .
\end{eqnarray*}

The discrete series $\hat G_d$ has been determined by Harish-Chandra. 
It is empty iff $X=H^n$, $n$
odd. In the other cases one can choose a Cartan subalgebra $\haaa$ of $\gaaa$ which is contained in $\kaaa$. An infinitesimal character $\chi_\lambda$, $\lambda\in\haaa^*_\C$, is called regular if no expression of the
form (\ref{jemi}) vanishes. Let $W(\kaaa_\C,\haaa_\C)$ be the Weyl group of $K$. Then $\hat G_d$ can be parametrized by the Harish-Chandra parameters
\begin{equation}\label{joyce}
\{\lambda\in\haaa^*_\C\:|\: \chi_\lambda \mbox{ is regular and integral}\}/W(\kaaa_\C,\haaa_\C)
\end{equation}
such that $\cZ$ acts on $\pi\in\hat G_d$ with infinitesimal character $\chi_\lambda$, and $\pi$ has the minimal $K$-type with highest weight
$\lambda+\rho_\gaaa-2\rho_\kaaa$ (see e.g. \cite{wallach88}, Ch. 6 and Ch. 8). Here
$\rho_\gaaa$ and $\rho_\kaaa$ are the half sums of the positive (w.r.t. $\lambda$) roots of $\haaa$
in $\gaaa$ and $\kaaa$, respectively. In particular, for any $\gamma\in\hat K$
there are only finitely many discrete series representations containing the
$K$-type $\gamma$. Strictly speaking, (\ref{joyce}) parametrizes the discrete series representations for the linear group $G$ with
Lie algebra $\gaaa$ which has a
simply connected complexification $G_\C$. In general, the parametrization remains valid, if one sharpens the notion of integrality.

The set $\hat G_u$ consists of the unitary principal series representations $\pi^{\sigma,\lambda}$, $\sigma\in\hat M$, $\Ree(\lambda)=0$. They are irreducible unless $\sigma=\sigma^w$, $p_\sigma(0)=0$ and $\lambda=0$. In the latter case we have
$\pi^{\sigma,0}=\pi^{\sigma,+}\oplus\pi^{\sigma,-}$, where the 
irreducible representations $\pi^{\sigma,\pm}$, called the non-degenerate limits
of discrete series, are the $\pm 1$-eigenspaces of $J_{\sigma,0}$. All
equivalences between these representations are induced by the intertwining operators $J_{\sigma,\lambda}$. For all that see \cite{knapp86}, Ch. XIV.

Though also $\hat G_c$ is completely known (see \cite{baldonisilvabarbasch83} and
the references therein) for our purposes less information is sufficient. The Langlands classification (see \cite{wallach88}, Ch. 5) associates to any irreducible non-tempered representation $(\pi,V_\pi)$ a unique Langlands parameter
$(\sigma,\lambda)$, $\sigma\in\hat M$, $\lambda\in\aca$, $\Ree(\lambda)>0$, such
that $V_{\pi,\pm\infty}$ is equivalent to the unique irreducible subrepresentation of the principal series representation $\pi^{\sigma,\lambda}$
acting on $C^{\pm\infty}(\partial X,V(\sigma_\lambda))$. We denote this subrepresentation by $(\bar\pi^{\sigma,\lambda},I^{\sigma,\lambda}_{\pm\infty})$.
It is the image of $J_{\sigma^w,-\lambda}$. It is also the unique
irreducible quotient of $C^{\infty}(\partial X,V(\sigma^w_{-\lambda}))$. If $\pi\in\hat G_c$, then $\sigma=\sigma^w$, $\lambda\in\aaaa^*$ (e.g. \cite{knapp86}, Thm. 16.6.) and 
$p_\sigma(0)\not=0$ (\cite{baldonisilva81}, Thm. 6.1). An invariant pre-Hilbert structure $(.,.)$ on $I^{\sigma,\lambda}_{\infty}$ realized as a quotient of
$C^{\infty}(\partial X,V(\sigma_{-\lambda}))$ can now  be described as follows: 
Let $(.,.)_0$ be the invariant sesquilinear
pairing between $C^{\infty}(\partial X,V(\sigma_{-\lambda}))$ and  $C^{\infty}(\partial X,V(\sigma_\lambda))$. Then
\begin{equation}\label{jam} 
([f],[g]):=(f,J_{\sigma,-\lambda}(g))_0\ , 
\end{equation}
where we have represented $[f],[g] \in I^{\sigma,\lambda}_\infty$ by $f,g\in C^{\infty}(\partial X,V(\sigma_{-\lambda}))$. Indeed, since ${}^tJ_{\sigma,-\lambda}=J_{\sigma,-\lambda}$ the pairing (\ref{jam}) 
is hermitian. By
$I^{\sigma,\lambda}$ we denote the Hilbert space completion of $I^{\sigma,\lambda}_{\infty}$ with respect to $(.,.)$.

We recall the relation between the Poisson transform $P_\lambda^T$ and the matrix coefficients of principal series representations.
For later reference we will state it as a lemma. Its verification
is a standard computation with integral formulas.

\begin{lem}\label{circe}
Let $\vp$, $\sigma$, $\gamma$ be  finite-dimensional representations of $\Gamma$, $M$, and $K$, respectively, $\lambda\in\aca$ and $T\in \Hom_M(V_\sigma,V_\gamma)$. We consider the Poisson transform
$$P^T_\lambda: C^{-\infty}(\partial X,V(\sigma_\lambda,\vp))\rightarrow
C^\infty(X,V(\gamma,\vp))\cong [C^\infty(G,V_\vp)\otimes V_\gamma]^K\ .$$
Then for any $v\in V_{\tilde\gamma}$,
$\phi\in C^{-\infty}(\partial X,V(\sigma_\lambda,\vp))$, $g\in G$ we have
$$
\langle P^T_\lambda\phi(g),v\rangle=c_{\phi,v_T}(g) \in V_\vp\ ,
$$
where $v_T\in C^{\infty}(\partial X,V(\tilde\sigma_{-\lambda}))$ is the element
defined by Frobenius reciprocity $v_T(k):={}^tT\tilde\gamma(k^{-1})v$, $k\in K$.
\hB
\end{lem}

Fix a finite-dimensional unitary representation $(\vp,V_\vp)$ of $\Gamma$. 
Recall the decomposition
${}^\Gamma C^{-\infty}(\Lambda,V(\sigma_\lambda,\vp))=
E_\Lambda(\sigma_\lambda,\varphi)\oplus U_\Lambda(\sigma_\lambda,\varphi)$
from Proposition \ref{engels}.
\begin{prop}\label{nonte}
Let $\sigma\in\hat M$, $\lambda\in\aca$ with $\Ree(\lambda)>0$. Then 
$${}^\Gamma(I^{\sigma,\lambda}_{-\infty}\otimes V_\vp)_d=
{}^\Gamma(I^{\sigma,\lambda}_{-\infty}\otimes V_\vp)_{temp}=
{}^\Gamma C^{-\infty}(\Lambda,V(\sigma_\lambda,\vp))=
E_\Lambda(\sigma_\lambda,\varphi)\oplus U_\Lambda(\sigma_\lambda,\varphi)\ .$$
If one of these spaces is non-zero, then $\bar\pi^{\sigma,\lambda}\in\hat G_c$.
\end{prop}
\proof
Observe that ${}^\Gamma(I^{\sigma,\lambda}_{-\infty}\otimes V_\vp)_{temp}\subset
{}^\Gamma C^{-\infty}(\partial X,V(\sigma_\lambda,\vp))_{temp}$. Fix a minimal $K$-type
$\gamma$ of $C^{\infty}(\partial X,V(\sigma_\lambda))$ together with a non-trivial $T\in \Hom_M(V_\sigma,V_\gamma)$. We consider the injective Poisson
transform
$$P=P^T_\lambda: C^{-\infty}(\partial X,V(\sigma_\lambda,\vp))\rightarrow
C^\infty(X,V(\gamma,\vp))\ .$$
Recall the definition of $v_T$ from Lemma \ref{circe}.
If $\phi\in {}^\Gamma C^{-\infty}(\partial X,V(\sigma_\lambda,\vp))_{temp}$, then by definition $c_{\phi,v_T}\in \cC^\prime(\Gamma\backslash G,\vp)$ for all
$v\in V_{\tilde\gamma}$, hence $P\phi\in\cC^\prime(Y,V_Y(\gamma,\vp))$. 

Let 
$f\in C^{\infty}(\partial X,V((\tilde\gamma_{|M})_{-\lambda},\tilde\vp))$ with
$\supp(f)\subset\Omega$. We want to show that $\langle \phi,{}^tTf\rangle=0$.
As in the proof of Theorem \ref{pokm} we extend $f$ to a section $\tilde f\in
C^\infty(X,V(\tilde\gamma,\tilde\vp))$. By Corollary \ref{lim} 
there is a constant $C$ such that
$$
\langle \phi,{}^tTf\rangle=
C\lim_{a\to\infty}a^{\rho-\lambda} \int_K  \langle P\phi(ka),f(k)\rangle\ dk
=C\lim_{n\to\infty} \int_G \langle P\phi(x),\tilde f_n(x)\rangle\ dx,
$$
where $\tilde f_n(x):=a_x^{-(\lambda+\rho)}\psi(\log(a_x-n))\tilde f(x)$ for some
$\psi\in C^\infty_c(0,1)$ satisfying $\int_0^1\psi(t)dt=1$. Define
$F_n\in C^\infty_c(Y,V_Y(\tilde\gamma,\tilde\vp))\subset\cC(Y,V_Y(\tilde\gamma,\tilde\vp))$
by
$$F_n(x):= \sum_{g\in\Gamma} \tilde\vp(g)\tilde f_n(g^{-1}x)\ .$$
We claim that $\lim_{n\to\infty} F_n=0$ in $\cC(Y,V_Y(\tilde\gamma,\tilde\vp))$.
Using that $\supp(f)\subset\Omega$ we find a finite subset $L\subset \Gamma$
such that 
$$\chi(x)F_n(x)= \chi(x) \sum_{g\in L} \tilde\vp(g)\tilde f_n(g^{-1}x)=\sum_{g\in L} \tilde\vp(g) ((g^{-1})^*\chi \tilde f_n)(g^{-1}x)\ ,$$ where $\chi$ is the cut-off function as in Lemma \ref{lll}. 
We have $\lim_{n\to\infty} \tilde f_n=0$ in $\cC(X,V(\tilde\gamma,\tilde\vp))$ and hence 
 $(g^{-1})^*\chi f_n\to 0$ in $\cC(X,V(\tilde\gamma,\tilde\vp))$.
This shows the claim.

 Since $P\phi\in
\cC^\prime(Y,V_Y(\gamma,\vp))$ we obtain
$$\langle \phi,{}^tTf\rangle=C\lim_{n\to\infty} \int_{\Gamma\backslash G} \langle P\phi(x),F_n(x)\rangle\ dx=C\lim_{n\to\infty}\langle P\phi,F_n\rangle=0\ .$$
This proves that ${}^\Gamma C^{-\infty}(\partial X,V(\sigma_\lambda,\vp))_{temp}\subset{}^\Gamma C^{-\infty}(\Lambda,V(\sigma_\lambda,\vp))$. 

On the other hand, we have by Proposition \ref{firtu} that 
$P\left({}^\Gamma C^{-\infty}(\Lambda,V(\sigma_\lambda,\vp))\right)\subset L^2(Y,V_Y(\gamma,\vp))$. Since the elements $v_T$, $v\in V_{\tilde\gamma}$, generate the $(\gaaa,K)$-module $C^{\infty}(\partial X,V(\tilde\sigma_{-\lambda}))_K$
it follows from Lemma \ref{cassel} and Lemma \ref{circe} that $c_{\phi,f}\in L^2(\Gamma\backslash G,\vp)$
for all $f\in C^{\infty}(\partial X,V(\tilde\sigma_{-\lambda}))$. Thus ${}^\Gamma C^{-\infty}(\Lambda,V(\sigma_\lambda,\vp))\subset {}^\Gamma C^{-\infty}(\partial X,V(\sigma_\lambda,\vp))_d$. It remains to show that
${}^\Gamma C^{-\infty}(\partial X,V(\sigma_\lambda,\vp))_d\subset{}^\Gamma(I^{\sigma,\lambda}_{-\infty}\otimes V_\vp)_d$. Indeed, let $\phi\in {}^\Gamma C^{-\infty}(\partial X,V(\sigma_\lambda,\vp))_d$
and consider the $G$-map
$$c_{\phi,.}: C^{\infty}(\partial X,V(\tilde\sigma_{-\lambda}))\rightarrow
L^2(\Gamma\backslash G,\vp)\ .$$
Since the target space is a unitary representation of $G$ the image  of $c_{\phi,.}$ decomposes into a direct
sum of irreducible representations. But $C^{\infty}(\partial X,V(\tilde\sigma_{-\lambda}))$ has the unique irreducible quotient $$C^{\infty}(\partial X,V(\tilde\sigma_{-\lambda}))/(I^{\sigma,\lambda}_{-\infty})^\perp\ .$$ 
Thus $c_{\phi,.}$ factorizes over this quotient, and hence $\phi\in
{}^\Gamma(I^{\sigma,\lambda}_{-\infty}\otimes V_\vp)_d$. Now, if $\phi$ is non-trival, then we can pull back the invariant pre-Hilbert structure from
$L^2(\Gamma\backslash G,\vp)$ to this quotient. By duality this induces
an invariant scalar product on $I^{\sigma,\lambda}_{\infty}$. Hence $\bar\pi^{\sigma,\lambda}\in\hat G_c$. In view of the chain of inclusions
\begin{eqnarray*}
{}^\Gamma(I^{\sigma,\lambda}_{-\infty}\otimes V_\vp)_{temp}&\subset& {}^\Gamma C^{-\infty}(\partial X,V(\sigma_\lambda,\vp))_{temp}\subset
{}^\Gamma C^{-\infty}(\Lambda,V(\sigma_\lambda,\vp))\subset\\
&\subset& {}^\Gamma C^{-\infty}(\partial X,V(\sigma_\lambda,\vp))_d
\ \subset\: {}^\Gamma(I^{\sigma,\lambda}_{-\infty}\otimes V_\vp)_d
\end{eqnarray*}
the proof of the proposition is now complete.
\hB  

\begin{kor}\label{xerxes}
The space ${}^\Gamma(I^{\sigma,\lambda}_{-\infty}\otimes V_\vp)_d$ is 
finite-dimensional. It is non-trivial iff $\lambda$ belongs to the finite set
$PS(\sigma,\vp)\setminus\{0\}\subset (0,\delta_\Gamma]$.
\end{kor}
\proof Combine Proposition \ref{nonte} with Proposition \ref{upperbound}.
\hB

The case $\sigma=1$ is particularly interesting. In the following table
we give the set of $\lambda>0$ with $I^{1,\lambda}\in\hat{G}$.\\
\centerline{\begin{tabular}{|c||c|c|c|c|}
\hline
$X$&$\R H^n$&$\C H^n$&$\HH H^n$&$\OO H^n$\\
\hline
$\lambda$ &$(0,\rho]$&$(0,\rho]$&$(0,\rho-2\alpha]\cup\{\rho\}$&$(0,\rho-6\alpha]\cup\{\rho\}$\\
\hline
\end{tabular}\ . }

If $\delta_\Gamma>0$, then by Proposition \ref{nonte} the representation $I^{1,\delta_\Gamma}$
is unitary.  
This leads to the restriction of the set of possible values of $\delta_\Gamma$ found by Corlette \cite{corlette90} (see Section \ref{fiert}).

\begin{lem}\label{zeus}
If $(\pi,V_\pi)\in\hat G$ is tempered, then 
${}^\Gamma(V_{\pi,-\infty}\otimes V_\vp)_{temp}={}^\Gamma(V_{\pi-\infty}\otimes V_\vp)$. Moreover, for any $v\in V_{\tilde\pi,\infty}$ the map
$${}^\Gamma(V_{\pi,-\infty}\otimes V_\vp)\ni\phi\mapsto\ c_{\phi,v}\in \cC^\prime(\Gamma\backslash G,\vp)$$
is continuous. If, in addition, $\pi=\pi^{\sigma,\lambda}$, $\lambda\not=0$ imaginary, then  
$ext$ identifies $C^{-\infty}(B,V_B(\sigma_\lambda,\vp))$ with ${}^\Gamma(V_{\pi,-\infty}\otimes V_\vp)_{temp}$.
\end{lem}
\proof 
Observe that there is a natural inclusion ${}^\Gamma\cC^\prime(G,V_\vp)\hookrightarrow \cC^\prime(\Gamma\backslash G,\vp)$
induced by the adjoint of the multiplication by $\chi$. The first assertions of the lemma now follow from Lemma \ref{arthur}. The
last one follows from Lemmas \ref{mainkor} and \ref{extregatim}.
\hB

\begin{prop}\label{hera}
Let $\sigma\in\hat M$. If $\lambda\not=0$ is imaginary, then ${}^\Gamma C^{-\infty}(\partial X,V(\sigma_\lambda,\vp))_{d}=0$. For $\lambda=0$ we have 
$${}^\Gamma C^{-\infty}(\partial X,V(\sigma_0,\vp))_{d}=U_\Lambda(\sigma_0,\varphi)\ .$$ 
If this finite-dimensional space is non-trivial, then
$\sigma\not=1$ and $0\in I_\sigma$.
\end{prop}
\proof
We claim that for $\Ree(\lambda)=0$
\begin{equation}\label{pene}
{}^\Gamma C^{-\infty}(\partial X,V(\sigma_\lambda,\vp))_{d}\subset
{}^\Gamma C^{-\infty}(\Lambda,V(\sigma_\lambda,\vp))\ .
\end{equation}
We now prove the proposition assuming the claim. If $\lambda\not=0$, then by Lemma \ref{extregatim} we have ${}^\Gamma C^{-\infty}(\Lambda,V(\sigma_\lambda,\vp))=0$ and hence
$C^{-\infty}(\partial X,V(\sigma_\lambda,\vp))_{d}=0$.
Now consider the case $\lambda=0$.
If $\phi\in C^{-\infty}(\partial X,V(\sigma_\lambda,\vp))_{d}$, then the claim implies $res_\Omega\phi=0$.
By the  asymptotic expansion Lemma \ref{poisasm}, 2., we see that $\phi$ is 
square-integrable iff $res_\Omega\circ \hat{J}_{0}\phi=0$. 
This proves that ${}^\Gamma C^{-\infty}(\partial X,V(\sigma_0,\vp))_{d}=U_\Lambda(\sigma_0,\varphi)$.
If the latter space is non-trivial, then $\sigma\not=1$ and $0\in I_\sigma$
by Corollary \ref{ujn}.

The proof of (\ref{pene}) is analogous to the proof of Proposition \ref{nonte}.
Let $\phi\in {}^\Gamma C^{-\infty}(\partial X,V(\sigma_\lambda,\vp))$. One
has to show that for some Poisson transform $P$ the condition $P\phi\in L^2(Y,V_Y(\gamma,\vp))$ implies $res(\phi)=0$. 
We give the argument for the
most involved case that $\sigma$ is Weyl-invariant, $\lambda=0$, and $\hat J_0$ regular.
The similar treatment of the remaining cases is left to the reader.

In fact, in this case we show at once that $res(\phi)=0$ and $res(\hat J_0\phi)=0$.
Recall that $\pi^{\sigma,0}$ splits into the $\pm1$-eigenspaces $\pi^{\sigma,\pm}$ of $J_0$. Let $\gamma_{\pm}$ be the minimal $K$-type 
of $\pi^{\sigma,\pm}$. Choose embeddings $T^\pm\in\Hom_M(V_\sigma,V_{\gamma_\pm})$ and set  $P^\pm:=P^{T^\pm}_0\otimes\id$. Define $t^\pm\in\Hom_M(V_{\tilde\sigma},V_{\tilde\gamma_\pm})$ by ${}^t(T^\pm)^w\circ t^\pm =\id_{V_{\tilde\sigma}}$. Using (\ref{murphy}) and Corollary \ref{lim} we
have for $f\in C^{\infty}(\partial X,V(\tilde\sigma_0,\tilde\vp))$
\begin{equation}\label{ghy6} \langle c_{\gamma_{\pm}}(0) T^{\pm} \phi+ (T^\pm)^w \hat{J}_0 \phi,t^\pm f\rangle=
\lim_{a\to\infty} a^{\rho}\int_K \langle P^\pm\phi(ka),t^\pm f(k)\rangle dk\ .\end{equation}
We can rewrite the left hand side as follows:
\begin{eqnarray*}
\langle c_{\gamma_{\pm}}(0) T^{\pm} \phi+ (T^\pm)^w \hat{J}_0 \phi,t^\pm f\rangle&=&
\langle \pm c_{\sigma}(0) (T^{\pm})^w \phi+ (T^\pm)^w \hat{J}_0 \phi,t^\pm f\rangle\\
&=& \langle \pm c_{\sigma}(0)  \phi+ \hat{J}_0 \phi, {}^t(T^\pm)^w t^\pm f\rangle\\
&=& \langle \pm c_{\sigma}(0)  \phi+ \hat{J}_0 \phi, f\rangle\ .
\end{eqnarray*}
We continue tranforming the right-hand side of (\ref{ghy6}).
\begin{eqnarray*}
\langle \pm c_\sigma(0)\phi +\hat J_0\phi,f\rangle&=&
\lim_{n\to\infty}\frac{1}{n}\int_0^n \ee^{t\rho} \int_K  \langle P^{\pm}\phi(k\exp(tH)),t^\pm f(k)\rangle\ dk\:dt\\
&=&\lim_{n\to\infty} \int_G \langle P\phi(x), f^\pm_n(x)\rangle\ dx,
\end{eqnarray*}
where $H\in\aaaa_+$ is the unit vector, and the compactly supported sections
$f^\pm_n\in L^2(X,V(\tilde\gamma_\pm,\tilde\vp))$ are defined by $f^\pm_n(k_1ak_2):=\frac{1}{n}a^{-\rho}\chi_{(0,n]}(\log(a))\tilde\gamma_\pm(k_2^{-1})t^\pm f(k_1)$. Here $\chi_{(0,n]}$ denotes the characteristic function of the interval
$(0,n]$. 

Assume now that $\supp(f)\subset\Omega$ and consider
$F^\pm_n\in L^2(Y,V_Y(\tilde\gamma_\pm,\tilde\vp))$
given by
$$F^\pm_n(x):= \sum_{g\in\Gamma} \tilde\vp(g) f^\pm_n(g^{-1}x)\ .$$
Then as in the proof of Proposition \ref{nonte} we see that
$\lim_{n\to\infty} F^\pm_n=0$ in $L^2(Y,V_Y(\tilde\gamma_\pm,\tilde\vp))$. Since $P^\pm\phi\in
L^2(Y,V_Y(\gamma_\pm,\vp))$ we obtain
$$\langle \pm c_\sigma(0)\phi+ \hat J_0\phi,f\rangle=\lim_{n\to\infty} \int_{\Gamma\backslash G} \langle P^\pm\phi(x),F^\pm_n(x)\rangle\ dx= \lim_{n\to\infty}\langle P^\pm\phi,F^\pm_n\rangle=0\ .$$
Since $c_\sigma(0)\not=0$ this proves
$res(\phi)=0$ and $res(\hat J_0\phi)=0$.
\hB 

The next lemma is independent of the theory
of tempered invariant distribution vectors.

\begin{lem}\label{skylla}
Let $(\sigma,\lambda),(\tau,\mu)\in\hat M\times\aca$, and let 
$$A:C^{-\infty}(\partial X,V(\sigma_{\lambda}))\rightarrow C^{-\infty}(\partial X,V(\tau_{\mu}))$$
be a $G$-intertwining operator. Let $(\vp,V_\vp)$ be a finite-dimensional (not necessarily unitary) representation of $\Gamma$. We consider the operator 
$$A_{\Gamma,\vp}:{}^\Gamma C^{-\infty}(\partial X,V(\sigma_{\lambda},\vp))\rightarrow {}^\Gamma C^{-\infty}(\partial X,V(\tau_{\mu},\vp))$$
induced by $A\otimes\id$. Then $\im(A_{\Gamma,\vp})$ is infinite-dimensional unless $\im(A)$ is 
finite-dimensional.
\end{lem}
\proof
The operator $A$ restricts to a continuous operator
$$A:C^\infty(\partial X,V(\sigma_{\lambda}))\rightarrow C^{\infty}(\partial X,V(\tau_{\mu}))\ .$$
In fact, $A$ induces an intertwining operator of the underlying $(\gaaa,K)$-modules which canonically extends to a continuous $G$-map 
between the spaces of smooth sections by the globalization theory of Casselman and Wallach \cite{wallach92}, Ch. 11. 

If $f$ is a distribution section of a vector bundle over some manifold $U$,
then let $\singsupp(f)\subset U$ denote the singular support of $f$.
We claim  that $\singsupp\: A(f)\subset \singsupp(f)$ for all
$f\in C^{-\infty}(\partial X,V(\sigma_{\lambda},\vp))$.

Consider the delta distribution $\delta v\in C^{-\infty}(\partial X,V(\sigma_{\lambda}))$ at $x_0:=MAN\in\partial X$ with vector part $v\in V_{\sigma}$. Since $\delta v$ is $N$-invariant we have  
$A(\delta v)\in {}^N C^{-\infty}(\partial X,V(\tau_{\mu}))$. The Bruhat
decomposition $\partial X=Nwx_0\cup \{x_0\}$ implies that $A(\delta v)$ is smooth outside $x_0$. 
Now let $f\in C^{-\infty}(\partial X,V(\sigma_{\lambda}))$. Then there
exists $\tilde f\in C_c^{-\infty}(G)$ such that 
$\singsupp(\tilde f) P=\singsupp (f)$ and
$\pi^{\sigma,\lambda}(\tilde f)\delta v=f$.
It follows that $A(f) = \pi^{\tau,\mu}(\tilde f) A(\delta v)$, and hence
$\singsupp \:A(f)\subset \singsupp(\tilde f) P \{x_0\}=\singsupp(f)$.
This shows the claim. In particular, 
$A_{\Gamma,\vp}$ maps ${}^\Gamma C^{-\infty}_\Omega(\partial X,V(\sigma_{\lambda},\vp))$ into ${}^\Gamma C^{-\infty}_\Omega(\partial X,V(\tau_{\mu},\vp))$ (see Lemma \ref{ddeenn} for notation).

Assume now that
$\dim\:\im(A_{\Gamma,\vp})<\infty$. By Lemma \ref{ddeenn} the space
$A_{\Gamma,\vp}\left({}^\Gamma C^{-\infty}_\Omega(\partial X,V(\sigma_{\lambda},\vp))\right)$ is dense in $\im(A_{\Gamma,\vp})$,
hence $A_{\Gamma,\vp}\left({}^\Gamma C^{-\infty}_\Omega(\partial X,V(\sigma_{\lambda},\vp))\right)=\im(A_{\Gamma,\vp})$. We conclude that\linebreak[4]
$\im(A_{\Gamma,\vp})\subset {}^\Gamma C^{-\infty}_\Omega(\partial X,V(\tau_{\mu},\vp))$.  

Without loss of generality we may assume that $x_0\in\Omega$. Choose $0\not=w\in V_\vp$ and consider the delta distribution 
$$T:=\sum_{g\in\Gamma}(\pi^{\sigma,\lambda}(g)\otimes\vp(g))(\delta v\otimes w)\in {}^\Gamma C^{-\infty}(\Omega,V(\sigma_\lambda,\vp))\cong  C^{-\infty}(B,V_B(\sigma_{\lambda},\vp))\ .$$ 
Since the singular
parts of $ext$ are finite-dimensional we find a smooth section $\phi\in C^{\infty}(B,V_B(\sigma_{\lambda},\vp))$ such that $f:=ext(T-\phi)\in {}^\Gamma C^{-\infty}(\partial X,V(\sigma_{\lambda},\vp))$ is defined. 
We decompose $\delta v\otimes w=f-(f-\delta v\otimes w)$.
Since $x_0\not\in \singsupp(f-\delta v\otimes w)$ and $x_0\not\in \singsupp A_{\Gamma,\vp}(f)$ we conlude that $A(\delta v)$ is smooth at $x_0$,
and hence $A(\delta v)\in {}^N C^\infty(\partial X,V(\tau_\mu))$.

Now 
$$ {}^N C^\infty(\partial X,V(\tau_\mu))=H^0(\naaa,C^\infty(\partial X,V(\tau_\mu)))\cong H^0(\naaa,C^\infty(\partial X,V(\tau_\mu))_K), $$
where the
second equality is a special case of Casselman's comparison
theorem for $\naaa$-cohomology (see e.g. \cite{bunkeolbrich961} or \cite{hechttaylor98}).
Thus $A(\delta v)$ is $K$-finite.
There exists a finite-dimensional subspace
$E\subset \cU(\gaaa)$ such that $\cU(\gaaa)=\cU(\kaaa) \cZ E\,\cU(\naaa)$
(see \cite{wallach88}, 3.7.1).
Thus the space $Z_K:=\cU(\gaaa)A(\delta v)=\cU(\kaaa)E A(\delta v)$
is finite-dimensional.

Since $\delta v$ generates $C^{-\infty}(\partial X,V(\sigma_\lambda))$,
the element $A(\delta v)$ generates $\im(A)$.
We conclude that $\im(A)_K=Z_K$ and thus $\dim\:\im(A)<\infty$.
This finishes the proof of the lemma.
\hB
 
\begin{kor}\label{leda}
Let $(\vp,V_\vp)$ be a finite-dimensional (not necessarily unitary) representation of $\Gamma$, and let $V_\pi$ be an irreducible admissible representation of $G$.
Then ${}^\Gamma(V_{\pi,-\infty}\otimes V_\vp)$ is infinite-dimensional unless $V_\pi$ is
finite-dimensional.
\end{kor} 
\proof 
By Casselman's subrepresentation theorem (\cite{wallach88}, 3.8.3.) in
conjunction with the functorial properties of the smooth globalization (\cite{wallach92}, 11.6.7.) we find elements 
$(\sigma,\lambda), (\tau,\mu)\in\hat M\times\aca$ such that $V_{\pi,-\infty}$ is
a quotient of $C^{-\infty}(\partial X,V(\sigma_{\lambda}))$ and a submodule
of $C^{-\infty}(\partial X,V(\tau_{\mu}))$. Thus there is a non-trivial $G$-intertwining operator
$$A:C^{-\infty}(\partial X,V(\sigma_{\lambda}))\rightarrow C^{-\infty}(\partial X,V(\tau_{\mu}))$$
satisfying $\im(A)\cong V_{\pi,-\infty}$. If $\dim(V_\pi)=\infty$, then by Lemma \ref{skylla} the subspace
$\im(A_{\Gamma,\vp})\subset {}^\Gamma(V_{\pi,-\infty}\otimes V_\vp)$
is infinite-dimensional. 
 \hB

The following proposition completes our description of the tempered and square-integrable invariant distribution vectors.

\begin{prop}\label{kles}
Let $(\pi,V_\pi)\in\hat G_d$. Then for any unitary representation $(\vp,V_\vp)$ of $\Gamma$ the space ${}^\Gamma(V_{\pi,-\infty}\otimes V_\vp)_d$ is  infinite-dimensional.
\end{prop}
\proof
Let $\gamma\in\hat K$ be the minimal $K$-type
of $V_\pi$. Casselman's subrepresentation theorem provides an embedding
$$\beta: V_{\pi,\infty}\rightarrow C^{\infty}(\partial X,V(\sigma_{-\lambda}))\ ,$$
where $\Hom_M(V_\sigma,V_\gamma)\not=0$, and $-(\lambda+\rho)$ is the leading exponent in the asymptotic
expansion for $a\to\infty$ of the matrix coefficients $c_{v,\tilde v}$, $v\in V_{\pi,\infty}$, $\tilde v\in V_{\pi^\prime,K}$. Since $V_\pi$ is a discrete
series representation we have $0<\lambda\in\aaaa^*$. Forming the adjoint with
respect to hermitian scalar products we obtain a projection
$$ q: C^{-\infty}(\partial X,V(\sigma_{\lambda}))\rightarrow
V_{\pi,-\infty}\ .$$
By functoriality we can extend $\beta$ to a map between the corresponding spaces of distribution
vectors and obtain a $G$-intertwining operator
$$A:=\beta\circ q: C^{-\infty}(\partial X,V(\sigma_{\lambda}))\rightarrow
C^{-\infty}(\partial X,V(\sigma_{-\lambda}))$$
satisfying $\im(A)\cong V_{\pi,-\infty}$.
As in Lemma \ref{skylla} we consider the operator
$$A_{\Gamma,\vp}:{}^\Gamma C^{-\infty}(\partial X,V(\sigma_{\lambda},\vp))\rightarrow {}^\Gamma C^{-\infty}(\partial X,V(\sigma_{-\lambda},\vp))\ .$$
This lemma combined with Lemma \ref{ddeenn} tells us that
$A_{\Gamma,\vp}\left({}^\Gamma C^{-\infty}_\Omega(\partial X,V(\sigma_{\lambda},\vp))\right)$
is infinite-dimensional. Hence $Z:=(q\otimes\id)\left({}^\Gamma C^{-\infty}_\Omega(\partial X,V(\sigma_{\lambda},\vp))\right)\subset{}^\Gamma(V_{\pi,-\infty}\otimes V_\vp)$ is infinite-dimensional, too. It remains to show that $Z\subset  
{}^\Gamma(V_{\pi,-\infty}\otimes V_\vp)_d$.

Choose an embedding $t\in\Hom_K(V_{\tilde\gamma},V_{\pi^\prime})$, and define $T\in\Hom_M(V_\sigma,V_\gamma)$ by 
$\langle T(w), v\rangle:=\langle w, [{}^tq\circ t( v)](e)\rangle$
for all $w\in V_\sigma$, $v\in V_{\tilde\gamma}$.
Recall the definition of $v_T$ from  Lemma  \ref{circe} and observe that $v_T={}^tq(t(v))$.
We consider the
Poisson transform
$$ P:=P^T_\lambda\otimes\id: C^{-\infty}(\partial X,V(\sigma_{\lambda},\vp))
\rightarrow C^\infty(X,V(\gamma,\vp))\ .$$
By Lemma \ref{circe} we find for all $\phi\in C^{-\infty}(\partial X,V(\sigma_{\lambda},\vp))$, $v\in V_{\tilde\gamma}$, $g\in G$
\begin{equation}\label{eos}
\langle P\phi(g),v\rangle=c_{\phi,v_T}(g)=c_{(q\otimes\id)(\phi),t(v)}(g)\in V_\vp\ .
\end{equation}
Since $t(v)$ generates $V_{\pi^\prime,K}$ by Lemma \ref{cassel} it suffices 
to show that $P\phi\in L^2(Y,V_Y(\gamma,\vp))$ for all
$\phi\in {}^\Gamma C^{-\infty}_\Omega(\partial X,V(\sigma_{\lambda},\vp))$.

Let $D$ be a closed neighbourhood of $\clo(F)\cap \Omega$ for
some fundamental domain $F\subset X$ of $\Gamma$, and let $Q\subset \partial X\setminus D$ be a closed neighbourhood of $\Lambda$.
Let $\chi\in C^\infty(\partial X)$ be a cut-off function with $\supp(\chi)
\subset Q$, $\supp(1-\chi)\cap\Lambda=\emptyset$.
Let $\phi\in {}^\Gamma C^{-\infty}_\Omega(\partial X,V(\sigma_{\lambda},\vp))$.
By Lemma \ref{poisasm}, 3., there exists 
a constant $C$ such that for $a\gg 0$, $k\in DM$
\begin{equation}\label{old} 
|P(\chi\phi)(ka)|\le C a^{-(\lambda+\rho)} \ .
\end{equation}
On the other hand $(1-\chi)\phi\in C^{\infty}(\partial X,V(\sigma_{\lambda},\vp))$. 
By (\ref{eos}) the function $P((1-\chi)\phi)$ can be expressed
in terms of matrix coefficients $c_{s,t}$
 of the  discrete series representation $V_\pi$, where
$s\in V_{\pi,\infty}$ and $t\in V_{\pi^\prime,K}$.
 Thus $P((1-\chi)\phi)$ has an asymptotic expansion with leading
exponent $-(\lambda+\rho)$. In particular, it satisfies an estimate of the form
(\ref{old}), too. We conclude that $P\phi\in L^2(Y,V_Y(\gamma,\vp))$. This finishes the proof of the proposition.
\hB

\section{Eisenstein series, wave packets, and scalar products}\label{wxa}

In this section we consider the Eisenstein series and the wave-packet transform. For a moment we can drop the unitarity condition on $\vp$. Let
$\gamma$ be a finite-dimensional unitary representation of $K$, 
$\sigma\in \hat{M}$,
$T\in\Hom_M(V_\sigma,V_\gamma)$, and  let $P^T_\lambda$, $\lambda\in\aca$,
be the associated Poisson transform (see Definition \ref{defofpoi}). 
\begin{ddd}
For $\phi\in C^{-\infty}(B,V_B(\sigma_\lambda,\vp))$ we define the Eisenstein series
$E(\lambda,\phi,T)\in C^\infty(Y,V_Y(\gamma,\vp))$ by  
$$E(\lambda,\phi,T):= P_\lambda^T\circ ext(\phi) \ .$$
\end{ddd}
The Eisenstein series $E(\lambda,\phi,T)$ is an eigenvector
of $\cZ$ for the infinitesimal character of 
the principal series representation $\pi^{\sigma,\lambda}$. 
Theorem \ref{part1} and the functional equation of the Poisson transform
(\ref{mi9}) have the following immediate corollary.
\begin{kor}\label{funeq}
The Eisenstein series gives rise to a meromorphic family defined on $\aca$ 
(or $\{\Ree(\lambda)>\delta_\Gamma+\delta_\vp\}$ in case that $X=\OO H^2$)
of continuous maps
$$E(\lambda,.,T):C^{-\infty}(B,V_B(\sigma_\lambda,\vp))\rightarrow C^\infty(Y,V_Y(\gamma,\vp))$$
with finite-dimensional singularities. It satisfies the functional equation 
\begin{equation}\label{cisss}
E(\lambda,\hat S_{-\lambda}\phi,T)=
E(-\lambda ,\phi,(c_{ \gamma}(\lambda)T)^w)\ .\end{equation}
\end{kor}
From Lemma \ref{poisasm}, 3., we gain detailed knowledge of the asymptotics
of $E(\lambda,\phi,T)(y)$ for $y\to b\in B$ (see the discussion of (\ref{diable})
in the proof of Proposition \ref{scalar} below).

We now return to our unitarity assumption on $\vp$.  
\begin{kor}\label{gogo}
\noindent
\begin{enumerate}\item
The Eisenstein series is regular on $\{\Ree(\lambda)=0,\lambda\not=0\}$.
\item
If $\Ree(\lambda)=0$, $\lambda\not=0$, then 
$E(\lambda,.,T)$ maps $C^{-\infty}(B,V_B(\sigma_\lambda,\vp))$
continuously to $\cC^\prime(Y,V_Y(\gamma,\vp))$.
\item
In the half plane  $\{\Ree(\mu)\ge 0\}$
the Eisenstein series has at most first-order poles which are located in the finite set $PS_\res(\sigma,\vp)\subset [0,\delta_\Gamma]$.
The residue at $\lambda\in PS_\res(\sigma,\vp)\setminus \{0\}$
$$\res_{\mu=\lambda} E(\mu,.,T)$$
maps $C^{-\infty}(B,V_B(\sigma_\lambda,\vp))$
to $L^2(Y,V_Y(\gamma,\vp))$.
\end{enumerate}
\end{kor}
\proof
1. follows from Lemma \ref{extregatim}.
To see 2. note that $\pi^{\sigma,\lambda}$ is
tempered if $\Ree(\lambda)=0$, that $E(\lambda,\phi,T)$
can be expressed in terms of matrix coefficients of $\pi^{\sigma,\lambda}$ (Lemma \ref{circe}),
and apply \ref{zeus}.
We have \begin{equation}
\label{formelop}
\res_{\mu=\lambda} E(\mu,.,T)(\phi) = P^T_\lambda\circ (\res_{\lambda} ext)(\phi)\ .\end{equation}
3. now follows from Propositions \ref{firtu}, \ref{firtu0}, and \ref{upperbound}.
\hB

By (\ref{formelop}) the following proposition can also
be considered as the determination of 
the $L^2$-scalar product between residues of Eisenstein series.

We introduce a scalar product 
on $\Hom_M(V_\sigma,V_\gamma)$ by $(T_1,T_2)\id_{V_\sigma}:=T_2^*T_1$.
This makes sense because of our standing assumption $\sigma\in\hat M$.   
By $(.,.)_B$ we denote the natural sesquilinear pairing between $C^{-\infty}(B,V_B(\sigma_\lambda,\vp))$ and  $C^{\infty}(B,V_B(\sigma_{-\bar\lambda},\vp))$.

\begin{prop}\label{volume}
Let $\lambda\in PS_\res(\sigma,\vp)\setminus\{0\}$, let $\gamma$ be a finite-dimensional representation of $K$, $T_1,T_2\in \Hom_M(V_\sigma,V_\gamma)$,
$\phi_1=(\res_\lambda ext)(f)\in E_\Lambda(\sigma_\lambda,\varphi)$ for some
$f\in C^{-\infty}(B,V_B(\sigma_\lambda,\vp))$, and $\phi_2\in {}^\Gamma C^{-\infty}(\Lambda,V(\sigma_\lambda,\vp))$. Then
$$ (P^{T_1}_\lambda\phi_1,P^{T_2}_\lambda\phi_2)_{L^2(Y,V_Y(\gamma,\vp))}
=\omega_X \:
(c_\gamma(\lambda)T_1,T_2^w)\:
(f,res\circ\hat J_\lambda\phi_2)_B\ .$$
Here $\omega_X=\frac{\omega_n}{2^{r}\alpha^{n-1}}$, where $n=\dim X$, 
$\omega_n=\vol(S^{n-1})=\frac{2\pi^\frac{n}{2}}{\Gamma(\frac{n}{2})}$ and
$r\in\nat$ is such that $r\alpha=2\rho$.
If $\gamma$ is the minimal $K$-type of $\pi^{\sigma,\lambda}$, then
$$(c_\gamma(\lambda)T_1,T_2^w)=c_\sigma(\lambda)(T_1,T_2)\ .$$
\end{prop}
\proof
First we can assume that $f$ is smooth. We extend $f$ to a holomorphic family
$\mu\mapsto f_\mu\in C^{\infty}(B,V_B(\sigma_\mu,\vp))$ defined in a neighbourhood of $\lambda$. Let $B_R\subset X$ be the ball of radius $R$,
and let $(. , .)_{B_R}$ denote the scalar product in $L^2(B_R,V(\gamma,\vp))$.
Let $\chi\in C^{\infty}(X)$ be the cut-off function constructed in Lemma
\ref{lll}. Set 
$$ S_R(\mu):=(\mu-\lambda)\: (P^{T_1}_\mu ext(f_\mu) , \chi P^{T_2}_\lambda\phi_2)_{B_R}\ .$$
Then
$$ (P^{T_1}_\lambda\phi_1,P^{T_2}_\lambda\phi_2)_{L^2(Y,V_Y(\gamma,\vp))}=\lim_{R\to\infty}\lim_{\mu\to\lambda} S_R(\mu) \ .$$
Let $D:=-\Omega_G+c(\sigma)+\lambda^2$ be the shifted Casimir operator. As in the proof of Proposition \ref{firtu} we obtain 
$$ S_R(\mu)=-\frac{1}{\lambda+\mu}\: 
(D P^{T_1}_\mu ext(f_\mu), \chi P^{T_2}_\lambda\phi_2)_{B_R}\ .$$
While $P^{T_1}_\mu ext(f_\mu)$ has a first order pole at $\mu=\lambda$,
$D P^{T_1}_\mu ext(f_\mu)$ is regular, and its
value at $\mu=\lambda$ equals $D F$, where
$F\in {}^\Gamma C^\infty(X,V(\gamma,\vp))$ is the constant term of the Laurent expansion of $P^{T_1}_\mu ext(f_\mu)$
at $\lambda$.

Note that $D=\nabla^*\nabla + \cR$ for some selfadjoint endomorphism $\cR$ 
of $V(\gamma,\vp)$, where $\nabla^*\nabla$ is the Bochner Laplacian associated to the 
invariant connection $\nabla$ of $V(\gamma,\vp)$. Thus we can apply Green's formula
in the same spirit as in the proof of Proposition \ref{green}:
\begin{eqnarray}\label{hubert}
\lim_{\mu\to\lambda} S_R(\mu)&=&
-\frac{1}{2\lambda}\left(
(DF,\chi P^{T_2}_\lambda\phi_2)_{B_R}-(F,D\chi P^{T_2}_\lambda\phi_2)_{B_R}
+(F,[D,\chi] P^{T_2}_\lambda\phi_2)_{B_R} \right)\nonumber\\
&=&\frac{1}{2\lambda}\left( 
(\nabla_n F,\chi P^{T_2}_\lambda\phi_2)_{\partial B_R}
-(F,\nabla_n\chi P^{T_2}_\lambda\phi_2)_{\partial B_R}
-(F,[D,\chi] P^{T_2}_\lambda\phi_2)_{B_R} \right) .
\end{eqnarray}
Applying Cauchy's integral formula to the asymptotic expansion of $P^{T_1}_\mu ext(f_\mu)$ given in Lemma \ref{poisasm}, 3. we find some $\epsilon>0$ such
that for $a\to\infty$
$$F(ka)=a^{\lambda-\rho} c_\gamma(\lambda)T_1 f(k)+O(a^{\lambda-\rho-\epsilon})$$
uniformly as $kM$ varies in compact subsets of $\Omega$.

As in the proof of Proposition \ref{green} we obtain 
$(F,[D,\chi] P^{T_2}_\lambda\phi_2)_{X}=0$. In order to perform the limit $R\to\infty$
in (\ref{hubert}) we use $\lim_{R\to\infty}\ee^{-2\rho R}\vol(\partial B_{R})=\omega_X$ and the asymptotic expansions of $F$,
$P^{T_2}_\lambda\phi_2$, and obtain 
\begin{eqnarray*}
(P^{T_1}_\lambda\phi_1,P^{T_2}_\lambda\phi_2)_{L^2(Y,V_Y(\gamma,\vp))}&=& 
\frac{\omega_X}{2\lambda}\left(
(\lambda-\rho) \int_{\partial X} (c_\gamma(\lambda)T_1 f(k), \chi_\infty(k)
T_2^w (\hat J_\lambda \phi_2)(k))\: dk \right.\\
&&\quad\quad\left. +(\lambda+\rho) \int_{\partial X} (c_\gamma(\lambda)T_1 f(k), \chi_\infty(k)
T_2^w (\hat J_\lambda \phi_2)(k))\: dk\right)\\
&=&\omega_X 
(T_2^{w*}c_\gamma(\lambda)T_1 f,res\circ\hat J_\lambda\phi_2)_B\ .
\end{eqnarray*}
This finishes the proof of the proposition.
\hB

On $I^{\sigma,\lambda}$ we consider consider the scalar product (\ref{jam}).  
Let $(.,.)$ be the scalar product on $E_\Lambda(\sigma_\lambda,\vp)\oplus U_\Lambda(\sigma_\lambda,\vp)={}^\Gamma (I^{\sigma,\lambda}_{-\infty}\otimes V_\vp)_d$ induced by the matrix coefficient map (see Proposition \ref{nonte}
and Corollary \ref{abel}).
\begin{kor}\label{otto1}
The decomposition $E_\Lambda(\sigma_\lambda,\vp)\oplus U_\Lambda(\sigma_\lambda,\vp)$ is orthogonal with respect to $(.,.)$.
If $\phi_1,\phi_2\in E_\Lambda(\sigma_\lambda,\vp)$, then
$$(\phi_1,\phi_2)=\frac{\omega_X c_\sigma(\lambda)}{\dim(V_\sigma)} (f,res\circ \hat{J}_\lambda \phi_2)\ ,$$
where  $f\in C^{-\infty}(B,V_B(\sigma_\lambda,\vp))$ such that $\phi_1=(\res_\lambda ext)(f)$.
\end{kor}
\proof
Let $\gamma$ be the minimal $K$-type of $\pi^{\sigma,\lambda}$.
For $T\in \Hom_M(V_\sigma,V_\gamma)$ and $v\in V_{\tilde{\gamma}}$  
let $v_T\in C^\infty(\partial X,V(\tilde{\sigma}_{\pm\lambda}))$
be given by $v_T(k)={}^t T\tilde{\gamma}(k^{-1}) v$.
We have $J_{\tilde{\sigma},-\lambda} v_T=v_T$. In fact,
$\sigma$ is Weyl-invariant and thus for all $\phi\in C^{-\infty}(\partial X,V(\sigma_{-\lambda}))$
\begin{eqnarray*}
\langle \phi,J_{\tilde{\sigma},-\lambda} v_T\rangle &=&
\langle J_{\sigma,-\lambda}\phi, v_T\rangle\\
&=&
\langle \hat{J}_{\sigma,-\lambda}c_\sigma(\lambda)^{-1}\phi, v_T\rangle\\
&\stackrel{\ref{circe}}{=}&c_\sigma(\lambda)^{-1}\langle  P^T_\lambda \circ \hat{J}_{\sigma,-\lambda}(\phi)(1), v \rangle\\
&\stackrel{(\ref{mi9})}{=}&c_\sigma(\lambda)^{-1}\langle  P^{(c_\gamma(\lambda)T)^w}_{-\lambda}(\phi)(1), v \rangle\\ 
&\stackrel{(\ref{mao})}{=}& \langle  P^{T}_{-\lambda}(\phi)(1), v \rangle\\
&=& \langle  \phi , v_T \rangle\ .
\end{eqnarray*}
By (\ref{jam}) we have
\begin{eqnarray}
([v_T],[v_T])&=&(v_T,J_{\tilde{\sigma},-\lambda} v_T)_0\nonumber\\
&=&\int_{K} \|v_T(k)\|^2 dk\nonumber\\
&=&\int_K ({}^t T \tilde{\gamma}(k^{-1})v,{}^t T \tilde{\gamma}(k^{-1})v) dk\nonumber\\
&=&\int_K (\tilde{\gamma}(k)({}^tT)^*{}^t T \tilde{\gamma}(k^{-1})v, v) dk\nonumber\\
&=&\frac{\dim(V_{\tilde{\sigma}})}{\dim(V_{\tilde{\gamma}})} \|v\|^2 \|T\|^2\ .\label{u8u8u}
\end{eqnarray}
Let $\{v^i\}_{i=1}^{\dim(V_{\tilde{\gamma}})}$ be an othonormal base of $V_{\tilde{\gamma}}$.
We compute using Proposition \ref{volume}, Lemma \ref{circe} 
\begin{eqnarray*}
\frac{\omega_X c_\sigma(\lambda)}{\dim(V_\sigma)}(f,res\circ \hat{J}_\lambda \phi_2)
&=&\frac{1}{\|T\|^2\dim(V_\sigma)} (P^T_\lambda \phi_1,P^T_\lambda \phi_2)_{L^2(Y,V_Y(\gamma,\vp))}\\
&=&\frac{1}{\|T\|^2\dim(V_\sigma)}\sum_{i=1}^{\dim(V_\gamma)} (c_{{\phi_1},v_T^i},c_{{\phi_2},v_T^i})_{L^2(\Gamma\backslash G,\vp)}\\
&=&\frac{1}{\|T\|^2\dim(V_\sigma)} \sum_{i=1}^{\dim(V_\gamma)} (\phi_1,\phi_2) \|v_T^i\|^2\\
&=&(\phi_1,\phi_2)\ .
\end{eqnarray*}
This proves the corollary.\hB

Now we turn to the definition of the wave packet transform.
Roughly speaking, a wave packet of Eisenstein series is
an average of the Eisenstein series over imaginary parameters
with, say, a smooth, compactly supported weight function with respect
to the expected Plancherel measure $p_\sigma(\ii\mu)d\mu$.
More precisely, let $\aaaa^*_+:=\{\lambda\in\aaaa^*\:|\:\langle\lambda,\alpha\rangle >0\}$ be the open positive chamber in $\aaaa^*$. Then the space of such weight functions
$\cH^\sigma_0(\vp)$ is the linear space of smooth families $\aaaa^*_+\ni \mu\mapsto \phi_{\imath\mu}\in C^\infty(B,V_B(\sigma_{\imath\mu},\vp))$ with compact support in $\aaaa^*_+$
with respect to $\mu$. 
Because of the functional equation  (\ref{cisss}) it will be sufficient to consider wave packets on the positive imaginary axis, only.
  
Let $\gamma$ be a finite-dimensional unitary representation of $K$.
We first define the wave packet transform on $\cH^\sigma_0(\vp)\otimes \Hom_M(V_\sigma,V_\gamma)$.
Later we will extend it by continuity to a Hilbert space closure.

\begin{ddd}\label{now}
The wave packet transform  is the  map 
$$E:\cH^\sigma_0(\vp)\otimes \Hom_M(V_\sigma,V_\gamma)\rightarrow C^\infty(Y,V_Y(\gamma,\vp))$$
given by 
$$E(\phi\otimes T) := E(\phi,T):= \int_{\aaaa^*_+} E(\imath\mu,\phi_{\imath\mu},T)p_\sigma(\ii\mu)\:d\mu\ ,$$
where $d\mu$ is the Lebesgue measure on
$\aaaa^*_+\cong (0,\infty)$.
The section $E(\phi,T)$, $\phi\in\cH^\sigma_0(\vp)$, $T\in\Hom_M(V_\sigma,V_\gamma)$ is called a wave packet (of Eisenstein series). 
\end{ddd}

\begin{lem}\label{swa}
If $T\in\Hom_M(V_\sigma,V_\gamma)$, $\phi\in \cH^\sigma_0(\vp)$, then $E(\phi,T)\in \cC(Y,V_Y(\gamma,\vp))$.
\end{lem}
\proof
Set $\psi_{\imath\mu}:=ext\:\phi_{\imath\mu}$ and define
$$P(\psi):=\int_{\aaaa^*_+} P^T_{\imath\mu}(\psi_{\imath\mu})p_\sigma(\ii\mu)\: d\mu\ .$$
Let $\chi$ be the cut-off function constructed in Lemma \ref{lll}.
In view of Definition \ref{sscchh} of the Schwartz space  we have to show that $\chi P(\psi)\in \cC(X,V(\gamma,\vp))$.
Let $\chi_0\in C_c^\infty(X)$  be some cut-off function which is equal to $1$
on some neighbourhood of  $eK\in X$. Obviously, we have $\chi\chi_0 P(\psi)\in C_c^\infty(X,V(\gamma,\vp))\subset\cC(X,V(\gamma,\vp))$. It remains to
show that 
$$\chi_1 P(\psi)\in \cC(X,V(\gamma,\vp))\ ,$$
where $\chi_1=\chi(1-\chi_0)$.
Observe that the seminorms $q_{D,N}$, $D\in\cU(\gaaa)$, $N\in N_0$ 
defined by
$$ q_{D,N}(f)^2:=\int_G |\log(a_g)^Nf(Dg)|^2\:dg $$
are sufficient in order to define the topology on $\cC(X,V(\gamma,\vp))$.
We thus have to show that
$g\mapsto \log(a_g)^N(\chi_1P(\psi))(Dg)$ is square-integrable.
Let $\Delta:\cU(\gaaa)\rightarrow \cU(\gaaa)\otimes\cU(\gaaa)$
denote the coproduct and write $\Delta(D)=\sum_{\alpha} D_\alpha\otimes D^\prime_\alpha$ for some $D_\alpha,D^\prime_\alpha\in\cU(\gaaa)$.
Then $(\chi_1P(\psi))(Dg)=\sum_{\alpha} \chi_1(D_\alpha g)P(\psi)(D^\prime_\alpha g)$. By Lemma \ref{lll}, 4., $g\mapsto \chi_1(D g)$ is bounded for any
$D\in\cU(\gaaa)$. It remains to show for any $D\in\cU(\gaaa)$ and cut-off function $\chi_2\in C^\infty(X)$ having compact support in $(X \cup\Omega)\setminus\{eK\}$ the function
\begin{equation}\label{hexe}
\log(a_g)^N\chi_2(g) P(\psi) (Dg)
=\log(a_g)^N\chi_2(ka)
\int_{\aaaa^*_+} P^T_{\imath\mu}(\pi^{\sigma,\imath\mu}(D)\psi_{\imath\mu})(g)p_\sigma(\ii\mu)\: d\mu
\end{equation}
is square integrable.
By Lemma \ref{poisasm}, 3., there exists $\epsilon>0$ such that for $D\in \cU(\gaaa)$, $ka\in \supp(\chi_2)$
\begin{eqnarray}\label{diable}
P^T_{\imath\mu}(\pi^{\sigma,\imath\mu}(D)\psi_{\imath\mu})(ka) 
&=& a^{\imath\mu-\rho}c_\gamma(\imath\mu)T\psi_{\imath\mu}(k)\nonumber\\
&&+  a^{-\imath\mu-\rho}T^w(\hat{J}_{\imath\mu}(\pi^{\sigma,\imath\mu}(D) \psi_{\imath\mu}))(k)\\
&&+ a^{-(\rho+\epsilon)} R(\imath\mu,\pi^{\sigma,\imath\mu}(D)\psi_{\imath\mu},ka)\ ,\nonumber
\end{eqnarray}
where the remainder function $(\mu,ka)\mapsto R(\imath\mu,\pi^{\sigma,\imath\mu}(D)\psi_{\imath\mu},ka)$ 
is uniformly bounded.

Since the families $\psi_{\imath\mu}$ and $\hat{J}_{\imath\mu}(\pi^{\sigma,\imath\mu}(D) \psi_{\imath\mu})$ have compact support
with respect to $\mu$ and  are smooth in $(\mu,k)$ (as long $kM\in\Omega$) by Lemma
\ref{off} each summand of (\ref{diable}) contributes to the function
(\ref{hexe}) a summand which is bounded by $C_{N^\prime} (1+\log(a_g))^{-N^\prime} a_g^{-\rho} $ for any $N^\prime\in\nat$.
It follows that  the function (\ref{hexe}) is square integrable. This implies
the lemma.
\hB

Let $\sigma^i \in\hat{M}$, $T_i\in \Hom_M(V_{\sigma^i},V_\gamma)$, 
$i=1,2$,  $\lambda\in\aaaa^*_+$, $\phi\in \cH_0^{\sigma^2}(\vp)$, and $\psi\in C^{-\infty}(B,V_B({\sigma^1}_{\imath\lambda},\vp))$.
By Lemma \ref{swa} and Corollary \ref{gogo}
the pairing $(E(\ii\lambda,\psi,T_1),E(\phi,T_2))_{L^2(Y,V_Y(\gamma,\vp))}$
between the Eisenstein series $E(\ii\lambda,\psi,T_1)$ and the wave packet $E(\phi,T_2)$ is well-defined.
The following proposition gives an explicit formula for this pairing and
is the crucial step in the determination of the absolute continuous
part of the Plancherel measure.

\begin{prop}\label{scalar}
We have 
$$(E(\ii\lambda,\psi,T_1), E(\phi,T_2))_{L^2(Y,V_Y(\gamma,\vp))}=\left\{
\begin{array}{cc} 2\pi\omega_X \: (T_1,T_2)\: (\psi,\phi_{\imath\lambda})_B & \sigma^1=\sigma^2 \\
0 & \sigma^1\not=\sigma^2
\end{array}\right.\ .$$
\end{prop}
\proof
For fixed $\psi$ and $T_i$, $i=1,2$, we consider the continuous linear functional
(which is in fact a distribution section of a bundle over $\imath\aaaa_+^*\times B$)
$$R:\cH^{\sigma^2}_0(\vp)\rightarrow \C$$ given by
$R(\phi):=(E(\ii\lambda,\psi,T_1), E(\phi,T_2))_{L^2(Y,V_Y(\gamma,\vp))}$.
If $D\in\cZ$ and $\sigma\in\hat{M}$ we consider the polynomial 
$\chi_{\mu_\sigma+\rho_m-\mu}(D)$ on $\aca$. Since $\chi_{\mu_\sigma+\rho_m-\mu}$ is the infinitesimal character
of $\pi^{\sigma,\mu}$ we have
\begin{eqnarray*}
0&=& ((D^*-\chi_{\mu_{\sigma^1}+\rho_m-\imath\lambda}(D^*))E(\ii\lambda,\psi,T_1), E(\phi,T_2))_{L^2(Y,V_Y(\gamma,\vp))}\\
&=& (E(\ii\lambda,\psi,T_1), (D-\bar{\chi}_{\mu_{\sigma^1}+\rho_m-\imath\lambda}(D^*))
E(\phi,T_2))_{L^2(Y,V_Y(\gamma,\vp))}\\
&=&(E(\ii\lambda,\psi,T_1), E(\tilde{\phi},T_2))_{L^2(Y,V_Y(\gamma,\vp))}\ ,
\end{eqnarray*}
where $\tilde{\phi}_{\imath\mu}=(\chi_{\mu_{\sigma^2}+\rho_m-\imath\mu} (D) -\chi_{\mu_{\sigma^1}+\rho_m-\imath\lambda} (D)) \phi_{\imath\mu}$.
We conclude that multiplication by the polynomial 
$(\chi_{\mu_{\sigma^2}+\rho_m-\imath\mu} (D) -\chi_{\mu_{\sigma^1}+\rho_m-\imath\lambda} (D))$
annihilates the functional $R$. Thus $R$ is supported on the zero set of this
polynomial. If $\sigma^1\not=\sigma^2$, then
$$\bigcap_{D\in \cZ} \{\mu\in\aaaa_+^*\:|\:(\chi_{\mu_{\sigma^2}+\rho_m-\imath\mu} (D) -\chi_{\mu_{\sigma^1}+\rho_m-\imath\lambda} (D))=0\} =\emptyset$$
and therefore $R=0$. This proves the proposition in case $\sigma^1\not=\sigma^2$.

Assume now that $\sigma^1=\sigma^2=:\sigma$.
Observe that
$$\bigcap_{D\in \cZ} \{\mu\in\aaaa_+^*\:|\:(\chi_{\mu_{\sigma }+\rho_m-\imath\mu} (D) -\chi_{\mu_{\sigma }+\rho_m-\imath\lambda} (D))=0\} =\{ \lambda\}$$
and that the functional $R$ is of the
form $R(\phi)= r(\phi_{\imath\lambda})$ for some $r\in C^{-\infty}(B,V_B({\sigma }_{\imath\lambda},\vp))$,
where $r$ remains to be determined.
We prefer to give a direct proof of the proposition in case $\sigma^1=\sigma^2$, which  
does not refer to this observation.

Because of the continuity of $E(\ii\lambda,.,T_1)$ (see Corollary \ref{gogo}, 2.)
we can assume that 
$\psi\in  C^{\infty}(B,V_B({\sigma}_{\imath\lambda},\vp))$.  
We apply Green's formula in a similar way as in the proofs of Propositions \ref{green} and \ref{volume}.

Let $\chi$ be the cut-off function as constructed in Lemma \ref{lll} and $B_R$ the ball of
radius $R$ around the origin of $X$.
For $\mu\in\aaaa^*_+$ 
we consider
$$S_R(\mu):=(E(\imath\lambda,\psi,T_1),
\chi E(\imath\mu,\phi_{\imath\mu},T_1) )_{B_R}\ .$$
We have thus to compute
$$ \lim_{R\to\infty} \int_0^\infty S_R(\mu) p_\sigma(\ii\mu)\: d\mu\ .$$ 
If we set $A:=-\Omega_G+c(\sigma)$, then we obtain
\begin{eqnarray}
(\lambda^2-\mu^2)S_R(\mu)
&=& 
(A\:E(\imath\lambda,\psi,T_1),\chi E(\imath\mu,\phi_{\imath\mu},T_2))_{B_R}
-(E(\imath\lambda,\psi,T_1),A\chi E(\imath\mu,\phi_{\imath\mu},T_2))_{B_R}
\nonumber\\
&&+ (E(\imath\lambda,\psi,T_1),
[A,\chi] E(\imath\mu,\phi_{\imath\mu},T_2))_{B_R}
\nonumber\\
&=&-(\nabla_n E(\imath\lambda,\psi,T_1), 
\chi E(\imath\mu,\phi_{\imath\mu},T_2))_{\partial B_R}\label{ad1}\\
&&+(E(\imath\lambda,\psi,T_1) , 
\nabla_n\chi E(\imath\mu,\phi_{\imath\mu},T_2) )_{\partial B_R} \label{ad2}\\
&&+ ([A,\chi] E(\imath\mu,\phi_{\imath\mu},T_1), E(\imath\lambda,\psi,T_2))_{B_R}\ .
\label{ad3}
\end{eqnarray}
We now apply the asymptotic expansion (\ref{diable}) which holds on the support of $\chi$  in case that  $D=1$, and for 
$ext(\phi_{\imath\mu})$ and $ext(\psi)$ in place of $\psi_{\imath\mu}$. 
Note that for $k\in\Omega M$ we have $ext(\phi_{\imath\mu})(k)=\phi_{\imath\mu}(k)$, $\hat J_{\ii\mu}(ext(\phi_{\imath\mu}))(k)=\hat S_{\ii\mu}(\phi_{\imath\mu})(k)$, 
etc. 
We obtain for large $R$ with $a_R:=\ee^{R}$, $\omega_R:=a_R^{-2\rho}\vol(\partial B_R)$
\begin{eqnarray*}
\frac{1}{\omega_R}((\ref{ad1})+(\ref{ad2}))
&=&
\ii(-\lambda-\mu) a_R^{\ii(\lambda-\mu)}
\int_{K} \chi_\infty(k) (c_\gamma(\ii\lambda)T_1 \psi)(k),
c_\gamma(\ii\mu)T_2 \phi_{\imath\mu}(k))\: dk\\
&&
+\ii(-\lambda+\mu) a_R^{\ii(\lambda+\mu)}
\int_{K} \chi_\infty(k) (c_\gamma(\ii\lambda)T_1 \psi(k),
T_2^w(\hat S_{\ii\mu}\phi_{\imath\mu})(k))\: dk\\
&&
+\ii(\lambda-\mu) a_R^{\ii(-\lambda-\mu)}
\int_{K} \chi_\infty(k)  (T_1^w\hat S_{\ii\lambda} \psi)(k),
c_\gamma(\ii\mu)T_2 \phi_{\imath\mu}(k))\: dk\\
&&
+\ii(\lambda+\mu) a_R^{\ii(-\lambda+\mu)}
\int_{K} \chi_\infty(k) (T_1^w (\hat S_{\ii\lambda} \psi)(k),
T_2^w(\hat S_{\ii\mu}\phi_{\imath\mu})(k))\: dk\\
&&
+o(1)\ .
\end{eqnarray*}
The remainder term $o(1)$ contains integrals over $K$ of terms involving the
normal derivative of $\chi$, the difference $\chi(ka_R)-\chi_\infty(k)$
and the function $a_R^{-\rho-\epsilon}R(\ii\mu,.,ka_R)$ appearing in (\ref{diable}).  
We combine this remainder with the term (\ref{ad3}) divided by $\omega_R$ to $F(\lambda,\mu,R)$. Since the asymptotic expansion 
(\ref{diable}) can be  differentiated with respect to $\mu$,  
there exists a constant $C\in\R$ such that
\begin{equation}\label{ceins}
|F(\lambda,\mu,R)|+|\frac{d}{d\mu}F(\lambda,\mu,R)|<C , \quad \forall R>0, \: \mu\in\aaaa^*_+\ .\end{equation}  
We can write for $\mu\not=\lambda$
\begin{eqnarray}\label{kunstwerk}
\frac{1}{\omega_R}S_R(\mu) 
&=&
\ii \frac{a_R^{\ii(\lambda-\mu)}}{-\lambda+\mu}\:
(T_2^*c_\gamma(\ii\mu)^*c_\gamma(\ii\lambda)T_1\psi,\phi_{\imath\mu})_B
+\ii \frac{a_R^{\ii(\lambda+\mu)}}{-\lambda-\mu}\:
(T_2^{w*}c_\gamma(\ii\lambda)T_1 \psi,\hat S_{\ii\mu}\phi_{\imath\mu})_B\nonumber\\
&&
+\ii \frac{a_R^{\ii(-\lambda-\mu)}}{\lambda+\mu}\:
((\hat S_{\ii\lambda}\psi,T_1^{w*}c_\gamma(\ii\mu)T_2\phi_{\imath\mu}))_B
+\ii \frac{a_R^{\ii(-\lambda+\mu)}}{\lambda-\mu}\:
(T_2^{w*}T_1^w \hat S_{\ii\lambda} \psi,\hat S_{\ii\mu}\phi_{\imath\mu})_B \quad\quad\\
&&
+\frac{F(\lambda,\mu,R)}{\lambda^2-\mu^2}\ . \nonumber
\end{eqnarray}
By the Lemma of Riemann-Lebesgue 
\begin{eqnarray*}
\lim_{R\to\infty}\int_0^\infty \frac{a_R^{\ii(\lambda+\mu)}}{-\lambda-\mu}\:
(T_2^{w*}c_\gamma(\ii\lambda)T_1 \psi,\hat S_{\ii\mu}\phi_{\imath\mu})_B\:
p_\sigma(\ii\mu)\:d\mu&=&0\ ,\\
\lim_{R\to\infty}\int_0^\infty \frac{a_R^{\ii(-\lambda-\mu)}}{\lambda+\mu}\:
(\hat S_{\ii\lambda}\psi,T_1^{w*}c_\gamma(\ii\mu)T_2\phi_{\imath\mu})_B\:
p_\sigma(\ii\mu)\:d\mu&=&0\ .
\end{eqnarray*}
We set $s:=\mu-\lambda$ and regroup the remaining terms of (\ref{kunstwerk}) to
\begin{eqnarray}
&&
\frac{a_R^{\imath s}-a_R^{-\imath s}}{\imath s}
(T_2^*c_\gamma(\ii\mu)^*c_\gamma(\ii\lambda)T_1\psi,\phi_{\imath\mu})_B
\label{lehmann}\\
&-&
a_R^{\imath s }  \frac{1}{\imath s}\left(
(T_2^*c_\gamma(\ii\mu)^*c_\gamma(\ii\lambda)T_1\psi,\phi_{\imath\mu})_B-
(T_1,T_2)(\hat S_{\ii\lambda} \psi,\hat S_{\ii\mu}\phi_{\imath\mu})_B \right)
\label{schulze}\\
&+&
\frac{F(\lambda,\mu,R)}{\lambda^2-\mu^2}\ . \nonumber
\end{eqnarray}
Note that (\ref{lehmann}) is smooth at $\mu=\lambda$.
We claim that (\ref{schulze}) is smooth at $\mu=\lambda$, too.

By $\hat{S}_{\imath\lambda}^*=\hat{S}_{-\imath\lambda}$ and the functional equation of the scattering matrix (\ref{forme}) we obtain
$$\hat S_{\ii\lambda}^*\hat S_{\ii\lambda}=\frac{1}{p_\sigma(\ii\lambda)} \id\ .$$
The claim now follows from   
\begin{equation}\label{meier} T_2^*c_\gamma(\ii\lambda)^*c_\gamma(\ii\lambda)T_1=
\frac{(T_1,T_2)}{p_\sigma(\ii\lambda)}
\end{equation}
which is a consequence of (\ref{cgcgw}) and (\ref{cgcgw1}).
 
Now (\ref{kunstwerk}) forces also $\frac{F(\lambda,\mu,R)}{\lambda^2-\mu^2}$
to be smooth at $\lambda=\mu$.  By (\ref{ceins}) and  
Lebesgue's theorem about dominated convergence we obtain
$$\lim_{R\to\infty} \int_0^\infty \frac{F(\lambda,\mu,R)}{\lambda^2-\mu^2}\: p_\sigma(\imath\mu) \: d\mu  = 0\ .$$
If we integrate (\ref{schulze}) with respect to $s$ and perform
the limit $R\to\infty$, then the result vanishes by the Riemann-Lebesgue lemma.

We now use (\ref{meier}) and the identity of distributions  
$\lim_{r\to\infty}\frac{\sin(rs)}{s }=\pi\delta_0(s)$ in order to compute
\begin{eqnarray*}
\lim_{R\to\infty} \int_0^\infty S_R(\mu) p_\sigma(\ii\mu)\: d\mu
&=& \lim_{R\to\infty} \omega_R\int_{-\infty}^\infty
\frac{a_R^{\imath s}-a_R^{-\imath s}}{\imath s}
(T_2^*c_\gamma(\ii\mu)^*c_\gamma(\ii\lambda)T_1\psi,\phi_{\imath\mu})_B\:
p_\sigma(\ii\mu)\: ds\\
&=&2\pi\: \omega_X \:(T_1,T_2)\: (\psi,\phi_{\imath\lambda})_B\ .
\end{eqnarray*}
This proves the proposition. 
\hB

\section{The Plancherel theorem and spectral decompositions}\label{didi}

In this final section we obtain our explicit Plancherel theorem, i.e., the decomposition
of $L^2(\Gamma\backslash G,\vp)$.
We use the scalar product formula of Proposition \ref{scalar} in order to show
that the subspace of $L^2(\Gamma\backslash G,\vp)$ spanned by the wave packets of Eisenstein series is the absolute continuous subspace, that 
its complement is the discrete subspace, and that there is no singular continuous subspace.
It is not surprising that the absolute continuous part of the Plancherel measure $\kappa$ coincides with the absolute continuous part for $L^2(G)$. 

As a consequence of the decomposition of the Plancherel theorem we derive the spectral decomposition of $L^2(Y,V_Y(\gamma,\vp))$ with respect to the invariant
differential operators.

We first introduce and describe certain subspaces of $\int_{\hat G}^\oplus M_\pi\hat\otimes V_\pi \: d\kappa(\pi)$ corresponding to the partition
$\hat G=\hat G_d\cup\hat G_u\cup\hat G_c$ (see the beginning of Section \ref{relsec}).  

For each $\pi\in\hat{G}$ we fix the scalar product on ${}^\Gamma(V_{\pi^\prime,-\infty}\otimes V_\vp)_d$ such that the 
matrix coefficient map $c_\pi$  (see Corollary \ref{abel}) is unitary.

 We define the Hilbert space associated to discrete series $\hat{G}_d$ by 
$$\cH_{cusp}(\vp) :=\bigoplus_{\pi\in \hat G_d}^{\mbox{\tiny Hilbert}} {}^\Gamma(V_{\pi^\prime,-\infty}\otimes V_\vp)_d\hat\otimes V_\pi\ .$$

According to Proposition \ref{nonte} and Corollary \ref{otto1} for
$\pi=\bar\pi^{\sigma,\lambda}\in\hat G_{c}$ we have an orthogonal  decomposition $${}^\Gamma(V_{\pi^\prime,-\infty}\otimes V_\vp)_d=
C^{-\infty}(\Lambda,V(\tilde\sigma_\lambda))= E_\Lambda(\tilde\sigma_\lambda,\varphi)\oplus U_\Lambda(\tilde\sigma_\lambda,\varphi)\ .$$
The same corollary gives an alternative expression for the restriction of the scalar product to $E_\Lambda(\sigma_\lambda,\vp)$ in terms of the boundary geometry. 
This space is  non-tivial iff $\lambda>0$ belongs to the finite index set
$PS(\tilde\sigma,\vp)= PS_\res(\tilde\sigma,\vp)\cup PS_U(\tilde\sigma,\vp)$ introduced in  Definition \ref{kabel}. 
We define the Hilbert spaces  
\begin{eqnarray*}
\cH_\res(\vp)&:=&
\bigoplus_{\{\sigma\in\hat M\:|\: p_\sigma(0)=0\}}^{\mbox{\tiny Hilbert}}
\bigoplus_{\lambda\in PS_\res(\tilde\sigma,\vp)\setminus\{0\}} E_\Lambda(\tilde\sigma_\lambda,\varphi)\otimes I^{\sigma,\lambda}
\\
{}_0\cH_{U}(\vp)&:=&
\bigoplus_{\{\sigma\in\hat M\:|\: p_\sigma(0)=0\}}^{\mbox{\tiny Hilbert}}
\bigoplus_{\lambda\in PS_U(\tilde\sigma,\vp)\setminus\{0\}} U_\Lambda(\tilde\sigma_\lambda,\varphi)\otimes I^{\sigma,\lambda}
\ .
\end{eqnarray*}
The sum $\cH_\res(\vp)\oplus {}_0\cH_{U}(\vp)$ is the Hilbert space
associated to the complementary series $\hat{G}_c$.

Now we consider the Hilbert space asssociated to the unitary principal series.
First we discuss the contribution of $\pi^{\sigma,0}$.
We decompose $\hat{M}=\cup_{i=1}^3 \hat{M}_i$ such that
\begin{enumerate}
\item $\sigma\in\hat{M}_1$ iff it is Weyl-invariant and $p_\sigma(0)=0$ (i.e. $\pi^{\sigma,0}$ irreducible), 
\item $\sigma\in\hat{M}_2$ iff it is Weyl-invariant and $p_\sigma(0)\not=0$ (i.e. $\pi^{\sigma,0}=\pi^{\sigma,+}\oplus\pi^{\sigma,-}$), and  
\item $\sigma\in\hat{M}_3$ iff is not Weyl-invariant.
\end{enumerate}
We define
$$ {}_1\cH_{U}(\vp):=\bigoplus_{\{\sigma\in\hat{M}_1 \:|\: 0\in PS_U(\tilde{\sigma},\vp)\}}^{\mbox{\tiny Hilbert}}
U_\Lambda(\tilde\sigma_0,\varphi)\otimes L^2(\partial X,V(\sigma_0))
\ .$$
If $\sigma\in\hat{M}_2$, then 
$L^2(\partial X,V(\sigma_0))=V_{\pi^{\sigma,+}}\oplus V_{\pi^{\sigma,+}}$,
and we define
\begin{eqnarray*}
U^\pm_\Lambda(\tilde{\sigma}_0,\vp)&:=&\{f\in U_\Lambda(\tilde{\sigma}_0,\vp)\:|\:  \langle f,g\rangle=0\ \forall g\in V_{\pi^{\sigma,\mp},\infty}\} \\ 
&=&\{f\in C^{-\infty}(\Lambda,V(\tilde{\sigma}_0,\vp))\:|\: J_0(f)=\pm f\}
\ .\end{eqnarray*}
We define
$$ {}_2\cH_{U}(\vp):=\bigoplus_{\{\sigma\in\hat{M}_2 \:|\: 0\in PS_U(\tilde{\sigma},\vp)\}}^{\mbox{\tiny Hilbert}}
(U^+_\Lambda(\tilde\sigma_0,\varphi)\otimes V_{\pi^{\sigma,+}})\oplus 
(U^-_\Lambda(\tilde\sigma_0,\varphi)\otimes V_{\pi^{\sigma,-}})\ .$$
Let $\hat{M}_4\subset \hat{M_3}$ be a set of representatives of $\hat{M}_3/W(\gaaa,\aaaa)$.
We define 
$${}_3\cH_{U}(\vp):=\bigoplus_{\{\sigma\in\hat{M}_4 \:|\: 0\in PS_U(\tilde{\sigma},\vp)\}}^{\mbox{\tiny Hilbert}}
U_\Lambda(\tilde\sigma_0,\varphi)\otimes L^2(\partial X,V(\sigma_0))
\ .$$
It seems to be natural to collect together all spaces connected with $U_\Lambda(\sigma_\lambda,\vp)$, $\lambda\in PS_U(\sigma,\vp)$, and to define
$$\cH_U(\vp):=\bigoplus_{i=0}^3\: {}_i \cH_U(\vp)\ .$$

The main contribution of unitary principal series is the Hilbert space
$$\cH_{ac}(\vp):=
\bigoplus_{\sigma\in\hat M}^{\mbox{\tiny Hilbert}} \cH_{ac}^\sigma(\vp)\ ,$$
where 
$$\cH_{ac}^\sigma(\vp):=\int_{\aaaa^*_+}^\oplus L^2(B,V_B(\tilde\sigma_{\ii\lambda},\vp))
\hat\otimes L^2(\partial X,V(\sigma_{-\ii\lambda}))\: \frac{2\pi\omega_X}{\dim(V_\sigma)}  p_\sigma(\ii\lambda)
 \: d\lambda\ .$$
For $\sigma\in\hat{M}$ we define the  Wave packet transform  
$$WP_\sigma:\cH_0^{\tilde{\sigma}}(\vp)\otimes L^2(\partial X,V(\sigma_{-\ii\lambda}))_K\rightarrow L^2(\Gamma\backslash G,\vp)$$ 
in the following way (the space  $\cH^{\tilde{\sigma}}_0(\vp)$ is defined just before Definition \ref{now}). 
Consider $T\in \Hom_M(V_{\tilde\sigma},V_{\tilde\gamma})$ for some $\gamma\in\hat K$, $v\in V_\gamma$, and $\phi\in \cH^{\tilde{\sigma}}_0(\vp)$.
We take the element $v_T\in 
L^2(\partial X,V(\sigma_{-\ii\lambda}))_K$ (see Lemma \ref{circe}) and form
$\phi\otimes v_T$. Then we define 
$$ WP_\sigma(\phi \otimes v_T)=\langle E(\phi,T),v\rangle 
\ .$$
We  employ the extension $ext$ in order to identify the space $L^2(B,V_B(\tilde\sigma_{\ii\lambda},\vp))$ with a subspace of $ {}^\Gamma(C^{-\infty}(\partial X,V(\tilde{\sigma}_{\imath\lambda}))\otimes V_\vp)_{temp}$ which is our candidate of $N_{\pi^{\sigma,-\imath\lambda}}$. 
Then by Lemma \ref{circe} the wave packet transform $WP_\sigma$ is related to the family of matrix coefficient maps $\{c_{\pi^{\sigma,-\imath\lambda}}\}_{\lambda\in\aaaa^*_+}$ by
\begin{equation}\label{willdoch}WP_\sigma(\phi\otimes v_T)=\frac{2\pi\omega_X}{\dim(V_\sigma)} \int_{\aaaa^*_+} c_{\pi^{\sigma,-\imath\lambda}}(ext(\phi)\otimes v_T)  
 p_\sigma(\ii\lambda)
\: d\lambda\ .\end{equation}
Note that the elements
of the form $\phi\otimes v_T$ span the dense subspace $\cH_0^{\tilde{\sigma}}(\vp)\otimes L^2(\partial X,V(\sigma_{-\ii\lambda}))_K$ of $\cH_{ac}^\sigma(\vp)$.

We are now able to state the Plancherel theorem for $L^2(\Gamma\backslash G,\vp)$. Recall that in case $X=\OO H^2$ we assume $\delta_\Gamma<0$.

\begin{theorem}\label{pl2} 
The direct sum of the matrix coefficient maps 
$c_\pi:{}^\Gamma (V_{\pi^\prime,-\infty}\otimes V_\vp)_d \otimes V_{\pi,\infty}\rightarrow L^2(\Gamma\backslash G,\vp)$, $\pi\in\hat{G}$, and the
wave packet transforms $WP_\sigma$, $\sigma\in \hat{M}$, extends to a unitary equivalence of $G$-representations
$$ \cH_{ac}(\vp)\oplus \cH_{cusp}(\vp)\oplus 
\cH_\res(\vp) \oplus \cH_{U}(\vp) \cong L^2(\Gamma\backslash G,\vp)\ .$$
It gives rise to a corresponding decomposition
$$ L^2(\Gamma\backslash G,\vp) = L^2(\Gamma\backslash G,\vp)_{ac}\oplus L^2(\Gamma\backslash G,\vp)_d
\ ,$$
where the discrete subspace
$$L^2(\Gamma\backslash G,\vp)_d:=L^2(\Gamma\backslash G,\vp)_{cusp}\oplus L^2(\Gamma\backslash G,\vp)_{\res}\oplus L^2(\Gamma\backslash G,\vp)_{U}
$$
is the sum of the cuspidal, the residual, and the "stable" part.

$L^2(\Gamma\backslash G,\vp)_d$ is the sum of all irreducible subrepresentations of $L^2(\Gamma\backslash G,\vp)$. $L^2(\Gamma\backslash G,\vp)_{cusp}$ decomposes into discrete series representations of $G$, each
discrete series representation of $G$ occurs with infinite multiplicity. It is empty
iff $X=\R H^n$, $n$ odd. The remaining part of $L^2(\Gamma\backslash G,\vp)_d$
consists of non-discrete series representations of $G$ with real infinitesimal
character occuring with finite multiplicity. If $\delta_\Gamma<0$, then it is
empty. If $\delta_\Gamma>0$ and $\vp=1$, then it contains the representation $I^{1,\delta_\Gamma}$ with multiplicity one.

$L^2(\Gamma\backslash G,\vp)_{ac}$ decomposes into a sum of direct integrals
corresponding to the unitary principal series representations of $G$, each
occuring with infinite multiplicity.
\end{theorem}
The notions $L^2(\Gamma\backslash G,\vp)_{cusp}$ and $L^2(\Gamma\backslash G,\vp)_{\res}$ are chosen in analogy with the case of groups $\Gamma$ with
finite covolume. Indeed, $L^2(\Gamma\backslash G,\vp)_{\res}$ is spanned
by the residues of Eisenstein series. However, the "cusp forms" forming
the space $L^2(\Gamma\backslash G,\vp)_{cusp}$ share the properties
of the cusp forms in the sense of Harish-Chandra associated to the 
trivial group and
not of those for groups with finite covolume. The appearance of $L^2(\Gamma\backslash G,\vp)_{U}$ does not seem to have an analogue. 

\noindent
{\it Proof of the theorem.$\:\:\:\:$} 
It follows from Corollary \ref{abel}, the determination of ${}^\Gamma(V_{\pi^\prime,-\infty}\otimes V_\vp)_d$ in Section \ref{relsec}, and
our definition of the scalar products on ${}^\Gamma(V_{\pi^\prime,-\infty}\otimes V_\vp)_d$
that the matrix coefficient maps induce a $G$-equivariant unitary map of $\cH_{cusp}(\vp)\oplus \cH_{\res}(\vp)\oplus \cH_{U}(\vp)$ onto  the discrete
subspace
$L^2(\Gamma\backslash G,\vp)_d$.
Next we show that the wave packet transform $WP_\sigma$
extends to a unitary embedding of $\cH_{ac}^\sigma(\vp)$
into $L^2(\Gamma\backslash G,\vp)$.
We compute using Proposition \ref{scalar} and (\ref{u8u8u})

\begin{eqnarray*}
 \lefteqn{\|WP_\sigma(\phi\otimes v_T)\|_{L^2(\Gamma\backslash G,\vp)}^2}\hspace{0.5cm}\\
&=&
\big( \int_{\aaaa^*_+} \langle E(\imath\lambda,\phi_{\imath\lambda},T),v\rangle  \:p_\sigma(\imath\lambda) d\lambda\:,\: \langle E(\phi,T), v\rangle \big)_{L^2(\Gamma\backslash G,\vp)}\\
&=&
\int_{\aaaa^*_+}\langle  ( E(\imath\lambda,\phi_{\imath\lambda},T),v\rangle\: ,\: \langle E(\phi,T),v\rangle )_{L^2(\Gamma\backslash G,\vp)} \: p_\sigma(\imath\lambda) d\lambda\\
&=&
\int_{\aaaa^*_+} \int_{\Gamma\backslash G} \int_K \langle E(\imath\lambda,\phi_{\imath\lambda},T)(gk),v\rangle   \overline{\langle  E(\phi,T)(gk),v   \rangle }\: dk \:dg  \:p_\sigma(\imath\lambda) d\lambda\\
&=&
 \int_{\aaaa^*_+} \int_{\Gamma\backslash G} \int_K \langle \tilde\gamma(k^{-1})E(\imath\lambda,\phi_{\imath\lambda},T)(g),v\rangle   \overline{\langle \tilde\gamma(k^{-1}) E(\phi,T)(g), v \rangle}  \:dk\: dg \: p_\sigma(\imath\lambda) d\lambda\\
&=&
\frac{\|v\|^2}{\dim(V_\gamma)} \int_{\aaaa^*_+} \int_{\Gamma\backslash G}   ( E(\imath\lambda,\phi_{\imath\lambda},T)(g) , E(\phi,T)(g) )  \:dg \:  p_\sigma(\imath\lambda) d\lambda\\
 &=&\frac{2\pi\omega_X \|T\|^2 \|v\|^2}{\dim(V_\gamma)}
\int_{\aaaa^*_+}    ( \phi_{\imath\lambda},\phi_{\imath\lambda} )_B  \:  p_\sigma(\imath\lambda) d\lambda\\
&=&\frac{2\pi\omega_X   \|v_T\|^2}{\dim(V_\sigma)}\int_{\aaaa^*_+}    ( \phi_{\imath\lambda},\phi_{\imath\lambda} )_B   \: p_\sigma(\imath\lambda) d\lambda\\
&=&   \|\phi\otimes v_T\|_{\cH^\sigma_{ac}(\vp)}^2 \ .      
\end{eqnarray*}
By a similar computation we see that the images of the $WP_\sigma$, $\sigma\in\hat{M}$, are pairwise orthogonal.
It follows from (\ref{willdoch}) that the extension of $WP_\sigma$ to $\cH^\sigma_{ac}(\vp)$ is $G$-equivariant.  
As a subspace of the continuous subspace of $L^2(\Gamma\backslash G,\vp)$
the image of $WP_\sigma$ is also  orthogonal to
$L^2(\Gamma\backslash G,\vp)_d$.

It remains to show that the image of $\cH_{ac}(\vp)\oplus\cH_{cusp}(\vp)\oplus \cH_{\res}(\vp)\oplus \cH_{U}(\vp)$ exhausts $L^2(\Gamma\backslash G,\vp)$. Let $f\in L^2(\Gamma\backslash G,\vp)$ be orthogonal to that
image. We may and shall assume that it belongs to a $K$-isotypic component 
$L^2(\Gamma\backslash G,\vp)(\gamma)$ for 
some $\gamma\in\hat K$. By Corollary \ref{abel} we can compute the
scalar product of $f$ with a Schwartz function $g\in \cC(\Gamma\backslash G,\vp)$ in the following way
\begin{eqnarray}\label{zorn}
(f,g)&=& \int_{\hat G}(\cF(f)(\pi),\cF (g)(\pi))_{N_\pi\otimes V_\pi}\:d\kappa(\pi)\nonumber\\
&=&\int_{\hat G}(\cF(f)(\pi),\cF_\pi(g))_{N_\pi\otimes V_\pi}\:d\kappa(\pi)\nonumber\\
&=& \int_{\hat G}(c_\pi(\cF(f)(\pi)),g)_{L^2(\Gamma\backslash G,\vp)}\:d\kappa(\pi)\ .
\end{eqnarray}
Since $f$ is orthogonal to $L^2(\Gamma\backslash G,\vp)_d$ we have
$\cF(f)(\pi)=0$ for all $\pi\in\hat G$ with $\kappa(\pi)\not=0$, 
i.e., ${}^\Gamma(V_{\pi^\prime,-\infty}\otimes V_\vp)_d\not=0$.
Thus, by the
results of Section \ref{relsec},
it remains
to show that $\cF(f)(\pi)=0$ for the unitary principal series representations
$\pi^{\sigma,-\ii\lambda}$. Because of the equivalence $\pi^{\sigma,\ii\lambda}\cong \pi^{\sigma^w,-\ii\lambda}$ we can assume $\lambda\in\aaaa^*_+$.
Clearly, $\cF(f)(\pi^{\sigma,-\ii\lambda})=0$ if
$[\gamma_{|M}:\sigma]=0$.

 By Corollary \ref{abel} and Lemma \ref{zeus} we have $$\cF(f)(\pi^{\sigma,-\ii\lambda})\in {}^\Gamma C^{-\infty}(\partial X,V(\tilde\sigma_{\ii\lambda},\vp))\otimes
L^2(\partial X,V(\sigma_{-\ii\lambda}))(\gamma)\ .$$ For each $\sigma\in \hat M$, with $[\gamma_{|M}:\sigma]\not=0$ let $\{T^\sigma_i\}_{i=1}^{\dim\:\Hom_M(V_{\tilde\sigma},V_{\tilde\gamma})}$ be a basis of $\Hom_M(V_{\tilde\sigma},V_{\tilde\gamma})$. Let $\{v^j\}_{j=1}^{\dim(V_\gamma)}$ be
a basis of $V_\gamma$. We conclude that there are sections
$\phi^{ij}_{\sigma,\ii\lambda}\in  C^{-\infty}(B,V_B(\tilde\sigma_{\ii\lambda},\vp))$, $i=1,\dots ,\dim\:\Hom_M(V_{\tilde\sigma},V_{\tilde\gamma})$, $j=1,\dots, \dim(V_\gamma)$,
such that
$$\cF(f)(\pi^{\sigma,-\ii\lambda})=\sum_{i,j} ext(\phi^{ij}_{\sigma,\ii\lambda})
\otimes v^j_{T_i^\sigma} \ .$$
Using Lemma \ref{circe} we find
$$c_{\pi^{\sigma,-\ii\lambda}}(\cF(f)(\pi^{\sigma,-\ii\lambda}))=\sum_{i,j} \langle E(\ii\lambda,\phi^{ij}_{\sigma,\ii\lambda},T^\sigma_i), v^j \rangle\ .
$$ 
We now evaluate the scalar product of $f$ with some wave packet of the
form $WP_{\sigma^\prime}(\psi\otimes v_T)$, where $\psi$ belongs to the space $\cH_0^{\tilde\sigma^\prime}(\vp)$, $v\in V_{\gamma}$, and $T\in\Hom(V_{\tilde{\sigma}^\prime},V_{\tilde{\gamma}})$. 

The map $\aaaa^*_+\ni\lambda\mapsto \pi^{\sigma,-\ii\lambda}\in
\hat G$ identifies $\aaaa^*_+$ with an open subset of $\hat{G}$.
Let $d\kappa(\sigma,\lambda)$ be the restriction of the Plancherel measure 
$d\kappa$ to this subset. 
Using (\ref{zorn}), Proposition \ref{scalar} we obtain 
 \begin{eqnarray*}
0&=& (f, WP_{\sigma^\prime}(\psi\otimes v_T))_{L^2(\Gamma\backslash G,\vp)}\\
&=&\int_{\hat G} (c_\pi(\cF(f)(\pi)), WP_{\sigma^\prime}(\psi\otimes v_T))\:d\kappa(\pi)\\
&=&\sum_{\sigma\in \hat{M},[\gamma_{|M}:\sigma]\not=0 }\sum_{i,j} \int_{\aaaa^*_+} ( \langle E(\ii\lambda,\phi^{ij}_{\sigma,\ii\lambda},T^\sigma_i), v^j \rangle, \langle E(\psi,T),v\rangle)\:d\kappa(\sigma,\lambda)\\
&=&2\pi\omega_X\sum_{i,j}(v_j,v)(T^{\sigma^\prime}_i,T) \int_{\aaaa^*_+} 
(\phi^{ij}_{\sigma^\prime,\ii\lambda},\psi_{\ii\lambda})_B \:d\kappa(\sigma^\prime,\lambda)\ .
\end{eqnarray*}
Varying $\sigma^\prime$, $T$ and $v$  we see that for all $i,j$ and
$\sigma^\prime$
$$\int_{\aaaa^*_+} 
(\phi^{ij}_{\sigma^\prime,\ii\lambda},\psi_{\ii\lambda})_B \:d\kappa(\sigma^\prime,\lambda)
=0\ .$$
Moreover, we are free to multiply $\psi$ with a function 
$h\in C_c^{\infty}(\aaaa^*_+)$. Hence for all such $h$
$$\int_{\aaaa^*_+} h(\lambda)
(\phi^{ij}_{\sigma^\prime,\ii\lambda},\psi_{\ii\lambda})_B \:d\kappa(\sigma^l,\lambda)
=0\ .$$ 
We conclude that $(\phi^{ij}_{\sigma^\prime,\ii\lambda},\psi_{\ii\lambda})_B
=0$ for almost all $\lambda$ (mod $\kappa$).
We choose a countable dense set $\{\psi_m\}\subset C^\infty(B,V_B(\tilde\sigma^\prime_{0},\vp))$. Using a holomorphic trivialization of the family of bundles $\{V_B(\tilde\sigma^\prime_\mu,\varphi)\}_{\mu\in\aca}$
we extend these sections to families $\aca\ni\mu\mapsto \psi_{m,\mu}\in C^\infty(B,V_B(\tilde\sigma^\prime_{\mu},\vp))$.
We form
$B_m:=\{\lambda|(\phi^{ij}_{\sigma^\prime,\ii\lambda},\psi_{m,\imath\lambda})\not=0\}\subset \aaaa^*_+$.
Then $\kappa(B_m)=0$. Moreover let $U:=\bigcup_m B_m$.
Then $\kappa(U)=0$, and we have 
$(\phi^{ij}_{\sigma^\prime,\ii\lambda},\psi_{m,\imath\lambda})=0$ for all $\lambda\in \aaaa_+^*\setminus U$
and all $m$.  Thus $\phi^{ij}_{\sigma^\prime,\ii\lambda}=0$ for $\lambda\in \aaaa_+^*\setminus U$.
Hence $\cF(f)(\pi)=0$ for almost all $\pi$ (mod $\kappa$).
Therefore  $f=0$. This proves that  the image of $\cH_{ac}(\vp)\oplus\cH_{cusp}(\vp)\oplus \cH_{\res}(\vp)\oplus \cH_{U}(\vp)$ is all of $L^2(\Gamma\backslash G,\vp)$ and hence the
theorem.
\hB

Let now $\gamma$ be a finite-dimensional unitary representation of $K$.
We want to draw consequences of Theorem \ref{pl2} for the spectral decomposition of the algebra $\cZ$ containing the Casimir operator $\Omega_G$ and acting by unbounded operators on $L^2(Y,V_Y(\gamma,\vp))$. 
Let $D(G,\gamma)$ be the algebra of $G$-invariant differential operators on $V(\gamma)$. This algebra 
is a finitely generated module over $\cZ$, where according to our convention
the homomorphism of $\cZ$ to $D(G,\gamma)$ is induced by the left regular
representation of $G$.
The algebra $D(G,\gamma)$  might
be noncommutative, even if $\gamma$ is irreducible (but then $X=\HH H^n$ or 
$\OO H^2$). 

As before, we identify  $C^\infty(Y,V_Y(\gamma,\vp))$ with $[C^\infty(\Gamma\backslash G,\vp)\otimes V_\gamma]^K$. This identification provides an isomorphism of $[\cU(\gaaa)\otimes_{\cU(\kaaa)} \End(V_\gamma)]^K$ with  $D(G,\gamma)$ , where the action of $\cU(\gaaa)$
on $C^\infty(\Gamma\backslash G,\vp)$ is induced by the right regular representation. In particular,  
the decomposition of $L^2(\Gamma\backslash G,\vp)$
with respect to $G$ induces a decomposition of $L^2(Y,V_Y(\gamma,\vp))$ with
respect to $D(G,\gamma)$.

If $\pi$ is any admissible representation of $G$, 
then the finite-dimensional space $\Hom_K(V_\pi,V_\gamma)$ has a natural structure of a $D(G,\gamma)$-module. In fact, the action of $\cU(\gaaa)\otimes \End(V_\gamma)$ on $\Hom(V_\pi,V_\gamma)$ induces an action
of $D(G,\gamma)\cong [\cU(\gaaa)\otimes_{\cU(\kaaa)} \End(V_\gamma)]^K$ on $\Hom_K(V_\pi,V_\gamma)$. Note that if $\pi$ is irreducible, then  $\Hom_K(V_\pi,V_\gamma)$ is an irreducible
$D(G,\gamma)$-module.  
Moreover, if $\gamma$ is irreducible, then  the functor 
$\pi\mapsto \Hom_K(V_\pi,V_\gamma)$ provides a one-to-one correspondence
between equivalence classes of irreducible $G$-representations (strictly speaking of irreducible $(\gaaa,K)$-modules) containing $\gamma$
and of irreducible $D(G,\gamma)$-modules (see \cite{wallach88}, 3.5.4.). 
The induced action of $\cZ$  on $\Hom_K(V_\pi,V_\gamma)$ is given by the infinitesimal character of $\pi$.

We define the Hilbert space and $D(G,\gamma)$-module
$$\cH^\gamma_{cusp}(\vp) :=\bigoplus_{\pi\in \hat G_d} {}^\Gamma(V_{\pi,-\infty}\otimes V_\vp)_d \otimes \Hom_K(V_\pi,V_\gamma)\ .$$
This sum is finite since every irreducible subrepresentation of $\gamma$ only occurs in finitely many discrete series representations. There is a 
natural unitary $D(G,\gamma)$-equivariant map of $\cH^\gamma_{cusp}(\vp)$ into $L^2(\Gamma\backslash G,\vp)$ which is given on each summand by
\begin{equation}\label{mamf} 
{}^\Gamma(V_{\pi,-\infty}\otimes V_\vp)_d \otimes \Hom_K(V_\pi,V_\gamma)
\ni \phi\otimes t\mapsto t(\pi(g^{-1})\otimes\id)\phi)\in [L^2(\Gamma\backslash G,\vp)\otimes V_\gamma]^K\ .
\end{equation}

By the Frobenius reciprocity $\Hom_M(V_\sigma,V_\gamma)\cong \Hom_K(V_{\pi^{\sigma,\lambda}},V_\gamma)$ the space
$\Hom_M(V_\sigma,V_\gamma)$ carries a structure of a 
$D(G,\gamma)$-module depending on $\lambda\in\aca$.
We denote this module by $\Hom_M(V_\sigma,V_\gamma)_\lambda$.
For $\Ree(\lambda)>0$ we consider $I^{\sigma,\lambda}$ as a submodule of 
$V_{\pi^{\sigma,\lambda}}$. We define the $D(G,\gamma)$-module $\overline{\Hom_M(V_\sigma,V_\gamma)}_\lambda$ to be the quotient 
of $\Hom_M(V_\sigma,V_\gamma)_\lambda$ which corresponds to the quotient
\linebreak[4]
$\Hom_K(I^{\sigma,\lambda},V_\gamma)$ of
$\Hom_K(V_{\pi^{\sigma,\lambda}},V_\gamma)$ by Frobenius reciprocity.
Using Lemma \ref{149u}, 1., one can check that
$\overline{\Hom_M(V_\sigma,V_\gamma)}_\lambda$ is the quotient of
$\Hom_M(V_\sigma,V_\gamma)_\lambda$ by the subspace which is annihilated by
multiplication by $c_\gamma(\lambda)$. 

Recall the definition of the scalar product $(T_1,T_2)\id_{V_{\sigma}}:=T_2^*T_1$ on
$\Hom_M(V_\sigma,V_\gamma)$. For
$\lambda\in \imath\aaaa^*$ this scalar product is equal to
the scalar product induced by Frobenius reciprocity $\Hom_M(V_\sigma,V_\gamma)\cong \Hom_K(L^2(\partial X,V(\sigma,\lambda)),V_\gamma)$ rescaled by
$\dim(V_\sigma)^{-1}$. Note that if $\sigma\in\hat{M}_2$, then
there is an orthogonal decomposition 
$$\Hom_M(V_\sigma,V_\gamma)_0=\Hom_M(V_\sigma,V_\gamma)_0^+\oplus\Hom_M(V_\sigma,V_\gamma)_0^-$$
corresponding to the decomposition $\pi^{\sigma,0}=\pi^{\sigma,+}\oplus\pi^{\sigma,-}$.

If $\Ree(\lambda)>0$ and $\bar{\pi}^{\sigma,\lambda}$ is unitary, then
we define the scalar product
on $\overline{\Hom_M(V_\sigma,V_\gamma)}_\lambda$ by
$$([T_1],[T_2]):=\frac{\dim(V_\sigma)}{c_\sigma(\lambda)}(c_\gamma(\lambda)T_1,T_2^w) \ ,$$ where $T_i\in \Hom_M(V_\sigma,V_\gamma)$ are representatives
of $[T_i]\in \overline{\Hom_M(V_\sigma,V_\gamma)}_\lambda$.
Note that this scalar product coincides with the scalar product 
induced by $\overline{\Hom_M(V_\sigma,V_\gamma)}_\lambda\cong \Hom_K(I^{\sigma,\lambda},V_\gamma)$.

We define the Hilbert spaces and $D(G,\gamma)$-modules

\begin{eqnarray*}
\cH^\gamma_\res(\vp)&:=&
\bigoplus_{\{\sigma\in\hat{M}\:|\:[\gamma_{|M}:\sigma]\not=0\:,\: p_\sigma(0)=0\}}
\bigoplus_{\lambda\in PS_\res(\sigma,\vp)\setminus\{0\}} E_\Lambda(\sigma_\lambda,\varphi)\otimes \overline{\Hom_M(V_\sigma,V_\gamma)}_\lambda \\
{}_0\cH^\gamma_{U}(\vp)&:=&
\bigoplus_{\{\sigma\in\hat{M}\:|\:[\gamma_{|M}:\sigma]\not=0\:,\: p_\sigma(0)=0\}} 
\bigoplus_{\lambda\in PS_U(\sigma,\vp)\setminus\{0\}} U_\Lambda(\sigma_\lambda,\varphi)\otimes \overline{\Hom_M(V_\sigma,V_\gamma)}_\lambda \\
{}_1\cH^\gamma_{U}(\vp)&:=&\bigoplus_{\{\sigma\in\hat{M}_1\:|\:[\gamma_{|M}:\sigma]\not=0 \:,\: 0\in PS_U({\sigma},\vp)\}} 
U_\Lambda(\sigma_0,\varphi)\otimes \Hom_M(V_\sigma,V_\gamma)_0\\
{}_2\cH^\gamma_{U}(\vp)&:=&\bigoplus_{\{\sigma\in\hat{M}_2\:|\:[\gamma_{|M}:\sigma]\not=0 \:,\: 0\in PS_U({\sigma},\vp)\}} U^+_\Lambda(\sigma_0,\varphi)\otimes \Hom_M(V_\sigma,V_\gamma)_0^+\oplus 
U^-_\Lambda(\sigma_0,\varphi)\otimes \Hom_M(V_\sigma,V_\gamma)_0^- \\ {}_3\cH^\gamma_{U}(\vp)&:=&\bigoplus_{\{\sigma\in\hat{M}_4\:,\:[\gamma_{|M}:\sigma]\not=0 \:|\: 0\in PS_U({\sigma},\vp)\}} 
U_\Lambda(\sigma_0,\varphi)\otimes \Hom_M(V_\sigma,V_\gamma)_0
\ .
\end{eqnarray*}
All these sums are finite.
We further define
$$\cH^\gamma_U(\vp):=\bigoplus_{i=0}^3\: {}_i \cH^\gamma_U(\vp)\ .$$
The matrix coefficient maps defined on
$\cH_{res}(\vp)\oplus\cH_U(\vp) $ induces a unitary map of $D(G,\gamma)$-modules
from $\cH_{res}^\gamma (\vp)\oplus\cH_U^\gamma(\vp)$
to $L^2(Y,V_Y(\gamma,\vp))$ which is given by the Poisson transform
on each summand (see Lemma \ref{circe}). 
In particular, if $ \lambda >0$,
then the Poisson transform factors over 
$(E_\Lambda(\sigma_\lambda,\vp)\oplus U_\Lambda(\sigma_\lambda,\vp)) \otimes
\overline{\Hom_M(V_\sigma,V_\gamma)}_\lambda$.
By (\ref{formelop}) the image under the Poisson transform of $\cH^\gamma_\res(\vp)$ consists exactly of the residues of Eisenstein series.

Last not least we introduce the absolute continuous part as the finite sum
$$\cH^\gamma_{ac}(\vp):=
\bigoplus_{\sigma\in\hat{M}\:,\:[\gamma_{|M}:\sigma]\not=0} \cH_{ac}^{\sigma,\gamma}(\vp)\ ,$$
where $\cH_{ac}^{\sigma,\gamma}$ is the direct integral of  Hilbert spaces
and $D(G,\gamma)$-modules  
  $$\cH_{ac}^{\sigma,\gamma}(\vp):=\int_{\aaaa^*_+}^\oplus L^2(B,V_B(\sigma,\ii\lambda))
\otimes \Hom_M(V_\sigma,V_\gamma)_{\imath\lambda}\: 2\pi\omega_X p_\sigma(\ii\lambda)
 \: d\lambda\ .$$
By Proposition \ref{scalar} the wave packet transforms extend to isometric 
$D(G,\gamma)$-equivariant embeddings
of $\cH_{ac}^\gamma(\vp)$ into $L^2(Y,V_Y(\gamma,\vp))$. Note
that $D(G,\gamma)$ acts on $\cH_{ac}^{\sigma,\gamma}(\vp)$ as an algebra
of unbounded operators via multiplication by $\End(\Hom_M(V_\sigma,V_\gamma))$-valued polynomials
on $\aaaa^*$.

Now the following theorem is an immediate consequence
of Theorem \ref{pl2}.

\begin{theorem}\label{pl3}
The maps (\ref{mamf}), the Poisson transforms and the wave packets of
Eisenstein series combine to a unitary equivalence of $D(G,\gamma)$-representations
$$\cH^\gamma_{ac}(\vp)\oplus \cH^\gamma_{cusp}(\vp)\oplus 
\cH^\gamma_\res(\vp) \oplus \cH^\gamma_{U}(\vp)\cong   L^2(Y,V_Y(\gamma,\vp)) \ .$$
It gives rise to a corresponding decomposition
$$ L^2(Y,V_Y(\gamma,\vp)) = L^2(Y,V_Y(\gamma,\vp))_{ac}\oplus L^2(Y,V_Y(\gamma,\vp))_d
\ ,$$
where the discrete subspace
$$L^2(Y,V_Y(\gamma,\vp))_d:=L^2(Y,V_Y(\gamma,\vp))_{cusp}\oplus L^2(Y,V_Y(\gamma,\vp))_{\res}\oplus L^2(Y,V_Y(\gamma,\vp))_{U}
$$
is this the sum of the cuspidal, the residual, and the "stable" part.

The algebra $\cZ$ has pure point spectrum on $L^2(Y,V_Y(\gamma,\vp))_d$, whereas its spectrum on \linebreak[4]  $L^2(Y,V_Y(\gamma,\vp))_{ac}$ consists of
finitely many branches  of absolute continuous spectrum  of infinite multiplicity. 
$L^2(Y,V_Y(\gamma,\vp))_{cusp}$ is a finite sum of infinite-dimensional
eigenspaces. The remaining part of the discrete subspace is finite-dimensional.
If $\delta_\Gamma<0$, then it is empty. 
\hB
\end{theorem}

We conclude this paper by some comments on Theorem \ref{pl3}.
\begin{itemize}
\item It is clear from Corollary \ref{abel} that 
integration against Eisenstein series gives a map 
$$\cC(Y,V_Y(\gamma,\vp))\rightarrow \cH^\gamma_{ac}(\vp)$$ which we call Eisenstein-Fourier transform. It is a left-inverse of the wave packet transform,
and it would be interesting to investigate its image.
 \item The decomposition of
$L^2(Y,V_Y(\gamma,\vp))$ given by Theorem \ref{pl3} is finer than the spectral decomposition with respect to the Casimir operator. 
In particular, the Casimir operator can have  eigenvalues embedded in the continuous spectrum. 
This is a kind of accident, because these
embedded eigenvalues can be separated from the continuous spectrum by additional
operators belonging to $D(G,\gamma)$
with one possible exception. Namely
an eigenspace arising from $\oplus_{i=1}^3 {}_i\cH^\gamma_U(\vp)$   
contributes an eigenvalue lying at the bottom of one branch of the absolute contiuous spectrum of $\Omega_G$.
\item So far we have not presented any example where the space $L^2(Y,V_Y(\gamma,\vp))_{U}$ is non-trivial. Here is one. It also
sheds some light on the previous remark.   
Let $\Gamma\subset SL(2,\R)$ be a cocompact Fuchsian group.
Consider $\Gamma\subset SL(2,\C)$ in the standard way. The limit set consists
of the equator of the sphere $S^2$. If we interpret the equator as a 1-current,
then it is not difficult to see that it defines an element of $U_\Lambda(\sigma_0,1)$, where $\sigma=\sigma^2\oplus\sigma^{2,w}$ is the  representation of $M\cong U(1)$
corresponding to $1$-forms. The corresponding square integrable $1$-form on
$Y$ is harmonic, i.e., it contributes to the $L^2$-cohomology of $Y$ (compare Mazzeo-Phillips \cite{mazzeophillips90}.
\item Even in the case $\delta_\Gamma>0$ it can happen that $L^2(Y,V_Y(\gamma,\vp))_d=0$. For instance if $Y$ is odd-dimensional, then Corollary \ref{treu} implies that the Dirac operator acting on spinors  
has pure absolute continuous spectrum $(-\infty,\infty)$.
\item Using the meromorphic continuation of the Eisenstein series to all of $\aca$  and Theorem
\ref{pl3} 
one can show that the resolvent kernel of $(\Omega-z)^{-1}$
on $L^2(Y,V_Y(\gamma,\vp))$ extends meromorphically to a finite-sheeted
branched cover of $\C$. 
\end{itemize}

\bibliographystyle{plain}

\end{document}